\documentclass[aos,preprint]{imsart}

\RequirePackage{amsthm,amsmath,amsfonts,amssymb,array}
\RequirePackage[authoryear]{natbib}

\usepackage{xcolor}

\usepackage{hyperref} 
\hypersetup{
    colorlinks=true,
    linkcolor=blue,      
    citecolor=blue,      
    urlcolor=blue,       
    pdfencoding=auto,
    pdfborder={0 0 0},   
    linktoc=all          
}
\usepackage{cleveref} 

\startlocaldefs
\numberwithin{equation}{section}

\theoremstyle{plain}
\newtheorem{theorem}{Theorem}[section]

\newtheorem{lemma}[theorem]{Lemma}
\newtheorem{corollary}[theorem]{Corollary}

\theoremstyle{definition}

\newtheorem{example}[theorem]{Example}

\newcolumntype{Y}[1]{>{\raggedright\arraybackslash}p{#1}}

\endlocaldefs

\begin{document}

\begin{frontmatter}

\title{Sharp Cram\'er-type moderate deviation for Hotelling's $T^2$ statistics}
\runtitle{Sharp moderate deviations for Hotelling's $T^2$-statistic}

\begin{aug}
\author[A]{\fnms{Song-hao}~\snm{Liu}\ead[label=e1]{liusonghao@dlut.edu.cn}},
\author[B]{\fnms{Qi-man}~\snm{Shao}\ead[label=e2]{shaoqm@sustech.edu.cn}},
\author[C]{\fnms{Hao}~\snm{Shi}\ead[label=e3]{Hao.Shi@link.cuhk.edu.hk}}
\and
\author[C]{\fnms{Jing-yu}~\snm{Xu}\ead[label=e4]{12131253@mail.sustech.edu.cn}}


\address[A]{School of Mathematical Sciences, Dalian University of Technology, Dalian, Liaoning, China \printead[presep={,\ }]{e1}}

\address[B]{Department of Statistics and Data Science, Shenzhen International Center for Mathematics, Southern University of Science and Technology, Shenzhen, Guangdong, 518055, China \printead[presep={,\ }]{e2}}

\address[C]{Department of Statistics and Data Science, Southern University of Science and Technology, Shenzhen, Guangdong, 518055, China \printead[presep={,\ }]{e3,e4}}

\end{aug}

\begin{abstract}
Let $T^2_{H,n}$ denote the one-sample Hotelling statistic based on
$n$ independent random vectors in $\mathbb R^d$ with a common mean,
where $d\ge2$ is a fixed integer. Assume that the average covariance
matrix is positive definite. We define the statistic to be zero when the centered sample covariance matrix is singular. Under the null hypothesis, we establish finite-sample Cram\'er-type moderate deviation results and uniform Berry--Esseen bounds under finite $p$th moments for $2<p\le3$. Explicit constructions show that the moderate deviation error order is optimal over independent triangular arrays. For identically distributed observations with a fixed distribution, we obtain
\[
\sup_{0\le x\le b_n}
\left|
\frac{\mathbb P(T^2_{H,n}\ge x^2)}
{\mathbb P(\chi_d^2\ge x^2)}-1
\right|
=
O\!\left((1+b_n)^p n^{-(p-2)/2}\right)
\]
for every nonnegative sequence $b_n=o(n^{(p-2)/(2p)})$.
In this identically distributed setting, the uniform Berry--Esseen
bound is $O(n^{-(p-2)/2})$.
The proof uses a radial Stein identity and a recursive error estimate.
\end{abstract}



\end{frontmatter}

\section{Introduction}\label{sec:intro}

Hotelling's \(T^2\)-statistic is a standard tool for inference on a
multivariate mean. It combines the sample mean and sample covariance
matrix into an affine-invariant statistic and is commonly calibrated
by a chi-squared distribution in fixed dimension.
Its two-sample version is used for inference on the
difference between two multivariate means.
The accuracy of this calibration determines the accuracy of the
resulting critical values, confidence regions, and \(p\)-values.
When many hypotheses are tested, the relevant tail probabilities can be
very small. An approximation that is accurate in absolute terms may then
give a substantial error in the ratio of the true and approximating tail
probabilities. This distinction motivates both Berry--Esseen bounds,
which control the difference between distribution functions, and
Cram\'er-type moderate deviation results, which control the ratio of
tail probabilities; see \citet{FanHallYao2007} and \citet{LiuShao2013}
for their role in simultaneous testing.

The scalar theory provides a natural starting point.
\citet{Shao1999} established a Cram\'er-type moderate deviation result
for Student's statistic under a finite third absolute moment.
\citet{JingShaoWang2003} obtained quantitative results for
self-normalized sums of independent random variables under finite
\((2+\delta)\)th moments, \(0<\delta\le1\). A systematic account is
given by \citet{delaPenaLaiShao2009}. For Hotelling's statistic,
\citet{DemboShShao2006} established large- and moderate-deviation
estimates for the logarithm of the tail probability.
\citet{LiuShao2013} subsequently proved Cram\'er-type moderate deviation
results for both one-sample and two-sample Hotelling statistics under
finite \((3+\delta)\)th moments for some \(\delta>0\). Their results
approximate the tail probabilities at \(x^2\) uniformly for
\(0\le x\le b_n\), where \(b_n=o(n^{1/6})\). They conjectured that
the same range holds under a finite third moment when \(d\ge2\).

We develop explicit finite-sample estimates for one- and two-sample
Hotelling statistics under finite \(p\)th moments for \(2<p\le3\),
allowing independent observations with different distributions in the
one-sample setting. Our results specify the dependence on the moment
sums and the range of \(x\), and include uniform Berry--Esseen bounds
and explicit constructions showing sharpness of the moderate
deviation error order over independent arrays.

To describe the main estimates, let
\(\mathbf X_1,\ldots,\mathbf X_n\) be independent random vectors in
\(\mathbb R^d\), where \(d\ge2\) is fixed, with common mean
\(\boldsymbol\mu_n\) and positive definite average covariance matrix
\[
 \Sigma_n=\frac1n\sum_{i=1}^n\operatorname{Cov}(\mathbf X_i).
\]
Let \(\bar{\mathbf X}_n\) denote the sample mean and
\(\widehat\Sigma_n\) the centered sample covariance matrix:
\[
 \bar{\mathbf X}_n=\frac1n\sum_{i=1}^n\mathbf X_i,\qquad
 \widehat\Sigma_n=\frac1n\sum_{i=1}^n
  (\mathbf X_i-\bar{\mathbf X}_n)
  (\mathbf X_i-\bar{\mathbf X}_n)^\top.
\]
For a given center \(\boldsymbol\mu_0\), define the
one-sample Hotelling statistic, when \(\widehat\Sigma_n\) is positive
definite, by
\[
 T^2_{H,n}=n(\bar{\mathbf X}_n-\boldsymbol\mu_0)^\top
          \widehat\Sigma_n^{-1}
          (\bar{\mathbf X}_n-\boldsymbol\mu_0).
\]
Set \(T^2_{H,n}=0\) when \(\widehat\Sigma_n\) is singular.
Fix \(2<p\le3\) and assume that all terms in the moment sum below
are finite. Under \(H_0:\boldsymbol\mu_n=\boldsymbol\mu_0\),
Theorem~\ref{thm:hotelling} gives the Cram\'er-type moderate
deviation result
\begin{equation}\label{eq:intro-CMD}
 \left|
 \frac{\mathbb P(T^2_{H,n}\ge x^2)}
      {\mathbb P(\chi_d^2\ge x^2)}-1
 \right|
 \le C_d(1+x)^p n^{-p/2}
 \sum_{i=1}^n\mathbb E
 \left|\Sigma_n^{-1/2}
       (\mathbf X_i-\boldsymbol\mu_n)\right|^p
\end{equation}
for every \(x\ge0\) satisfying
\[
 (1+x)^p n^{-p/2}\sum_{i=1}^n
 \mathbb E\left|\Sigma_n^{-1/2}
       (\mathbf X_i-\boldsymbol\mu_n)\right|^p
 \le c_d.
\]
Here \(\chi_d^2\) has the chi-squared distribution with \(d\)
degrees of freedom, \(|\cdot|\) is the Euclidean norm, and
\(c_d>0\) and \(C_d<\infty\) depend only on \(d\), uniformly for \(2<p\le3\).

For the two-sample problem, let
\(\mathbf X_1,\ldots,\mathbf X_{n_1}\) and
\(\mathbf Y_1,\ldots,\mathbf Y_{n_2}\) be independent samples,
each internally i.i.d., with means
\(\boldsymbol\mu_X,\boldsymbol\mu_Y\) and positive definite
covariance matrices \(\Sigma_X,\Sigma_Y\), respectively.
Write \(\bar{\mathbf X},\bar{\mathbf Y}\) for the sample means and
\(\widehat\Sigma_X,\widehat\Sigma_Y\) for the centered sample
covariance matrices computed with divisors \(n_1,n_2\).
When \(\widehat\Sigma_X/n_1+\widehat\Sigma_Y/n_2\)
is positive definite, define
\[
 T_{n_1,n_2}^2=(\bar{\mathbf X}-\bar{\mathbf Y})^\top
 \left(\frac{\widehat\Sigma_X}{n_1}
       +\frac{\widehat\Sigma_Y}{n_2}\right)^{-1}
 (\bar{\mathbf X}-\bar{\mathbf Y}).
\]
Set \(T_{n_1,n_2}^2=0\) when this sample denominator
is singular.
The population covariances may differ, and no balance condition
on \(n_1,n_2\) is imposed. Assume that the two standardized
\(p\)th moments below are finite. Under
\(H_0:\boldsymbol\mu_X=\boldsymbol\mu_Y\), Theorem~\ref{thm:hotelling} gives
\begin{equation}\label{eq:intro-CMD-two}
 \begin{aligned}
 \left|
 \frac{\mathbb P(T_{n_1,n_2}^2\ge x^2)}
      {\mathbb P(\chi_d^2\ge x^2)}-1
 \right|
 \le C_d(1+x)^p
 \left\{
 \frac{\mathbb E|\Sigma_X^{-1/2}
                   (\mathbf X_1-\boldsymbol\mu_X)|^p}
      {n_1^{(p-2)/2}}
 +
 \frac{\mathbb E|\Sigma_Y^{-1/2}
                   (\mathbf Y_1-\boldsymbol\mu_Y)|^p}
      {n_2^{(p-2)/2}}
 \right\}
 \end{aligned}
\end{equation}
for every \(x\ge0\) satisfying
\[
 (1+x)^p
 \left\{
 \frac{\mathbb E|\Sigma_X^{-1/2}
                   (\mathbf X_1-\boldsymbol\mu_X)|^p}
      {n_1^{(p-2)/2}}
 +
 \frac{\mathbb E|\Sigma_Y^{-1/2}
                   (\mathbf Y_1-\boldsymbol\mu_Y)|^p}
      {n_2^{(p-2)/2}}
 \right\}
 \le c_d.
\]

Berry--Esseen bounds for Student's statistic under finite third
moments were established by \citet{BentkusGotze1996} and
\citet{BentkusBloznelisGotze1996}; \citet{Shao2005} obtained an explicit
bound using Stein's method. Distributional expansions for Hotelling's
statistic under nonnormality were studied by \citet{Kano1995} and
\citet{Fujikoshi1997}, while \citet{PinelisMolzon2016} treated
noncentral Hotelling statistics in a multivariate delta-method
framework.
For the one- and two-sample Hotelling statistics considered here,
uniform Berry--Esseen bounds follow directly from the moderate
deviation estimates above; see Corollary~\ref{cor:unif}.

For i.i.d. observations with uniformly bounded standardized \(p\)th
moments, the one-sample Berry--Esseen bound is of order
\(n^{-(p-2)/2}\), and the Cram\'er-type moderate deviation
approximation holds uniformly for \(0\le x\le b_n\) whenever
\[
 b_n=o\!\left(n^{(p-2)/(2p)}\right).
\]
For two samples, the Berry--Esseen bound is of order
\((n_1\wedge n_2)^{-(p-2)/2}\), and the Cram\'er-type moderate
deviation approximation holds uniformly for
\(0\le x\le b_{n_1,n_2}\) whenever
\[
 b_{n_1,n_2}=o\!\left((n_1\wedge n_2)^{(p-2)/(2p)}\right).
\]
At \(p=3\), the one-sample Berry--Esseen rate is \(n^{-1/2}\)
and the moderate deviation range is \(x=o(n^{1/6})\).
For two samples, the rate is \((n_1\wedge n_2)^{-1/2}\) and the
range is \(x=o((n_1\wedge n_2)^{1/6})\);
see Corollary~\ref{cor:iid}.

We construct examples showing that the error orders in the
one- and two-sample Cram\'er-type moderate deviation results are optimal.
See Examples~\ref{thm:sharpness} and~\ref{thm:two-sample-sharpness}.

We also apply the general result to simultaneous covariance testing
motivated by phenome-wide association studies. For bounded genetic
markers, disjoint pairing permits inference on untransformed phenotype
counts under finite \(p\)th moments, with a quantitative bound on the
familywise error rate.

The remainder of this paper is organized as follows.
Section~\ref{sec:main} presents the Cram\'er-type moderate deviation
results, Berry--Esseen bounds and sharpness results.
Section~\ref{sec:apps} gives the application and a simulation study.
Section~\ref{sec:proof} contains the proof of the main results.
Section~\ref{sec:sec4-lemma-proofs} proves the lemmas
used in the Stein argument and the application result.

\section{Main results}\label{sec:main}

For random vectors \(\boldsymbol\xi_1,\ldots,\boldsymbol\xi_n\)
in \(\mathbb R^d\), let
\begin{equation}\label{def:Sn-Vn}
 S_n=\sum_{i=1}^n\boldsymbol\xi_i,\qquad
 V_n=\sum_{i=1}^n\boldsymbol\xi_i\boldsymbol\xi_i^\top.
\end{equation}
We define the associated self-normalized quadratic form by
\begin{equation}\label{def:Rn}
 R_n^2=\sup_{\theta\in\mathbb R^d}
 \left\{2\theta^\top S_n-\theta^\top V_n\theta\right\}.
\end{equation}
The following theorem gives a Cram\'er-type moderate deviation
result for \(R_n^2\).

\begin{theorem}\label{thm:main-CMD-T}
Fix \(d\ge2\) and \(2<p\le3\). Let
\(\boldsymbol\xi_1,\ldots,\boldsymbol\xi_n\) be independent
random vectors in \(\mathbb R^d\) such that
\[
 \mathbb E\boldsymbol\xi_i=0,\qquad
 \sum_{i=1}^n\mathbb E
       (\boldsymbol\xi_i\boldsymbol\xi_i^\top)=I_d,\qquad
 \sum_{i=1}^n\mathbb E|\boldsymbol\xi_i|^p<\infty.
\]
Here \(|\cdot|\) denotes the Euclidean norm.
Then there exist positive constants
\(c_d\) and \(C_d\), depending only on \(d\), uniformly for \(2<p\le3\), such that
\begin{equation}\label{eq:CMD-Rn-1}
 \left|\frac{\mathbb P(R_n^2\ge x^2)}{Q_d(x^2)}-1\right|
 \le C_d(1+x)^p
 \sum_{i=1}^n\mathbb E|\boldsymbol\xi_i|^p,
\end{equation}
for every \(x\ge0\) satisfying
\begin{equation}\label{eq:CMD-x-range}
 (1+x)^p\sum_{i=1}^n\mathbb E|\boldsymbol\xi_i|^p
 \le c_d.
\end{equation}
Here \(Q_d(t)=\mathbb P(\chi_d^2\ge t)\), where \(\chi_d^2\)
has the chi-squared distribution with \(d\) degrees of freedom.
\end{theorem}

Theorem~\ref{thm:main-CMD-T} yields Cram\'er-type moderate
deviation results for the one- and two-sample Hotelling statistics.
For one sample, let
\(\mathbf X_1,\ldots,\mathbf X_n\) be random vectors in
\(\mathbb R^d\) with common mean \(\boldsymbol\mu_n\).
For a specified vector \(\boldsymbol\mu_0\), we consider testing
\[
 H_0:\boldsymbol\mu_n=\boldsymbol\mu_0.
\]
Let \(\bar{\mathbf X}_n\) and \(\widehat\Sigma_n\) denote the
sample mean and centered sample covariance matrix:
\[
 \bar{\mathbf X}_n=\frac1n\sum_{i=1}^n\mathbf X_i,\qquad
 \widehat\Sigma_n=\frac1n\sum_{i=1}^n
  (\mathbf X_i-\bar{\mathbf X}_n)
  (\mathbf X_i-\bar{\mathbf X}_n)^\top.
\]
When \(\widehat\Sigma_n\) is positive definite, the one-sample
Hotelling statistic for testing this hypothesis is
\begin{equation}\label{def:Hotelling-T}
 T^2_{H,n}=\sup_{|\theta|=1}
 \frac{n\{\theta^\top(\bar{\mathbf X}_n-\boldsymbol\mu_0)\}^2}
      {\theta^\top\widehat\Sigma_n\theta}.
\end{equation}
Set \(T^2_{H,n}=0\) when \(\widehat\Sigma_n\) is singular.
Write $\Sigma_n=\frac1n\sum_{i=1}^n\operatorname{Cov}(\mathbf X_i)$
for the average population covariance matrix.
When \(\Sigma_n\) is positive definite, define the standardized
observations by
\begin{equation}\label{def:xi}
 \boldsymbol\xi_i=n^{-1/2}\Sigma_n^{-1/2}
                  (\mathbf X_i-\boldsymbol\mu_n).
\end{equation}
These vectors are centered and their covariance matrices sum
to \(I_d\). For \(R_n^2\) formed from these vectors as in
\eqref{def:Rn}, the one-sample statistic satisfies, under \(H_0\)
and on the event \(\{\widehat\Sigma_n\succ0\}\),
\[
 T^2_{H,n}=\frac{R_n^2}{1-R_n^2/n};
\]
see the proof of Theorem~\ref{thm:hotelling}\textup{(i)}.

For two samples, let
\(\mathbf X_1,\ldots,\mathbf X_{n_1}\) and
\(\mathbf Y_1,\ldots,\mathbf Y_{n_2}\) have population means
\(\boldsymbol\mu_X,\boldsymbol\mu_Y\) and covariance matrices
\(\Sigma_X,\Sigma_Y\), respectively. We consider testing
\[
 H_0:\boldsymbol\mu_X=\boldsymbol\mu_Y.
\]
Write \(\bar{\mathbf X},\bar{\mathbf Y}\) for the sample means
and define the centered sample covariance matrices by
\[
 \widehat\Sigma_X=\frac1{n_1}\sum_{i=1}^{n_1}
  (\mathbf X_i-\bar{\mathbf X})(\mathbf X_i-\bar{\mathbf X})^\top,
 \qquad
 \widehat\Sigma_Y=\frac1{n_2}\sum_{j=1}^{n_2}
  (\mathbf Y_j-\bar{\mathbf Y})(\mathbf Y_j-\bar{\mathbf Y})^\top.
\]
When \(\widehat\Sigma_X/n_1+\widehat\Sigma_Y/n_2\) is positive
definite, the two-sample Hotelling statistic for testing this
hypothesis is
\begin{equation}\label{def:two-sample-T}
 T_{n_1,n_2}^2=\sup_{|\theta|=1}
 \frac{\{\theta^\top(\bar{\mathbf X}-\bar{\mathbf Y})\}^2}
      {\theta^\top\widehat\Sigma_X\theta/n_1
       +\theta^\top\widehat\Sigma_Y\theta/n_2}.
\end{equation}
Set \(T_{n_1,n_2}^2=0\) when
\(\widehat\Sigma_X/n_1+\widehat\Sigma_Y/n_2\) is singular.
The following theorem gives Cram\'er-type moderate deviation
results for both Hotelling statistics.

\begingroup
\begin{theorem}\label{thm:hotelling}
Fix \(d\ge2\) and \(2<p\le3\). There exist positive constants
\(c_d\) and \(C_d\), depending only on \(d\), uniformly for \(2<p\le3\), such that
the following statements hold.

\textup{(i)} Suppose that \(\mathbf X_1,\ldots,\mathbf X_n\)
are independent random vectors in \(\mathbb R^d\) with common
mean \(\boldsymbol\mu_n\) and finite \(p\)th moments.
Assume that \(\Sigma_n\) is positive definite.
Under \(H_0:\boldsymbol\mu_n=\boldsymbol\mu_0\), we have
\begin{equation}\label{eq:CMD-T-1}
 \begin{aligned}
 \left|\frac{\mathbb P(T^2_{H,n}\ge x^2)}{Q_d(x^2)}-1\right|
 \le C_d(1+x)^p n^{-p/2}
 \sum_{i=1}^n\mathbb E
 \left|\Sigma_n^{-1/2}(\mathbf X_i-\boldsymbol\mu_n)\right|^p,
 \end{aligned}
\end{equation}
for every \(x\ge0\) satisfying
\[
 (1+x)^p n^{-p/2}\sum_{i=1}^n\mathbb E
 \left|\Sigma_n^{-1/2}(\mathbf X_i-\boldsymbol\mu_n)\right|^p
 \le c_d.
\]

\textup{(ii)} Suppose that the two samples are independent and
each is internally i.i.d., with respective means
\(\boldsymbol\mu_X,\boldsymbol\mu_Y\), positive definite
covariance matrices \(\Sigma_X,\Sigma_Y\), and finite \(p\)th
moments. Under
\(H_0:\boldsymbol\mu_X=\boldsymbol\mu_Y\), we have
\begin{equation}\label{eq:CMD-two-1}
 \begin{aligned}
 \left|\frac{\mathbb P(T_{n_1,n_2}^2\ge x^2)}{Q_d(x^2)}-1\right|
 \le C_d(1+x)^p
 \left\{
 \frac{\mathbb E|\Sigma_X^{-1/2}
              (\mathbf X_1-\boldsymbol\mu_X)|^p}{n_1^{(p-2)/2}}
 +
 \frac{\mathbb E|\Sigma_Y^{-1/2}
              (\mathbf Y_1-\boldsymbol\mu_Y)|^p}{n_2^{(p-2)/2}}
 \right\},
 \end{aligned}
\end{equation}
for every \(x\ge0\) satisfying
\[
 (1+x)^p
 \left\{
 \frac{\mathbb E|\Sigma_X^{-1/2}
              (\mathbf X_1-\boldsymbol\mu_X)|^p}{n_1^{(p-2)/2}}
 +
 \frac{\mathbb E|\Sigma_Y^{-1/2}
              (\mathbf Y_1-\boldsymbol\mu_Y)|^p}{n_2^{(p-2)/2}}
 \right\}
 \le c_d.
\]
\end{theorem}
\endgroup

The moderate deviation results imply the following uniform
Berry--Esseen bounds.

\begin{corollary}\label{cor:unif}
For \(R_n^2\), retain the assumptions of
Theorem~\ref{thm:main-CMD-T}. For the two Hotelling statistics,
retain the corresponding model and moment assumptions and null
hypotheses of Theorem~\ref{thm:hotelling}. The following bounds hold with a constant \(C_d\) depending only on \(d\):
\begin{equation}\label{eq:BE-all}
 \sup_{t\ge0}\left|\mathbb P(R_n^2\ge t)-Q_d(t)\right|
 \le C_d\sum_{i=1}^n\mathbb E|\boldsymbol\xi_i|^p,
\end{equation}
\begin{equation}\label{eq:BE-one}
 \begin{aligned}
 \sup_{t\ge0}\left|\mathbb P(T^2_{H,n}\ge t)-Q_d(t)\right|
 \le C_dn^{-p/2}\sum_{i=1}^n\mathbb E
 \left|\Sigma_n^{-1/2}(\mathbf X_i-\boldsymbol\mu_n)\right|^p,
 \end{aligned}
\end{equation}
and
\begin{equation}\label{eq:BE-two}
 \begin{aligned}
 \sup_{t\ge0}\left|\mathbb P(T_{n_1,n_2}^2\ge t)-Q_d(t)\right|
 \le C_d
 \left\{
 \frac{\mathbb E|\Sigma_X^{-1/2}
              (\mathbf X_1-\boldsymbol\mu_X)|^p}{n_1^{(p-2)/2}}
 +
 \frac{\mathbb E|\Sigma_Y^{-1/2}
              (\mathbf Y_1-\boldsymbol\mu_Y)|^p}{n_2^{(p-2)/2}}
 \right\}.
 \end{aligned}
\end{equation}
\end{corollary}

For i.i.d. samples, the preceding bounds yield the following rates.

\begin{corollary}\label{cor:iid}
Under the model assumptions and corresponding null hypothesis
of Theorem~\ref{thm:hotelling}, suppose each sample is i.i.d.,
with standardized \(p\)th moments uniformly bounded in the sample sizes.

\textup{(i)} For one sample, as \(n\to\infty\),
\begin{equation}\label{eq:iid-BE-rate}
 \sup_{t\ge0}|\mathbb P(T^2_{H,n}\ge t)-Q_d(t)|
 =O\!\left(n^{-(p-2)/2}\right),
\end{equation}
and
\begin{equation}\label{eq:iid-CMD-rate}
 \sup_{0\le x\le b_n}
 \left|\frac{\mathbb P(T^2_{H,n}\ge x^2)}{Q_d(x^2)}-1\right|
 =O\!\left((1+b_n)^p n^{-(p-2)/2}\right)=o(1),
\end{equation}
for every nonnegative sequence \(b_n=o(n^{(p-2)/(2p)})\).

\textup{(ii)} For two samples, as \(n_1\wedge n_2\to\infty\),
\begin{equation}\label{eq:iid-BE-two-rate}
 \sup_{t\ge0}|\mathbb P(T_{n_1,n_2}^2\ge t)-Q_d(t)|
 =O\!\left((n_1\wedge n_2)^{-(p-2)/2}\right),
\end{equation}
and
\begin{equation}\label{eq:iid-CMD-two-rate}
 \begin{aligned}
 &\sup_{0\le x\le b_{n_1,n_2}}
 \left|\frac{\mathbb P(T_{n_1,n_2}^2\ge x^2)}{Q_d(x^2)}-1\right|
 =O\!\left((1+b_{n_1,n_2})^p
             (n_1\wedge n_2)^{-(p-2)/2}\right)=o(1),
 \end{aligned}
\end{equation}
for every nonnegative \(b_{n_1,n_2}
=o((n_1\wedge n_2)^{(p-2)/(2p)})\).
The implied constants depend only on \(d\) and the uniform moment bound.
\end{corollary}

At \(p=3\), the one-sample Berry--Esseen rate is \(n^{-1/2}\)
and the Cram\'er-type moderate deviation range is \(x=o(n^{1/6})\).
For two samples, the Berry--Esseen rate is \((n_1\wedge n_2)^{-1/2}\)
and the range is \(x=o((n_1\wedge n_2)^{1/6})\).

The following two examples show that the error orders in the
one- and two-sample Cram\'er-type moderate deviation results are sharp.

\begingroup
\begin{example}\label{thm:sharpness}
Fix \(d\ge2\) and \(2<p\le3\). There exist \(n_k,x_k\to\infty\)
and a triangular array of independent centered vectors
\(\mathbf X_{k1},\ldots,\mathbf X_{k n_k}\), with
\(n_k^{-1}\sum_{i=1}^{n_k}
\mathbb E(\mathbf X_{ki}\mathbf X_{ki}^{\top})=I_d\)
and almost surely positive definite centered sample covariance.
Let \(T^2_{H,n_k}\) be their one-sample Hotelling statistic with
null center zero. Then, for some \(c_{d,p}\in(0,\infty)\),
\begin{equation}\label{eq:one-sample-sharp-rate}
 \left|
  \frac{\mathbb P(T^2_{H,n_k}\ge x_k^2)}{Q_d(x_k^2)}-1
 \right|
 \sim c_{d,p}(1+x_k)^p n_k^{-p/2}
 \sum_{i=1}^{n_k}\mathbb E|\mathbf X_{ki}|^p
 \longrightarrow0.
\end{equation}
Thus the moderate deviation error order is sharp over independent
arrays.
\end{example}
\endgroup

\begingroup
\begin{example}\label{thm:two-sample-sharpness}
Fix \(d\ge2\) and \(2<p\le3\). There exist \(n_k,x_k\to\infty\)
and two triangular arrays
\(\mathbf X_{k1},\ldots,\mathbf X_{k n_k}\) and
\(\mathbf Y_{k1},\ldots,\mathbf Y_{k n_k}\), whose \(2n_k\)
observations in each row are i.i.d., with mean zero, covariance
\(\Sigma_k\succ0\), and
\(\sup_k\mathbb E|\Sigma_k^{-1/2}\mathbf X_{k1}|^p<\infty\).
Both centered sample covariance matrices are positive definite
almost surely. Let \(T_{n_k,n_k}^2\) be their two-sample Hotelling
statistic. Then, for some \(c_{d,p}\in(0,\infty)\),
\begin{equation}\label{eq:two-sample-sharp-rate}
 \left|
  \frac{\mathbb P(T_{n_k,n_k}^2\ge x_k^2)}{Q_d(x_k^2)}-1
 \right|
 \sim c_{d,p}(1+x_k)^p n_k^{-(p-2)/2}\longrightarrow0.
\end{equation}
Thus the two-sample moderate deviation error order is sharp even
for equal sample sizes and identical populations.
\end{example}
\endgroup

\section{Applications}\label{sec:apps}

Phenome-wide association studies (PheWAS) examine associations
between genetic markers and many disease phenotypes.
Following \citet{CaiCaiLiaoLiu2019}, let
\((Y_i,X_{i1},\ldots,X_{iM})\), \(1\le i\le n\), be independent
and identically distributed (i.i.d.) observations from \(n\ge2\)
subjects. Here \(Y_i\in\mathbb R^d\) records \(d\ge2\) genetic
markers, with \(d\) fixed, and \(X_{ij}\) counts medical encounters
for disease \(j\), \(1\le j\le M\), with \(M\ge1\).
A single-nucleotide polymorphism (SNP) is a variation at one genetic
position. For SNP markers, each coordinate of \(Y_i\) counts
\(0,1,2\) copies of a specified variant (allele) at that position,
so \(Y_i\) is bounded.
Assuming finite second moments, we test the covariance hypotheses
\[
 H_{0j}:\operatorname{Cov}(Y,X_j)=0,\qquad 1\le j\le M.
\]

To remove the unknown means while retaining independence, fix
disjoint pairs of subjects before observing the data and define,
with \(N=\lfloor n/2\rfloor\),
\begin{equation}\label{eq:apps-paired-vectors}
 Z_{kj}=\frac12(Y_{2k}-Y_{2k-1})(X_{2k,j}-X_{2k-1,j}),
 \qquad 1\le k\le N.
\end{equation}
For each \(j\), these vectors are i.i.d.\ and satisfy
\(\mathbb EZ_{kj}=\operatorname{Cov}(Y,X_j)\).
Thus each covariance hypothesis becomes a mean-zero hypothesis.
Use the statistic
\begin{equation}\label{eq:apps-statistic}
 U_j=\sup_{\theta\in\mathbb R^d}
 \sum_{k=1}^N\{2\theta^\top Z_{kj}-(\theta^\top Z_{kj})^2\}.
\end{equation}
This is the statistic in \eqref{def:Rn} applied to the paired vectors
and is defined even when their empirical second-moment matrix is singular.
Theorem~\ref{thm:main-CMD-T} and Bonferroni correction give the
following bound for the family-wise error rate (FWER), the probability
of rejecting at least one true null hypothesis.
Let \(\mathcal H_0\) index the true null hypotheses.

\begin{corollary}\label{cor:phewas-bonferroni}
Let \(2<p\le3\) and \(K>0\). Under the sampling model above, suppose that
\(\Gamma_j=\operatorname{Cov}(Z_{1j})\succ0\) and
\[
 \mathbb E|\Gamma_j^{-1/2}Z_{1j}|^p\le K
 \qquad\text{for every }j\in\mathcal H_0.
\]
For \(0<\alpha<1\), let \(q_{\alpha,M}\) satisfy
\(Q_d(q_{\alpha,M})=\alpha/M\), and reject \(H_{0j}\) when
\(U_j\ge q_{\alpha,M}\). Then
\begin{equation}\label{eq:apps-fwer}
 \mathbb P\{\text{at least one true null is rejected}\}
 \le\alpha\left[
 1+C_dK N^{-(p-2)/2}
       \{1+\log(M/\alpha)\}^{p/2}\right],
\end{equation}
provided
\[
 K N^{-(p-2)/2}\{1+\log(M/\alpha)\}^{p/2}\le c_d.
\]
The positive constants depend only on \(d\), uniformly for
\(2<p\le3\). No independence assumption across phenotypes is needed.
\end{corollary}

For bounded markers whose covariance matrix has smallest eigenvalue
bounded away from zero, the corollary's moment and nondegeneracy
assumptions hold when the null phenotypes have positive variances
and uniformly bounded standardized \(p\)th moments.
Thus it permits inference on the original counts
even when third moments are infinite. The proof is in
Section~\ref{sec:sec4-lemma-proofs}; it also shows that, for fixed
\(\alpha\), rejection probability under each fixed alternative tends
to one if \(\log M=o(N)\), \(\mathbb E|Z_{1j}|^2<\infty\), and
\(\mathbb E(Z_{1j}Z_{1j}^\top)\succ0\).

For comparison, the full-sample statistic in equation~(2.2) of
\citet{CaiCaiLiaoLiu2019} is constructed from
\[
 \widehat Z_{ij}=(Y_i-\bar Y)(X_{ij}-\bar X_j),\qquad
 \widehat{\boldsymbol\sigma}_j=\frac1n\sum_{i=1}^n\widehat Z_{ij},
\]
where \(\bar Y\) and \(\bar X_j\) are the sample means. Define
\[
 \widehat\Sigma_{Zj}=\frac1n\sum_{i=1}^n
 (\widehat Z_{ij}-\widehat{\boldsymbol\sigma}_j)
 (\widehat Z_{ij}-\widehat{\boldsymbol\sigma}_j)^\top.
\]
Their statistic is
\[
 T_j=n\widehat{\boldsymbol\sigma}_j^\top
       \widehat\Sigma_{Zj}^{-1}\widehat{\boldsymbol\sigma}_j.
\]
They study the false discovery proportion (FDP),
the fraction of rejected hypotheses that are true nulls (zero if
none are rejected), and its expectation, the false discovery rate (FDR).
Theorem~1 of \citet{CaiCaiLiaoLiu2019} establishes asymptotic
control of the false discovery rate and false discovery proportion
under conditions including moments of order \(8b+4+\epsilon\)
when \(M\le n^b\).
For bounded genetic markers, Corollary~\ref{cor:phewas-bonferroni}
gives a quantitative family-wise error bound for our paired
procedure under standardized \(p\)th moment assumptions,
\(2<p\le3\), subject to the conditions stated there.
The logarithmic transformation \(\log(1+X_j)\) used in their data
analysis generally changes the covariance target; our paired
procedure uses the original counts.

We illustrate performance using Poisson mixtures, which model
variation in count intensities and have been used for healthcare
utilization data \citep{GoncalvesBarretoSouzaOmbao2025}.
In the simulations, we compare Bonferroni tests based on the
paired statistic \(U_j\) and the full-sample statistic \(T_j\)
of \citet{CaiCaiLiaoLiu2019}, using the same chi-squared threshold.
This comparison does not use their original FDR procedure.
For the simulations, \(T_j\) is set to zero when
\(\widehat\Sigma_{Zj}\) is numerically singular.
Take \(n=1000\), \(d=2\), \(M=100\),
\(\alpha=0.05\), and independent marker counts
\(Y_{ir}\sim\operatorname{Binomial}(2,0.1)\).
Conditional on the markers and mixing variables, generate independent
counts
\[
 X_{ij}\mid Y_i,W_{ij}\sim\operatorname{Poisson}(\lambda_{ij}W_{ij}),
 \qquad
 \lambda_{ij}=\begin{cases}
 5e^{0.6Y_{i1}}/(0.9+0.1e^{0.6})^2,&j\le10,\\
 5,&j>10.
 \end{cases}
\]
The \(W_{ij}\) are i.i.d.\ and independent of the markers.
We use three illustrative mixing laws, each with shape \(5/2\)
and mean one: Pareto with minimum \(3/5\), inverse-gamma with
scale \(3/2\), and log-logistic with scale
\(5\sin(2\pi/5)/(2\pi)\). Their variances differ.
The resulting counts have finite moments of order \(9/4\) but infinite
third moments; the first \(10\) diseases have nonzero covariance
with the markers, and the remaining \(90\) are true nulls.

After exploratory selection, we evaluate these settings with
\(10{,}000\) new independent replications. Table~\ref{tab:phewas-simulation}
reports empirical FWER and power, the mean fraction of the \(10\)
alternatives rejected. In the three heavy-tailed settings considered,
the paired tests have empirical FWER below \(0.05\) and power between
\(0.68\) and \(0.95\). The Bonferroni tests based on \(T_j\) have
higher empirical power, but their empirical FWER ranges from
\(0.10\) to \(0.14\). A matched Poisson control
\((W_{ij}=1)\) gives FWER \(0.0305\) and \(0.0417\), respectively,
with power one for both. All exploratory results, including settings
with lower power, are retained with the simulation code.

\begin{table}[!ht]
\centering
\begin{tabular}{lrrrr}
\hline
 & \multicolumn{2}{c}{Paired $U$} & \multicolumn{2}{c}{Full sample} \\
Mixing distribution & FWER & Power & FWER & Power \\
\hline
Pareto & 0.015 & 0.949 & 0.101 & 0.985 \\
Inverse-gamma & 0.012 & 0.680 & 0.140 & 0.913 \\
Log-logistic & 0.014 & 0.804 & 0.114 & 0.954 \\
\hline
\end{tabular}
\caption{Bonferroni tests at level \(0.05\), with \(1000\) subjects,
\(100\) phenotypes, \(10\) signals and \(10{,}000\) independent
confirmation replications per setting. Monte Carlo standard errors
are below \(0.0035\) for FWER and \(0.0016\) for power.}
\label{tab:phewas-simulation}
\end{table}


\section{Proof of the main results}\label{sec:proof}

\begin{proof}[Proof of Theorem~\ref{thm:main-CMD-T}]
We use Stein's method to prove \eqref{eq:CMD-Rn-1}.
\begingroup
The positive constants \(C_d,c_d\) depend only on \(d\) and may
vary from line to line unless otherwise specified.
\endgroup
Let \(\chi_d=(\chi_d^2)^{1/2}\) and \(R_n=(R_n^2)^{1/2}\).
For \(x\ge0\), we have
\(\mathbb P(R_n^2\ge x^2)=\mathbb P(R_n\ge x)\) and
\(\mathbb P(\chi_d^2\ge x^2)=\mathbb P(\chi_d\ge x)\).
It therefore suffices to compare the tails of \(R_n\) and \(\chi_d\).
The density of \(\chi_d\) is
\[
 p_d(y)=\frac{y^{d-1}e^{-y^2/2}}
 {2^{d/2-1}\Gamma(d/2)},\qquad y>0.
\]
To derive a Stein identity, let \(a(\cdot)\) be a continuously
differentiable function with compact support in \((0,\infty)\).
Since \(p_d'(y)/p_d(y)=(d-1)/y-y\), we have
\[
 \frac{\{y p_d(y)a(y)\}'}{p_d(y)}
 =a(y)+ya'(y)+(d-1-y^2)a(y)
 =ya'(y)+(d-y^2)a(y).
\]
Multiplying by \(p_d(y)\) and integrating over \((0,\infty)\), we obtain
\[
\begin{aligned}
 \mathbb E\{\chi_d a'(\chi_d)+(d-\chi_d^2)a(\chi_d)\}
 &=\int_0^\infty\{ya'(y)+(d-y^2)a(y)\}p_d(y)\,dy\\
 &=\int_0^\infty\{y p_d(y)a(y)\}'\,dy\\
 &=\bigl[y p_d(y)a(y)\bigr]_0^\infty=0.
\end{aligned}
\]
The last equality follows from the compact support of \(a(\cdot)\):
\(a(y)=0\) both near zero and for all sufficiently large \(y\),
so both boundary terms vanish.
To compare expectations of a smooth function
\(H:[0,\infty)\to[0,1]\), this identity leads to the radial Stein equation
\begin{equation}\label{eq:sec4-radial-equation}
ya_H'(y)+(d-y^2)a_H(y)=H(y)-\mathbb EH(\chi_d).
\end{equation}
Multiplying by \(p_d(y)\) and integrating from zero to \(y\), with
\(a_H\) continuous at zero, gives
\[
 y p_d(y)a_H(y)
 =\int_0^y\{H(u)-\mathbb EH(\chi_d)\}p_d(u)\,du.
\]
Dividing by \(y p_d(y)\) and substituting the density formula yields
\begin{equation}\label{eq:sec4-radial-solution}
a_H(y)=y^{-d}e^{y^2/2}
\int_0^y\{H(u)-\mathbb EH(\chi_d)\}u^{d-1}e^{-u^2/2}\,du,
\qquad y>0.
\end{equation}
We extend \(a_H\) continuously to \(y=0\) by defining
\[
 a_H(0)=\lim_{y\downarrow0}a_H(y)
 =\frac{H(0)-\mathbb EH(\chi_d)}{d}.
\]
For a vector \(s\in\mathbb R^d\) and a symmetric positive definite
matrix \(v\in\mathbb R^{d\times d}\), write
\(y=(s^\top v^{-1}s)^{1/2}\) and define
\begin{equation}\label{eq:sec4-stein-field}
f_H(s,v)=v^{-1}s\,
a_H(y).
\end{equation}
Let \(\partial_s f_H(s,v)\) denote the Jacobian matrix with respect
to \(s\). For \(s\ne0\), we have \(\nabla_s y=v^{-1}s/y\), so the
product and chain rules give
\[
 \partial_s f_H(s,v)
 =a_H(y)v^{-1}
 +\frac{a_H'(y)}{y}\,v^{-1}ss^\top v^{-1}.
\]
Since \(s^\top v^{-1}s=y^2\), taking the trace and using
\eqref{eq:sec4-radial-equation} yield
\begin{equation}\label{eq:sec4-stein-identity-matrix}
\begin{aligned}
\operatorname{tr}\{v\,\partial_sf_H(s,v)\}-s^\top f_H(s,v)
&=d\,a_H(y)+y a_H'(y)-y^2a_H(y)\\
&=y a_H'(y)+(d-y^2)a_H(y)\\
&=H(y)-\mathbb EH(\chi_d).
\end{aligned}
\end{equation}
At \(s=0\), continuity of \(a_H\) gives
\(\partial_s f_H(0,v)=a_H(0)v^{-1}\).
The left-hand side is therefore
\(d\,a_H(0)=H(0)-\mathbb EH(\chi_d)\), so
\eqref{eq:sec4-stein-identity-matrix} also holds at \(s=0\).
Recall from \eqref{def:Sn-Vn} that
$V_n=\sum_{j=1}^n\boldsymbol\xi_j\boldsymbol\xi_j^\top$.
Removing the contribution of the \(i\)th observation gives
$V_n-\boldsymbol\xi_i\boldsymbol\xi_i^\top
 =\sum_{j\ne i}\boldsymbol\xi_j\boldsymbol\xi_j^\top$.
Both matrices are positive semidefinite and may be singular.
Since \(f_H(s,v)\) in \eqref{eq:sec4-stein-field} is defined only
for positive definite \(v\), we first construct a smooth modification
that is also defined for positive semidefinite matrix arguments.
For this construction, we use smooth cutoff functions.
For \(u\ge0\) and \(0<\eta\le1\), choose a smooth non-increasing
function \(H_{u,\eta}:\mathbb R\to[0,1]\) such that
\begin{equation}\label{eq:sec4-cutoff}
H_{u,\eta}(y)=1\quad(y\le u),\qquad
H_{u,\eta}(y)=0\quad(y\ge u+\eta),\qquad
|H_{u,\eta}'(y)|\le C_d/\eta.
\end{equation}
Here \(u\) is the start of the transition and \(\eta\) its width.
For a symmetric matrix \(v\), taking $u =1/16$  and $ \eta = 3/16$ in \eqref{eq:sec4-cutoff}, define
\begin{equation}\label{eq:sec4-extension}
\begin{aligned}
K(v)&=I_d+H_{1/16,\,3/16}(\|v-I_d\|_{\mathrm F}^2)(v-I_d),\\
F(s,v)&=f_H(s,K(v)).
\end{aligned}
\end{equation}
Here \(\|A\|_{\mathrm F}^2=\operatorname{tr}(A^\top A)\).
The eigenvalues of \(K(v)\) lie in \([1/2,3/2]\), and
\(K(v)=v\) if \(\|v-I_d\|_{\mathrm F}\le1/4\).
Since \(K(v)\succ0\), \(F(s,v)\) is well defined for every
\(s\in\mathbb R^d\) and \(v\succeq0\).
Substituting \(K(v)\) into \eqref{eq:sec4-stein-identity-matrix} gives,
for every \(v\succeq0\),
\begin{equation}\label{eq:sec4-stein-identity-regularized}
\operatorname{tr}\{K(v)\,\partial_sF(s,v)\}-s^\top F(s,v)
=H\bigl((s^\top K(v)^{-1}s)^{1/2}\bigr)-\mathbb EH(\chi_d).
\end{equation}
Substituting \((s,v)=(S_n,V_n)\), with \(S_n,V_n\) defined in
\eqref{def:Sn-Vn}, into \eqref{eq:sec4-stein-identity-regularized} gives
\[
\begin{aligned}
H\bigl((S_n^\top K(V_n)^{-1}S_n)^{1/2}\bigr)-\mathbb EH(\chi_d)
=\operatorname{tr}\{K(V_n)\,\partial_sF(S_n,V_n)\}
          -S_n^\top F(S_n,V_n).
\end{aligned}
\]
To obtain \(H(R_n)-\mathbb EH(\chi_d)\) on the left-hand side,
we add
\[
 H(R_n)-H\bigl((S_n^\top K(V_n)^{-1}S_n)^{1/2}\bigr)
\]
to both sides. Writing
\(K(V_n)=V_n+\{K(V_n)-V_n\}\) and using the linearity of the trace,
we obtain
\begin{equation}\label{eq:sec4-stein-identity-sample}
\begin{aligned}
H(R_n)-\mathbb EH(\chi_d)
&=\operatorname{tr}\{V_n\,\partial_sF(S_n,V_n)\}-S_n^\top F(S_n,V_n)\\
&\quad+H(R_n)-H\bigl((S_n^\top K(V_n)^{-1}S_n)^{1/2}\bigr)\\
&\quad+\operatorname{tr}\{(K(V_n)-V_n)\,\partial_sF(S_n,V_n)\}.
\end{aligned}
\end{equation}
The first term retains the sum representations of \(S_n\) and
\(V_n\), which will allow us to estimate its expectation by a
leave-one-out expansion. The last two terms arise from replacing
\(V_n\) by \(K(V_n)\). When \(\|V_n-I_d\|_{\mathrm F}\le1/4\),
we have \(K(V_n)=V_n\) and
\(R_n=(S_n^\top V_n^{-1}S_n)^{1/2}\), so both terms vanish.
We therefore need to estimate these two corrections only on the
event \(\|V_n-I_d\|_{\mathrm F}>1/4\).

Lemma~\ref{lem:sec4-field-global} below gives bounds for \(F\)
and its derivatives.
These bounds will be used in the leave-one-out expansion of the first
term and in estimating the two correction terms in
\eqref{eq:sec4-stein-identity-sample}.
The proofs of the lemmas and corollary below are given in
Section~\ref{sec:sec4-lemma-proofs}.

\begin{lemma}\label{lem:sec4-field-global}
Let $F$ be defined by \eqref{eq:sec4-extension}.
There is a constant \(C_d>0\), depending only on \(d\), such that
the following estimates hold. For every smooth
\(H:[0,\infty)\to[0,1]\) and all \(s\in\mathbb R^d\),
\(v\succeq0\), we have
\begin{equation}\label{eq:sec4-truncated-field}
|F(s,v)|+|s^\top F(s,v)|+\|\partial_sF(s,v)\|_{\mathrm F}
\le C_d,
\end{equation}
and
\begin{equation}\label{eq:sec4-global-matrix-derivative}
\|\partial_vF(s,v)\|
=\sup_{\substack{A=A^\top\\\|A\|_{\mathrm F}\le1}}
 |\partial_vF(s,v)[A]|
\le C_d(1+|s|).
\end{equation}
Here \(\partial_vF(s,v)[A]\) denotes the directional derivative
with respect to \(v\) in the symmetric matrix direction \(A\).
Since \(K(I_d)=I_d\), we have
\(F(\cdot,I_d)=f_H(\cdot,I_d)\).
For \(H=H_{u,\eta}\) in \eqref{eq:sec4-cutoff} and
\(w,w'\in\mathbb R^d\), put \(y=|w'|\) and \(\delta=|w-w'|\).
Then the difference of the Jacobians satisfies
\begin{equation}\label{eq:sec4-global-jacobian-modulus}
\begin{aligned}
\|\partial_sf_H(w,I_d)-\partial_sf_H(w',I_d)\|_{\mathrm F}
&\le C_d\delta(1+y+\delta)
 +C_d\min\left\{1,\frac{\delta}{y+\eta}\right\}\\
&\quad+C_d\min\left\{1,\frac{\delta}{\eta}\right\}
 \mathbf1_{\{|y-u|\le\delta+\eta\}}.
\end{aligned}
\end{equation}
\end{lemma}

We next estimate the expectations of the three terms in
\eqref{eq:sec4-stein-identity-sample}. The moment condition
\(\sum_{i=1}^n\mathbb E|\boldsymbol\xi_i|^p<\infty\) in
Theorem~\ref{thm:main-CMD-T} does not require the observations to be
bounded. Each observation contributes \(\boldsymbol\xi_i\) to
\(S_n\) and \(\boldsymbol\xi_i\boldsymbol\xi_i^\top\) to \(V_n\).
We first restrict the size of every observation so that both
contributions are controlled.

For \(x>4\), consider a change \(h\) in the radius from \(x\)
to \(x+h\). The exponential factor in the chi density satisfies
\[
 \frac{e^{-(x+h)^2/2}}{e^{-x^2/2}}
 =e^{-xh-h^2/2}.
\]
Thus this ratio stays close to one when \(x|h|\) is small.
This motivates choosing the observation cutoff on the scale
\(x^{-1}\); the resulting changes in the radius will be
estimated below. For \(0\le x\le4\), we use a fixed small cutoff.
The threshold \(4\) is chosen for convenience; another sufficiently
large fixed constant would work as well, after adjusting the constants
in the estimates below.
Fix \(\varepsilon=\varepsilon(d)\in(0,1]\) sufficiently small for the estimates
below. Let \(b\) be the cutoff and \(\mathcal N\) the event that
all observations satisfy it. We write \(\mathbb P^*\) and
\(\mathbb E^*\) for probability and expectation conditional on this
event, respectively:
\begin{equation}\label{eq:sec4-conditioning}
\begin{aligned}
b=\frac{\varepsilon}{\max\{4,x\}},\qquad
\mathcal N=\bigcap_{i=1}^n\{|\boldsymbol\xi_i|\le b\},
\\
\mathbb P^*(\cdot)=\mathbb P(\cdot\mid\mathcal N),\qquad
\mathbb E^*(\cdot)=\mathbb E(\cdot\mid\mathcal N).
\end{aligned}
\end{equation}
This choice gives \(b\le\varepsilon/4\) and \(xb\le\varepsilon\).
On \(\mathcal N\), every summand of \(S_n\) has norm at most \(b\),
and every summand of \(V_n\) has Frobenius norm at most \(b^2\),
since \(\|\boldsymbol\xi_i\boldsymbol\xi_i^\top\|_{\mathrm F}
=|\boldsymbol\xi_i|^2\).
For the conditional estimates, assume
\(b^{-p}\sum_i\mathbb E|\boldsymbol\xi_i|^p\le1/4\);
this condition will be checked when applying the estimates below.
The union bound gives
\begin{equation}\label{eq:sec4-conditioning-bound}
\mathbb P(\mathcal N^c)
\le b^{-p}\sum_i\mathbb E|\boldsymbol\xi_i|^p\le\frac14,
\qquad\mathbb P(\mathcal N)\ge\frac34.
\end{equation}
Thus the conditional law is well defined. Each event
\(\{|\boldsymbol\xi_i|\le b\}\) in the definition of \(\mathcal N\)
depends only on \(\boldsymbol\xi_i\). By independence under
\(\mathbb P\), for any Borel sets \(B_1,\ldots,B_n\subseteq\mathbb R^d\),
\[
\begin{aligned}
\mathbb P^*(\boldsymbol\xi_1\in B_1,\ldots,\boldsymbol\xi_n\in B_n)
=
\frac{\prod_{i=1}^n\mathbb P(\boldsymbol\xi_i\in B_i,\,
                              |\boldsymbol\xi_i|\le b)}
     {\prod_{i=1}^n\mathbb P(|\boldsymbol\xi_i|\le b)}
=\prod_{i=1}^n\mathbb P(\boldsymbol\xi_i\in B_i
                         \mid |\boldsymbol\xi_i|\le b).
\end{aligned}
\]
Thus, under \(\mathbb P^*\), the observations remain independent,
and each \(\boldsymbol\xi_i\) has its original distribution
conditioned on \(|\boldsymbol\xi_i|\le b\).
We first estimate the
three terms under this law and will subsequently control the effect
of conditioning on the tail probability when returning to
\(\mathbb P\). Taking \(\mathbb E^*\) in
\eqref{eq:sec4-stein-identity-sample} gives
\begin{equation}\label{eq:sec4-stein-expectation-split}
\mathbb E^*H(R_n)-\mathbb EH(\chi_d)=I_1+I_2+I_3,
\end{equation}
where
\begin{align*}
I_1&=\mathbb E^*\left[
\operatorname{tr}\{V_n\partial_sF(S_n,V_n)\}-S_n^\top F(S_n,V_n)\right],\\
I_2&=\mathbb E^*\left[H(R_n)
-H\bigl((S_n^\top K(V_n)^{-1}S_n)^{1/2}\bigr)\right],\\
I_3&=\mathbb E^*\left[
\operatorname{tr}\{(K(V_n)-V_n)\partial_sF(S_n,V_n)\}\right].
\end{align*}
We begin with \(I_1\). Using the definitions of \(S_n\) and \(V_n\)
in \eqref{def:Sn-Vn}, together with the linearity and cyclic property
of the trace, we obtain
\[
\begin{aligned}
\operatorname{tr}\{V_n\partial_sF(S_n,V_n)\}
&=\sum_i\operatorname{tr}\{
\boldsymbol\xi_i\boldsymbol\xi_i^\top\partial_sF(S_n,V_n)
\}
=\sum_i\boldsymbol\xi_i^\top
\partial_sF(S_n,V_n)\boldsymbol\xi_i,\\
S_n^\top F(S_n,V_n)
&=\sum_i\boldsymbol\xi_i^\top F(S_n,V_n).
\end{aligned}
\]
Substituting these identities into the definition of \(I_1\) and
using the linearity of conditional expectation gives
\[
\begin{aligned}
I_1&=\sum_i\mathbb E^*\left[
\boldsymbol\xi_i^\top\partial_sF(S_n,V_n)\boldsymbol\xi_i
-\boldsymbol\xi_i^\top F(S_n,V_n)\right].
\end{aligned}
\]
In the second term of \(I_1\), \(F(S_n,V_n)\) depends on
\(\boldsymbol\xi_i\), since both \(S_n\) and \(V_n\) involve
this observation. To use independence,
we split this function value into a part depending only on the other
observations and the change caused by adding \(\boldsymbol\xi_i\).
Define
\[
 S_{-i}=\sum_{j\ne i}\boldsymbol\xi_j,\qquad
 V_{-i}=\sum_{j\ne i}\boldsymbol\xi_j\boldsymbol\xi_j^\top.
\]
Since the observations are independent under \(\mathbb P^*\),
\((S_{-i},V_{-i})\) is independent of \(\boldsymbol\xi_i\).
The conditional mean \(\mathbb E^*\boldsymbol\xi_i\) need not be
zero, so we retain it in the identity
\[
\mathbb E^*[\boldsymbol\xi_i^\top F(S_{-i},V_{-i})]
=(\mathbb E^*\boldsymbol\xi_i)^\top
\mathbb E^*F(S_{-i},V_{-i}).
\]
Adding and subtracting \(F(S_{-i},V_{-i})\) in the second term
therefore gives
\[
\begin{aligned}
\mathbb E^*[\boldsymbol\xi_i^\top F(S_n,V_n)]
&=(\mathbb E^*\boldsymbol\xi_i)^\top
  \mathbb E^*F(S_{-i},V_{-i})
+\mathbb E^*\left[
\boldsymbol\xi_i^\top
\{F(S_n,V_n)-F(S_{-i},V_{-i})\}\right].
\end{aligned}
\]
Substitution into \(I_1\) yields
\[
\begin{aligned}
I_1
&=-\sum_i(\mathbb E^*\boldsymbol\xi_i)^\top
       \mathbb E^*F(S_{-i},V_{-i})\\
&\quad+\sum_i\mathbb E^*\left[
\boldsymbol\xi_i^\top\partial_sF(S_n,V_n)\boldsymbol\xi_i
-\boldsymbol\xi_i^\top
\{F(S_n,V_n)-F(S_{-i},V_{-i})\}\right].
\end{aligned}
\]
To expand \(F(S_n,V_n)-F(S_{-i},V_{-i})\), note that
\(S_n-S_{-i}=\boldsymbol\xi_i\) and
\(V_n-V_{-i}=\boldsymbol\xi_i\boldsymbol\xi_i^\top\).
Integrating the derivatives of \(F\) along the segment from
\((S_{-i},V_{-i})\) to \((S_n,V_n)\) gives
\begin{equation}\label{eq:sec4-linear-expansion}
\begin{aligned}
F(S_n,V_n)-F(S_{-i},V_{-i})
&=\int_0^1\partial_sF\left(
S_{-i}+t\boldsymbol\xi_i,
V_{-i}+t\boldsymbol\xi_i\boldsymbol\xi_i^\top\right)
\boldsymbol\xi_i\,dt\\
&\quad+\int_0^1\partial_vF\left(
S_{-i}+t\boldsymbol\xi_i,
V_{-i}+t\boldsymbol\xi_i\boldsymbol\xi_i^\top\right)
[\boldsymbol\xi_i\boldsymbol\xi_i^\top]\,dt.
\end{aligned}
\end{equation}
Substituting this expansion into \(I_1\) and collecting terms,
we obtain
\[
\begin{aligned}
I_1
&=-\sum_i(\mathbb E^*\boldsymbol\xi_i)^\top
       \mathbb E^*F(S_{-i},V_{-i})\\
&\quad-\sum_i\mathbb E^*\int_0^1
\boldsymbol\xi_i^\top\partial_vF\left(
S_{-i}+t\boldsymbol\xi_i,
V_{-i}+t\boldsymbol\xi_i\boldsymbol\xi_i^\top\right)
[\boldsymbol\xi_i\boldsymbol\xi_i^\top]\,dt\\
&\quad+\sum_i\mathbb E^*\int_0^1
\boldsymbol\xi_i^\top\biggl\{\partial_sF(S_n,V_n)
-\partial_sF\left(
S_{-i}+t\boldsymbol\xi_i,
V_{-i}+t\boldsymbol\xi_i\boldsymbol\xi_i^\top\right)
\biggr\}\boldsymbol\xi_i\,dt.
\end{aligned}
\]
Taking absolute values, using
\(\|\boldsymbol\xi_i\boldsymbol\xi_i^\top\|_{\mathrm F}
=|\boldsymbol\xi_i|^2\) and
\(|u^\top A u|\le |u|^2\|A\|_{\mathrm F}\),
and applying the definition of \(\|\partial_vF\|\) in
\eqref{eq:sec4-global-matrix-derivative}, we obtain
\begin{equation}\label{eq:sec4-three-terms}
|I_1|\le I_{11}+I_{12}+I_{13},
\end{equation}
where
\begin{align*}
I_{11}&=\sum_{i=1}^n|\mathbb E^*\boldsymbol\xi_i|\,
\mathbb E^*|F(S_{-i},V_{-i})|,\\
I_{12}&=\sum_{i=1}^n\mathbb E^*\biggl[|\boldsymbol\xi_i|^3
\int_0^1\left\|\partial_vF\left(
S_{-i}+t\boldsymbol\xi_i,
V_{-i}+t\boldsymbol\xi_i\boldsymbol\xi_i^\top\right)\right\|\,dt\biggr],\\
I_{13}&=\sum_{i=1}^n\mathbb E^*\biggl[|\boldsymbol\xi_i|^2
\int_0^1\left\|\partial_sF(S_n,V_n)
-\partial_sF\left(S_{-i}+t\boldsymbol\xi_i,
V_{-i}+t\boldsymbol\xi_i\boldsymbol\xi_i^\top\right)\right\|_{\mathrm F}\,dt\biggr].
\end{align*}
Combining \eqref{eq:sec4-stein-expectation-split} with
\eqref{eq:sec4-three-terms}, we obtain
\[
|\mathbb E^*H(R_n)-\mathbb EH(\chi_d)|
\le I_{11}+I_{12}+I_{13}+|I_2|+|I_3|.
\]
We therefore need to bound the five terms on the right-hand side.
We first bound \(|I_2|+|I_3|\). The integrands defining \(I_2\)
and \(I_3\) vanish when \(\|V_n-I_d\|_{\mathrm F}\le1/4\), so define
\begin{equation}\label{eq:sec4-good-matrix}
\mathcal G=\{\|V_n-I_d\|_{\mathrm F}\le1/4\}.
\end{equation}
On \(\mathcal G\), we have \(K(V_n)=V_n\) and
\(R_n=(S_n^\top V_n^{-1}S_n)^{1/2}\), so the integrands in
\(I_2\) and \(I_3\) are zero. For \(I_2\), the bound
\(0\le H\le1\) therefore gives
\(|I_2|\le\mathbb P^*(\mathcal G^c)\).
For \(I_3\), the Frobenius norm inequality for the trace bounds
the integrand in absolute value by
\(\|K(V_n)-V_n\|_{\mathrm F}
\|\partial_sF(S_n,V_n)\|_{\mathrm F}\).
Combining these estimates with \eqref{eq:sec4-truncated-field}
and the spectral bounds for \(K(V_n)\), we obtain
\begin{equation}\label{eq:sec4-matrix-corrections}
\begin{aligned}
|I_2|+|I_3|
&\le\mathbb P^*(\mathcal G^c)
+\mathbb E^*\left[
\|K(V_n)-V_n\|_{\mathrm F}\,
\|\partial_sF(S_n,V_n)\|_{\mathrm F};\mathcal G^c\right]\\
&\le\mathbb P^*(\mathcal G^c)
+C_d\mathbb E^*\left[
\|K(V_n)\|_{\mathrm F}+\|V_n\|_{\mathrm F};\mathcal G^c\right]\\
&\le C_d\mathbb E^*[(1+\|V_n\|_{\mathrm F})
\mathbf1_{\mathcal G^c}].
\end{aligned}
\end{equation}
Substituting \eqref{eq:sec4-matrix-corrections} and \eqref{eq:sec4-three-terms} into
\eqref{eq:sec4-stein-expectation-split}, we obtain
\begin{equation}\label{eq:sec4-stein-decomposition}
\begin{aligned}
|\mathbb E^*H(R_n)-\mathbb EH(\chi_d)|
&\le I_{11}+I_{12}+I_{13}
+C_d\mathbb E^*[(1+\|V_n\|_{\mathrm F})
\mathbf1_{\mathcal G^c}].
\end{aligned}
\end{equation}
We first estimate the right-hand side of
\eqref{eq:sec4-stein-decomposition} for \(0\le x\le4\).
For \(u\ge0\) and \(0<\eta\le1\), take \(H=H_{u,\eta}\)
in the definitions of \(F\) and \(I_1,I_2,I_3\) in
\eqref{eq:sec4-extension} and
\eqref{eq:sec4-stein-expectation-split}.
Throughout this part, \(b=\varepsilon/4\), and we work under
the condition in
\eqref{eq:sec4-conditioning-bound}.
We first complete the bound for \(|I_2|+|I_3|\).
On \(\mathcal G^c\),
\(1+\|V_n\|_{\mathrm F}\le C_d\|V_n-I_d\|_{\mathrm F}^2\),
so \eqref{eq:sec4-matrix-corrections} gives
\[
|I_2|+|I_3|\le C_d\mathbb E^*\|V_n-I_d\|_{\mathrm F}^2.
\]
To bound this expectation, write $V_n-I_d=(V_n-\mathbb E^*V_n)+(\mathbb E^*V_n-I_d)$.
Using independence under \(\mathbb P^*\), we obtain
\[
\begin{aligned}
\mathbb E^*\|V_n-I_d\|_{\mathrm F}^2
&=\mathbb E^*\|V_n-\mathbb E^*V_n\|_{\mathrm F}^2
  +\|\mathbb E^*V_n-I_d\|_{\mathrm F}^2
\le\sum_i\mathbb E^*|\boldsymbol\xi_i|^4
  +\|\mathbb E^*V_n-I_d\|_{\mathrm F}^2.
\end{aligned}
\]
Thus we need bounds for
\(\sum_i\mathbb E^*|\boldsymbol\xi_i|^4\) and
\(\|\mathbb E^*V_n-I_d\|_{\mathrm F}\).
The following lemma provides the required bounds.

\begin{lemma}\label{lem:sec4-conditional-moments}
Fix $2<p\le3$. Let $\boldsymbol\xi_1,\ldots,\boldsymbol\xi_n$
be independent centered random vectors in $\mathbb R^d$ such that
$\sum_{i=1}^n\mathbb E (\boldsymbol\xi_i\boldsymbol\xi_i^\top)=I_d$.
Use the cutoff $b$ and the product conditional law in
\eqref{eq:sec4-conditioning}, and suppose that
$b^{-p}\sum_{i=1}^n\mathbb E|\boldsymbol\xi_i|^p\le\frac14$.
Then $\mathbb P(\mathcal N)\ge3/4$, the observations remain
independent under $\mathbb P^*$, and
\begin{equation}\label{eq:sec4-conditional-moments}
\begin{aligned}
\sum_{i=1}^n|\mathbb E^*\boldsymbol\xi_i|
 &\le2b^{1-p}\sum_{i=1}^n\mathbb E|\boldsymbol\xi_i|^p,
\qquad
\|\mathbb E^*V_n-I_d\|_{\mathrm F}
 \le3b^{2-p}\sum_{i=1}^n\mathbb E|\boldsymbol\xi_i|^p,\\
\sum_{i=1}^n\mathbb E^*|\boldsymbol\xi_i|^j
 &\le2b^{j-p}\sum_{i=1}^n\mathbb E|\boldsymbol\xi_i|^p, \quad j=3,4,
\qquad
\sum_{i=1}^n\mathbb E^*|\boldsymbol\xi_i|^2
 \le4d/3.
\end{aligned}
\end{equation}
\end{lemma}

Returning to \(|I_2|+|I_3|\), we apply
\eqref{eq:sec4-conditional-moments} with the fixed
\(b=\varepsilon/4\) to the two quantities above. This gives
\[
\sum_i\mathbb E^*|\boldsymbol\xi_i|^4
+\|\mathbb E^*V_n-I_d\|_{\mathrm F}
\le C_d\sum_i\mathbb E|\boldsymbol\xi_i|^p.
\]
The condition in \eqref{eq:sec4-conditioning-bound} also implies
\(\sum_i\mathbb E|\boldsymbol\xi_i|^p\le b^p/4\le1\),
so \(\|\mathbb E^*V_n-I_d\|_{\mathrm F}^2
\le C_d\sum_i\mathbb E|\boldsymbol\xi_i|^p\) as well.
Substituting these bounds into the preceding estimate yields
\begin{equation}\label{eq:sec4-small-matrix}
\begin{aligned}
|I_2|+|I_3|
&\le C_d\mathbb E^*[(1+\|V_n\|_{\mathrm F})\mathbf1_{\mathcal G^c}]
\le C_d\mathbb E^*\|V_n-I_d\|_{\mathrm F}^2
\le C_d\sum_i\mathbb E|\boldsymbol\xi_i|^p.
\end{aligned}
\end{equation}

We next estimate \(I_{11}\). Since \(V_{-i}\succeq0\), applying
\eqref{eq:sec4-truncated-field} at \(s=S_{-i}\) and \(v=V_{-i}\)
gives \(|F(S_{-i},V_{-i})|\le C_d\).
Consequently, \(\mathbb E^*|F(S_{-i},V_{-i})|\le C_d\),
and the definition of \(I_{11}\) yields
\[
I_{11}
=\sum_i|\mathbb E^*\boldsymbol\xi_i|\,
       \mathbb E^*|F(S_{-i},V_{-i})|
\le C_d\sum_i|\mathbb E^*\boldsymbol\xi_i|.
\]
The conditional-mean estimate in \eqref{eq:sec4-conditional-moments},
with the fixed \(b\), therefore yields
\[
I_{11}\le C_d\sum_i\mathbb E|\boldsymbol\xi_i|^p.
\]

For \(I_{12}\), use \(|S_{-i}+t\boldsymbol\xi_i|\le|S_{-i}|+|\boldsymbol\xi_i|\)
for \(0\le t\le1\) in the matrix-derivative bound
\eqref{eq:sec4-global-matrix-derivative}. This gives
\begin{equation}\label{eq:sec4-small-first-two-preliminary}
I_{12}\le C_d\sum_i\mathbb E^*\left[
|\boldsymbol\xi_i|^3(1+|S_{-i}|+|\boldsymbol\xi_i|)\right].
\end{equation}
To bound the right-hand side, we need the conditional third and
fourth moments of \(\boldsymbol\xi_i\), as well as a bound on
\(\mathbb E^*|S_{-i}|\).
Applying Lemma~\ref{lem:sec4-conditional-moments}
with \(b=\varepsilon/4\) gives
\begingroup
the following bounds with constants depending only on \(d\), since
\(b^{1-p}\le b^{-2}\), \(b^{2-p}\le b^{-1}\), and
\(b^{j-p}\le1\) for \(j=3,4\):
\endgroup
\begin{equation}\label{eq:sec4-small-conditional-moments}
\begin{aligned}
\sum_i|\mathbb E^*\boldsymbol\xi_i|
+\|\mathbb E^*V_n-I_d\|_{\mathrm F}
&\le C_d\sum_i\mathbb E|\boldsymbol\xi_i|^p,\\
\sum_i\bigl(\mathbb E^*|\boldsymbol\xi_i|^3
+\mathbb E^*|\boldsymbol\xi_i|^4\bigr)
&\le C_d\sum_i\mathbb E|\boldsymbol\xi_i|^p.
\end{aligned}
\end{equation}
To control \(\mathbb E^*|S_{-i}|\), we first bound its second moment.
Independence under \(\mathbb P^*\), the second-moment estimate in
\eqref{eq:sec4-conditional-moments}, and the conditional-mean bound
in \eqref{eq:sec4-small-conditional-moments} give
\begin{equation}\label{eq:sec4-small-deleted-second-moment}
\sup_i\mathbb E^*|S_{-i}|^2
\le\sum_j\mathbb E^*|\boldsymbol\xi_j|^2
+\left(\sum_j|\mathbb E^*\boldsymbol\xi_j|\right)^2\le C_d.
\end{equation}
Since \(\boldsymbol\xi_i\) and \(S_{-i}\) are independent under
\(\mathbb P^*\), the expectation of
\(|\boldsymbol\xi_i|^3|S_{-i}|\) factors into two expectations.
Using \(\mathbb E^*|S_{-i}|\le(\mathbb E^*|S_{-i}|^2)^{1/2}\),
we can now substitute \eqref{eq:sec4-small-deleted-second-moment}
and \eqref{eq:sec4-small-conditional-moments} into
\eqref{eq:sec4-small-first-two-preliminary} to obtain
\[
I_{12}\le C_d\sum_i\bigl(\mathbb E^*|\boldsymbol\xi_i|^3
+\mathbb E^*|\boldsymbol\xi_i|^4\bigr)
\le C_d\sum_i\mathbb E|\boldsymbol\xi_i|^p.
\]
Combining this with the bound already proved for \(I_{11}\), we obtain
\begin{equation}\label{eq:sec4-small-first-two}
I_{11}+I_{12}\le C_d\sum_i\mathbb E|\boldsymbol\xi_i|^p.
\end{equation}

It remains to estimate \(I_{13}\), defined after
\eqref{eq:sec4-three-terms}.
Bounding the two Jacobians in its integrand separately by
\eqref{eq:sec4-truncated-field} would give only
\[
I_{13}\le C_d\sum_i\mathbb E^*|\boldsymbol\xi_i|^2\le C_d.
\]
\begingroup
To obtain a bound that tends to zero with the moment sum, we
separate the event that \(V_{-i}\) is close to \(I_d\) from its
complement. Define
\endgroup
\begin{equation}\label{eq:sec4-good-deleted-matrix}
\mathcal G_i=\{\|V_{-i}-I_d\|_{\mathrm F}\le1/8\}.
\end{equation}
\begingroup
On \(\mathcal G_i^c\), we will use the global Jacobian bound
and control \(\mathbb P^*(\mathcal G_i^c)\).
On \(\mathcal G_i\), we will estimate the change in the Jacobian
along the segment. Split \(I_{13}\) accordingly as
\endgroup
\begin{equation}\label{eq:sec4-small-third-split}
\begin{aligned}
I_{13}
&=\underbrace{\sum_i\mathbb E^*\biggl[
 |\boldsymbol\xi_i|^2\int_0^1
 \left\|\partial_sF(S_n,V_n)
 -\partial_sF\left(S_{-i}+t\boldsymbol\xi_i,
 V_{-i}+t\boldsymbol\xi_i\boldsymbol\xi_i^\top\right)
 \right\|_{\mathrm F}\,dt;\mathcal G_i^c\biggr]}_{I_{131}}\\
&\quad+\underbrace{\sum_i\mathbb E^*\biggl[
 |\boldsymbol\xi_i|^2\int_0^1
 \left\|\partial_sF(S_n,V_n)
 -\partial_sF\left(S_{-i}+t\boldsymbol\xi_i,
 V_{-i}+t\boldsymbol\xi_i\boldsymbol\xi_i^\top\right)
 \right\|_{\mathrm F}\,dt;\mathcal G_i\biggr]}_{I_{132}}.
\end{aligned}
\end{equation}
We first estimate \(I_{131}\).
By \eqref{eq:sec4-truncated-field}, the norm of the Jacobian
difference is at most \(2C_d\).
Moreover, \(\mathcal G_i\) depends only on the observations with
indices different from \(i\), so it is independent of
\(\boldsymbol\xi_i\) under \(\mathbb P^*\). Hence
\begin{equation}\label{eq:sec4-small-bad-matrix-target}
\begin{aligned}
I_{131}
&\le C_d\sum_i\mathbb E^*\left[
 |\boldsymbol\xi_i|^2\mathbf1_{\mathcal G_i^c}\right]
=C_d\sum_i\mathbb E^*|\boldsymbol\xi_i|^2
                    \mathbb P^*(\mathcal G_i^c).
\end{aligned}
\end{equation}
Since \(\sum_i\mathbb E^*|\boldsymbol\xi_i|^2\le4d/3\), it
remains to bound \(\mathbb P^*(\mathcal G_i^c)\), uniformly in \(i\).
Markov's inequality reduces this to bounding
\(\mathbb E^*\|V_{-i}-I_d\|_{\mathrm F}^2\), which decomposes as
\[
\mathbb E^*\|V_{-i}-I_d\|_{\mathrm F}^2
=\mathbb E^*\|V_{-i}-\mathbb E^*V_{-i}\|_{\mathrm F}^2
 +\|\mathbb E^*V_{-i}-I_d\|_{\mathrm F}^2.
\]
For the second term, use
\(\mathbb E^*V_{-i}=\mathbb E^*V_n
-\mathbb E^*(\boldsymbol\xi_i\boldsymbol\xi_i^\top)\).
The conditional expectation formula and H\"older's inequality give
\[
\mathbb E^*|\boldsymbol\xi_i|^2
\le\frac43\mathbb E|\boldsymbol\xi_i|^2
\le\frac43(\mathbb E|\boldsymbol\xi_i|^p)^{2/p},
\]
where the factor \(4/3\) follows from
\(\mathbb P(|\boldsymbol\xi_i|\le b)\ge3/4\).
Together with \eqref{eq:sec4-small-conditional-moments}, this gives
\begin{equation}\label{eq:sec4-small-deleted-matrix-bias}
\begin{aligned}
\|\mathbb E^*V_{-i}-I_d\|_{\mathrm F}
&\le\|\mathbb E^*V_n-I_d\|_{\mathrm F}
 +\mathbb E^*|\boldsymbol\xi_i|^2
\le C_d\left\{\sum_j\mathbb E|\boldsymbol\xi_j|^p
 +(\mathbb E|\boldsymbol\xi_i|^p)^{2/p}\right\}.
\end{aligned}
\end{equation}
For the first term in the decomposition, independence under
\(\mathbb P^*\) and the fourth-moment bound in
\eqref{eq:sec4-small-conditional-moments} give
\begin{equation}\label{eq:sec4-small-deleted-matrix-variance}
\mathbb E^*\|V_{-i}-\mathbb E^*V_{-i}\|_{\mathrm F}^2
\le\sum_{j\ne i}\mathbb E^*|\boldsymbol\xi_j|^4
\le C_d\sum_j\mathbb E|\boldsymbol\xi_j|^p.
\end{equation}
Since \(\sum_j\mathbb E|\boldsymbol\xi_j|^p\le1\) and \(4/p>1\),
both the square of this sum and
\((\mathbb E|\boldsymbol\xi_i|^p)^{4/p}\) are bounded by the sum.
Thus the squared bias in
\eqref{eq:sec4-small-deleted-matrix-bias},
together with \eqref{eq:sec4-small-deleted-matrix-variance}, bounds
\(\mathbb E^*\|V_{-i}-I_d\|_{\mathrm F}^2\) by
\(C_d\sum_j\mathbb E|\boldsymbol\xi_j|^p\), uniformly in \(i\).
Markov's inequality now gives
\begin{equation}\label{eq:sec4-small-deleted-bad-probability}
\begin{aligned}
\sup_i\mathbb P^*(\mathcal G_i^c)
&=\sup_i\mathbb P^*(\|V_{-i}-I_d\|_{\mathrm F}>1/8)
\le64\sup_i\mathbb E^*\|V_{-i}-I_d\|_{\mathrm F}^2
\le C_d\sum_j\mathbb E|\boldsymbol\xi_j|^p.
\end{aligned}
\end{equation}
Substituting \eqref{eq:sec4-small-deleted-bad-probability} and
\eqref{eq:sec4-conditional-moments} into
\eqref{eq:sec4-small-bad-matrix-target} gives
\begin{equation}\label{eq:sec4-small-bad-matrix-contribution}
I_{131}\le C_d\sum_i\mathbb E^*|\boldsymbol\xi_i|^2
 \mathbb P^*(\mathcal G_i^c)
\le C_d\sum_i\mathbb E|\boldsymbol\xi_i|^p.
\end{equation}

We now estimate \(I_{132}\). In its definition in
\eqref{eq:sec4-small-third-split}, the integral over \([0,1]\)
is bounded by the supremum of its integrand. Hence
\[
\begin{aligned}
I_{132}\le\sum_i\mathbb E^*\biggl[
 |\boldsymbol\xi_i|^2\sup_{0\le t\le1}
 \Bigl\|\partial_sF(S_n,V_n)
 -\partial_sF\left(S_{-i}+t\boldsymbol\xi_i,
 V_{-i}+t\boldsymbol\xi_i\boldsymbol\xi_i^\top\right)
 \Bigr\|_{\mathrm F};\mathcal G_i\biggr].
\end{aligned}
\]
We therefore need a bound for this supremum on \(\mathcal G_i\).
The following lemma gives such a bound in terms of
\(|\boldsymbol\xi_i|\) and the radius formed from the remaining
observations, defined by
\[
 R_{-i}=(S_{-i}^\top V_{-i}^+S_{-i})^{1/2},
\]
where \(V_{-i}^+\) is the Moore--Penrose inverse of \(V_{-i}\).
The proof of the lemma is given in
Section~\ref{sec:sec4-lemma-proofs}.

\begin{lemma}\label{lem:sec4-small-increments}
Let $F$ be defined by \eqref{eq:sec4-extension} using
$H_{u,\eta}$ from \eqref{eq:sec4-cutoff}, and use the cutoff
$b$ from \eqref{eq:sec4-conditioning}. On $\mathcal G_i$, whenever
$|\boldsymbol\xi_i|\le b$, we have
\begin{equation}\label{eq:sec4-small-path}
\begin{aligned}
&\sup_{0\le t\le1}\left\|\partial_sF(S_n,V_n)
 -\partial_sF\left(S_{-i}+t\boldsymbol\xi_i,
                  V_{-i}+t\boldsymbol\xi_i\boldsymbol\xi_i^\top\right)
  \right\|_{\mathrm F}\\
&\quad\le C_d|\boldsymbol\xi_i|(1+R_{-i}^2)
 +C_d\min\left\{1,\frac{C_d|\boldsymbol\xi_i|}{R_{-i}+\eta}\right\}\\
&\qquad+C_d\min\left\{1,\frac{C_d|\boldsymbol\xi_i|(1+u)}{\eta}\right\}
 \mathbf1_{\{|R_{-i}-u|\le C_d|\boldsymbol\xi_i|(1+u)+\eta\}}.
\end{aligned}
\end{equation}
The constants are uniform in $u$, $\eta$, and $i$.
\end{lemma}

\begingroup
Under \(\mathbb P^*\), the truncation gives
\(|\boldsymbol\xi_i|\le b=\varepsilon/4\) almost surely.
Applying \eqref{eq:sec4-small-path} to the supremum in the bound
for \(I_{132}\), we obtain
\begin{equation}\label{eq:sec4-small-third-three-contributions}
\begin{aligned}
I_{132}
&\le C_d\sum_i\mathbb E^*\left[
 |\boldsymbol\xi_i|^3(1+R_{-i}^2)\mathbf1_{\mathcal G_i}\right]\\
&\quad+C_d\sum_i\mathbb E^*\left[
 |\boldsymbol\xi_i|^2
 \min\left\{1,\frac{C_d|\boldsymbol\xi_i|}{R_{-i}+\eta}\right\}
 \mathbf1_{\mathcal G_i}\right]\\
&\quad+C_d\sum_i\mathbb E^*\left[
 |\boldsymbol\xi_i|^2
 \min\left\{1,\frac{C_d|\boldsymbol\xi_i|(1+u)}{\eta}\right\}
 \mathbf1_{\{|R_{-i}-u|\le C_d|\boldsymbol\xi_i|(1+u)+\eta\}}
 \mathbf1_{\mathcal G_i}\right].
\end{aligned}
\end{equation}
We first estimate the first term on the right-hand side.
On \(\mathcal G_i\), the smallest eigenvalue of \(V_{-i}\) is
at least \(7/8\), so
\[
R_{-i}^2=S_{-i}^\top V_{-i}^{-1}S_{-i}
\le\frac87|S_{-i}|^2,
\qquad \text{on }\mathcal G_i.
\]
Substituting this bound, dropping the indicator, and using the
independence of \(\boldsymbol\xi_i\) and \(S_{-i}\) under
\(\mathbb P^*\), we obtain
\begin{equation}\label{eq:sec4-small-easy-increments}
\begin{aligned}
&\sum_i\mathbb E^*\left[
 |\boldsymbol\xi_i|^3(1+R_{-i}^2)\mathbf1_{\mathcal G_i}\right]
\le C_d\sum_i\mathbb E^*|\boldsymbol\xi_i|^3
 (1+\mathbb E^*|S_{-i}|^2)
\le C_d\sum_i\mathbb E|\boldsymbol\xi_i|^p,
\end{aligned}
\end{equation}
where the last inequality uses the conditional third-moment
bound in \eqref{eq:sec4-small-conditional-moments} and the
second-moment bound \eqref{eq:sec4-small-deleted-second-moment}.
This controls the first term of
\eqref{eq:sec4-small-third-three-contributions}.

For the second term of
\eqref{eq:sec4-small-third-three-contributions}, split the
expectation according to whether \(R_{-i}>1\) or \(R_{-i}\le1\).
On \(R_{-i}>1\), we have
\[
\min\left\{1,\frac{C_d|\boldsymbol\xi_i|}{R_{-i}+\eta}\right\}
\le C_d|\boldsymbol\xi_i|,
\qquad R_{-i}>1.
\]
Consequently, the conditional third-moment bound gives
\begin{equation}\label{eq:sec4-small-second-large-radius}
\begin{aligned}
&\sum_i\mathbb E^*\left[
 |\boldsymbol\xi_i|^2
 \min\left\{1,\frac{C_d|\boldsymbol\xi_i|}{R_{-i}+\eta}\right\}
 \mathbf1_{\mathcal G_i}\mathbf1_{\{R_{-i}>1\}}\right]
\le C_d\sum_i\mathbb E^*|\boldsymbol\xi_i|^3
\le C_d\sum_i\mathbb E|\boldsymbol\xi_i|^p.
\end{aligned}
\end{equation}
The part with \(R_{-i}\le1\) remains to be estimated, as does
the third term of \eqref{eq:sec4-small-third-three-contributions}.
Substituting \eqref{eq:sec4-small-easy-increments} and
\eqref{eq:sec4-small-second-large-radius} into
\eqref{eq:sec4-small-third-three-contributions}, and dropping
\(\mathbf1_{\mathcal G_i}\) from these two remaining nonnegative
expectations, gives
\endgroup
\begin{equation}\label{eq:sec4-small-third-remainder}
\begin{aligned}
I_{132}
&\le C_d\sum_i\mathbb E|\boldsymbol\xi_i|^p\\
&\quad+C_d\sum_i\mathbb E^*\left[
 |\boldsymbol\xi_i|^2
 \min\left\{1,\frac{C_d|\boldsymbol\xi_i|}{R_{-i}+\eta}\right\}
 \mathbf1_{\{R_{-i}\le1\}}\right]\\
&\quad+C_d\sum_i\mathbb E^*\left[
 |\boldsymbol\xi_i|^2\min\left\{1,
 \frac{C_d|\boldsymbol\xi_i|(1+u)}{\eta}\right\}
 \mathbf1_{\{|R_{-i}-u|\le C_d|\boldsymbol\xi_i|(1+u)+\eta\}}\right].
\end{aligned}
\end{equation}
\begingroup
The two remaining expectations in
\eqref{eq:sec4-small-third-remainder} depend on the distribution
of \(R_{-i}\) near zero and in a short interval around \(u\).
We estimate them by comparing this distribution with that of
\(\chi_d\). Since the approximation error is not yet known,
we retain it in the following quantity. For \(\lambda>0\), define
\endgroup
\begin{equation}\label{eq:sec4-small-error}
\mathcal B(\lambda)
=\sup\frac{|\mathbb P(R_n\le t)-\mathbb P(\chi_d\le t)|}{1+t},
\end{equation}
where the supremum is over all \(n\ge1\), \(t\ge0\), and centered
independent arrays with total covariance \(I_d\) and
\(\sum_i\mathbb E|\boldsymbol\xi_i|^p\le\lambda\).
In particular, \(0\le\mathcal B(\lambda)\le1\).

\begingroup
Here \(\mathcal B(\lambda)\) measures the distributional error
under the original law \(\mathbb P\), over all thresholds and
arrays with moment sum at most \(\lambda\).
To apply it to the deleted sample, we first normalize the
covariance under \(\mathbb P\), and then account for conditioning.
Choose \(0<\lambda_0<1\), depending only on \(d\), so that
\endgroup
\begingroup
\begin{equation}\label{eq:sec4-small-moment-threshold}
(4/\varepsilon)^3\lambda_0\le1/4,\qquad
(1-\lambda_0^{2/3})^{-3/2}\le2.
\end{equation}
For every \(2<p\le3\) and \(0<\lambda\le\lambda_0\),
\[
 (4/\varepsilon)^p\le(4/\varepsilon)^3,\qquad
 (1-\lambda^{2/p})^{-p/2}
 \le(1-\lambda^{2/3})^{-3/2}\le2.
\]
\endgroup
\begingroup
For \(\sum_i\mathbb E|\boldsymbol\xi_i|^p\le\lambda\le\lambda_0\),
the first condition ensures
\(b^{-p}\sum_i\mathbb E|\boldsymbol\xi_i|^p\le1/4\), as required
for the conditional moment estimates.
The second ensures that normalizing the covariance after deleting
one observation increases the moment sum to at most \(2\lambda\).
This gives a comparison using \(\mathcal B(2\lambda)\);
passing from \(\mathbb P\) to \(\mathbb P^*\) adds at most
\(C_d\lambda\) to this distributional comparison.
The following lemma carries out this comparison and bounds
the two remaining expectations in
\eqref{eq:sec4-small-third-remainder}.
\endgroup

\begin{lemma}\label{lem:sec4-deleted-expectations}
Under the assumptions of Lemma~\ref{lem:sec4-conditional-moments},
let \(d\ge2\), \(0\le x\le4\), and suppose that
\(\sum_i\mathbb E|\boldsymbol\xi_i|^p\le\lambda\le\lambda_0\).
For the parameters \(u,\eta\) in \eqref{eq:sec4-cutoff},
\begin{equation}\label{eq:sec4-small-inverse-contribution}
\begin{aligned}
&\sum_i\mathbb E^*\left[
 |\boldsymbol\xi_i|^2
 \min\left\{1,\frac{C_d|\boldsymbol\xi_i|}{R_{-i}+\eta}\right\}
 \mathbf1_{\{R_{-i}\le1\}}\right]
\le C_d\lambda
 +\frac{C_d\lambda}{\eta}\{\mathcal B(2\lambda)+\lambda\}.
\end{aligned}
\end{equation}
\begin{equation}\label{eq:sec4-small-shell-contribution}
\begin{aligned}
&\sum_i\mathbb E^*\left[
 |\boldsymbol\xi_i|^2\min\left\{1,
 \frac{C_d|\boldsymbol\xi_i|(1+u)}{\eta}\right\}
 \mathbf1_{\{|R_{-i}-u|\le C_d|\boldsymbol\xi_i|(1+u)+\eta\}}\right]\\
&\quad\le C_d\lambda(1+u)
 +\frac{C_d\lambda(1+u)}{\eta}
 \{(1+u)\mathcal B(2\lambda)+\lambda\}.
\end{aligned}
\end{equation}
\end{lemma}

\begingroup
The bounds \eqref{eq:sec4-small-inverse-contribution} and
\eqref{eq:sec4-small-shell-contribution} control, respectively,
the expectation on \(R_{-i}\le1\) and the expectation over the
short interval around \(u\) in
\eqref{eq:sec4-small-third-remainder}.
Substituting them there bounds \(I_{132}\).
Adding the bound for \(I_{131}\) from
\eqref{eq:sec4-small-bad-matrix-contribution}, we obtain
\endgroup
\begin{equation}\label{eq:sec4-small-increments}
\begin{aligned}
I_{13}&=I_{131}+I_{132}
 \le C_d\sum_i\mathbb E|\boldsymbol\xi_i|^p+I_{132}\\
&\le C_d\lambda(1+u)
 +\frac{C_d\lambda(1+u)^2}{\eta}\mathcal B(2\lambda)
 +\frac{C_d\lambda^2(1+u)}{\eta}.
\end{aligned}
\end{equation}

\begingroup
The estimate for \(I_{13}\) still contains
\(\mathcal B(2\lambda)\), which will enter the recursive
inequality below. We now return to the conditional expectation
in \eqref{eq:sec4-stein-decomposition}.
For \(0<\lambda\le\lambda_0\) and
\(\sum_i\mathbb E|\boldsymbol\xi_i|^p\le\lambda\),
substituting \eqref{eq:sec4-small-matrix},
\eqref{eq:sec4-small-first-two}, and
\eqref{eq:sec4-small-increments} into
\eqref{eq:sec4-stein-decomposition} gives
\endgroup
\begin{equation}\label{eq:sec4-small-smooth-recursion}
\begin{aligned}
&\left|\mathbb E^*H_{u,\eta}(R_n)
       -\mathbb EH_{u,\eta}(\chi_d)\right|\\
&\quad\le C_d\lambda(1+u)
 +C_d\frac{\lambda(1+u)^2}{\eta}\mathcal B(2\lambda)
 +C_d\frac{\lambda^2(1+u)}{\eta},
\qquad u\ge0,\quad0<\eta\le1.
\end{aligned}
\end{equation}

\begingroup
We next use \eqref{eq:sec4-small-smooth-recursion} to estimate
\(\left|\mathbb P(R_n\le t)-\mathbb P(\chi_d\le t)\right|\).
We must return from \(\mathbb P^*\) to \(\mathbb P\) and replace
the smooth functions by the indicator of \(\{y\le t\}\).
For \(t\ge0\), decompose the difference as
\begin{equation}\label{eq:sec4-small-cdf-decomposition}
\begin{aligned}
&\left|\mathbb P(R_n\le t)-\mathbb P(\chi_d\le t)\right|
\le\left|\mathbb P(R_n\le t)-\mathbb P^*(R_n\le t)\right|
 +\left|\mathbb P^*(R_n\le t)-\mathbb P(\chi_d\le t)\right|.
\end{aligned}
\end{equation}
For the first term, the definition
\(\mathbb P^*(\cdot)=\mathbb P(\cdot\mid\mathcal N)\) gives
\[
\begin{aligned}
\mathbb P(R_n\le t)-\mathbb P^*(R_n\le t)
&=\mathbb P(R_n\le t,\mathcal N^c)
 -\mathbb P^*(R_n\le t)\mathbb P(\mathcal N^c).
\end{aligned}
\]
Both terms on the right lie between zero and
\(\mathbb P(\mathcal N^c)\). Hence, by
\eqref{eq:sec4-conditioning-bound} and the fixed
\(b=\varepsilon/4\),
\[
\begin{aligned}
\left|\mathbb P(R_n\le t)-\mathbb P^*(R_n\le t)\right|
&\le\mathbb P(\mathcal N^c)
\le b^{-p}\sum_i\mathbb E|\boldsymbol\xi_i|^p
\le C_d\lambda.
\end{aligned}
\]
This bounds the first term of
\eqref{eq:sec4-small-cdf-decomposition}.

For the second term, we use the smooth functions to bound the
indicator from above and below. If \(t\ge\eta\),
\eqref{eq:sec4-cutoff} gives
$H_{t-\eta,\eta}(y)\le\mathbf1_{\{y\le t\}}
\le H_{t,\eta}(y)$, $y\ge0$,
and therefore
\[
\mathbb E^*H_{t-\eta,\eta}(R_n)
\le\mathbb P^*(R_n\le t)
\le\mathbb E^*H_{t,\eta}(R_n).
\]
For \(\chi_d\), these two functions differ from the indicator
only in \([t-\eta,t+\eta]\). Comparing the two endpoints of the
last display with the corresponding expectations under the chi
law yields
\[
\begin{aligned}
&\left|\mathbb P^*(R_n\le t)-\mathbb P(\chi_d\le t)\right|
\le\max_{u\in\{t-\eta,t\}}
 \left|\mathbb E^*H_{u,\eta}(R_n)
       -\mathbb EH_{u,\eta}(\chi_d)\right|
 +\mathbb P(t-\eta\le\chi_d\le t+\eta).
\end{aligned}
\]
To bound the last probability, use the bounded density \(p_d\)
of \(\chi_d\):
\[
\mathbb P(t-\eta\le\chi_d\le t+\eta)
=\int_{\max\{0,t-\eta\}}^{t+\eta}p_d(y)\,dy
\le2\eta\sup_{y>0}p_d(y)\le C_d\eta.
\]
Thus the comparison with the smooth functions contributes at
most \(C_d\eta\) in addition to their expectation differences.

If \(0\le t<\eta\), we still have
\(\mathbf1_{\{y\le t\}}\le H_{t,\eta}(y)\) for the upper bound.
For the lower bound,
\[
\mathbb P(\chi_d\le t)-\mathbb P^*(R_n\le t)
\le\mathbb P(\chi_d\le t)\le C_d\eta.
\]
Consequently, in this case,
\[
\left|\mathbb P^*(R_n\le t)-\mathbb P(\chi_d\le t)\right|
\le\left|\mathbb E^*H_{t,\eta}(R_n)
          -\mathbb EH_{t,\eta}(\chi_d)\right|+C_d\eta.
\]
We can now apply \eqref{eq:sec4-small-smooth-recursion} with
\(u=t\) and, when \(t\ge\eta\), with \(u=t-\eta\).
In both cases \(1+u\le1+t\).
Substituting the resulting bound and the conditioning estimate
into \eqref{eq:sec4-small-cdf-decomposition}, and absorbing
\(C_d\lambda\) into \(C_d\lambda(1+t)\), gives
\[
\begin{aligned}
&\left|\mathbb P(R_n\le t)-\mathbb P(\chi_d\le t)\right|
\le C_d\left\{
\lambda(1+t)+\eta
+\frac{\lambda(1+t)^2}{\eta}\mathcal B(2\lambda)
+\frac{\lambda^2(1+t)}{\eta}\right\}.
\end{aligned}
\]
Dividing by \(1+t\), as in the definition
\eqref{eq:sec4-small-error} of \(\mathcal B(\lambda)\), we obtain
\endgroup

\begin{equation}\label{eq:sec4-small-cdf-recursive-bound}
\begin{aligned}
&\frac{\left|\mathbb P(R_n\le t)-\mathbb P(\chi_d\le t)\right|}{1+t}
\le C_d\left\{
\lambda+\frac{\eta}{1+t}
+\frac{\lambda(1+t)}{\eta}\mathcal B(2\lambda)
+\frac{\lambda^2}{\eta}\right\},
\qquad t\ge0,\quad0<\eta\le1.
\end{aligned}
\end{equation}

\begingroup
The right-hand side still contains \(\mathcal B(2\lambda)\).
To obtain a recursive bound for \(\mathcal B(\lambda)\), we
choose \(\eta\) so that the coefficient of
\(\mathcal B(2\lambda)\) is a small constant.
Taking \(\eta=m\lambda(1+t)\) makes
\(\lambda(1+t)/\eta=1/m\).
Fix an integer \(m\ge1\), depending only on \(d\), large enough for both recursions below,
with
\endgroup

\begin{equation}\label{eq:sec4-recursion-choice}
C_d/m\le 1/4,
\end{equation}

\begingroup
where \(C_d\) bounds their coefficients and is independent of \(m\).
If \(m\lambda(1+t)\le1\), this choice of \(\eta\) satisfies the
required restriction \(0<\eta\le1\). Substituting it into
\eqref{eq:sec4-small-cdf-recursive-bound} gives
\endgroup

\[
\begin{aligned}
&\frac{\left|\mathbb P(R_n\le t)-\mathbb P(\chi_d\le t)\right|}{1+t}
\le C_d\left\{
\lambda+m\lambda+\frac{\mathcal B(2\lambda)}m
+\frac{\lambda}{m(1+t)}\right\}
\le C_dm\lambda+\frac{C_d}{m}\mathcal B(2\lambda).
\end{aligned}
\]
\begingroup
If \(m\lambda(1+t)>1\), we bound the difference of the two
probabilities by one:
\endgroup

\[
\frac{\left|\mathbb P(R_n\le t)-\mathbb P(\chi_d\le t)\right|}{1+t}
\le\frac1{1+t}<m\lambda.
\]
\begingroup
Both cases therefore give the bound
\(C_dm\lambda+(C_d/m)\mathcal B(2\lambda)\).
Taking the supremum over all \(t\ge0\) and all arrays with
moment sum at most \(\lambda\), as in
\eqref{eq:sec4-small-error}, yields
\endgroup

\begin{equation}\label{eq:sec4-small-cdf-recursion}
\mathcal B(\lambda)\le C_dm\lambda+\frac{C_d}{m}\mathcal B(2\lambda),
\qquad 0<\lambda\le\lambda_0,
\end{equation}

\begingroup
We now iterate this inequality to control the remaining
\(\mathcal B(2\lambda)\) term.
For \(0<\lambda<\lambda_0\), choose \(k\ge1\) with
\(2^{k-1}\lambda<\lambda_0\le2^k\lambda\).
Then \eqref{eq:sec4-small-cdf-recursion} applies successively at
\(\lambda,2\lambda,\ldots,2^{k-1}\lambda\).
At the final argument, we use
\(\mathcal B(2^k\lambda)\le1\).
Since \eqref{eq:sec4-recursion-choice} gives
\(2C_d/m\le1/2\) and \(C_d/m\le1/4\), iteration yields
\endgroup

\[
\begin{aligned}
\mathcal B(\lambda)
&\le C_dm\lambda\sum_{j=0}^{k-1}(2C_d/m)^j
 +(C_d/m)^k\mathcal B(2^k\lambda)\\
&\le2C_dm\lambda+4^{-k}
\le2C_dm\lambda+(\lambda/\lambda_0)^2
\le C_d\lambda.
\end{aligned}
\]

\begingroup
Here \(4^{-k}\le(\lambda/\lambda_0)^2\) follows from
\(2^k\lambda\ge\lambda_0\); \(m\) and \(\lambda_0\) are fixed
and are absorbed into the constant.
For \(\lambda\ge\lambda_0\), the bound
\(\mathcal B(\lambda)\le1\le\lambda/\lambda_0\) gives the same
order. Thus
\endgroup

\begin{equation}\label{eq:sec4-small-cdf-linear}
\mathcal B(\lambda)\le C_d\lambda,\qquad\lambda>0.
\end{equation}

\begingroup
We now return to the Cram\'er-type moderate deviation result
for \(0\le x\le4\). Set
\(\lambda=\sum_i\mathbb E|\boldsymbol\xi_i|^p\).
For \(x>0\), the identity
$\mathbb P(R_n\ge x)
=1-\lim_{t\uparrow x}\mathbb P(R_n\le t)$
converts the distribution function bound into a bound for the
tail event with \(\ge\).
Using \eqref{eq:sec4-small-error} and
\eqref{eq:sec4-small-cdf-linear}, we obtain
\endgroup

\[
\begin{aligned}
\left|\mathbb P(R_n^2\ge x^2)-\mathbb P(\chi_d^2\ge x^2)\right|
&=\lim_{t\uparrow x}
 \left|\mathbb P(R_n\le t)-\mathbb P(\chi_d\le t)\right|\\
&\le(1+x)\mathcal B(\lambda)
\le C_d(1+x)\sum_i\mathbb E|\boldsymbol\xi_i|^p.
\end{aligned}
\]

\begingroup
At \(x=0\), both tail probabilities equal one.
For \(0\le x\le4\), the approximating tail is bounded below by
\[
\mathbb P(\chi_d^2\ge x^2)
\ge\mathbb P(\chi_d^2\ge16)>0.
\]
Dividing the preceding tail probability bound by
\(\mathbb P(\chi_d^2\ge x^2)\) therefore gives
\endgroup

\begin{equation}\label{eq:sec4-small-cmd}
\begin{aligned}
\left|\frac{\mathbb P(R_n^2\ge x^2)}
{\mathbb P(\chi_d^2\ge x^2)}-1\right|
&\le\frac{C_d(1+x)}{\mathbb P(\chi_d^2\ge16)}
 \sum_i\mathbb E|\boldsymbol\xi_i|^p
\le C_d\sum_i\mathbb E|\boldsymbol\xi_i|^p,
\qquad0\le x\le4.
\end{aligned}
\end{equation}
This proves \eqref{eq:CMD-Rn-1} for \(0\le x\le4\).

\begingroup
We now prove the moderate deviation result for \(x>4\).
Our goal is
\endgroup
\[
\left|\mathbb P(R_n^2\ge x^2)-\mathbb P(\chi_d^2\ge x^2)\right|
\le C_d(1+x)^p\sum_i\mathbb E|\boldsymbol\xi_i|^p\,
       \mathbb P(\chi_d^2\ge x^2).
\]
\begingroup
The desired bound contains the tail probability
\(\mathbb P(\chi_d^2\ge x^2)\), which becomes small as \(x\) increases.
We must retain this factor in the estimates below.
With the fixed \(\varepsilon\) in \eqref{eq:sec4-conditioning}
and \(m\) in \eqref{eq:sec4-recursion-choice}, choose
\(\lambda_*>0\), depending only on \(d\), sufficiently small for the auxiliary lemmas
used below, with
\endgroup
\begingroup
\begin{equation}\label{eq:sec4-large-moment-threshold}
\begin{gathered}
m\lambda_*\le1,\qquad
\lambda_*\le\min\{\varepsilon^3/4,1/8\},\\
(1-2\lambda_*^{2/3})^{-3/2}(21/20)^3\le2.
\end{gathered}
\end{equation}
\endgroup
\begingroup
In this part, \(\lambda\) bounds the moment sum multiplied by
\((1+x)^p\): suppose that
\endgroup
\begin{equation}\label{eq:sec4-budget}
(1+x)^p\sum_i\mathbb E|\boldsymbol\xi_i|^p
\le\lambda\le\lambda_*/2.
\end{equation}
\begingroup
For \(x>4\), the cutoff in \eqref{eq:sec4-conditioning} is
\(b=\varepsilon/x\); we continue to use the conditional law
\(\mathbb P^*\) and expectation \(\mathbb E^*\) from that equation.
For the cutoffs \(H_{u,\eta}\) in \eqref{eq:sec4-cutoff}, take
\(0<\eta\le1/(10x)\) and \(u\in\{x-\eta,x\}\).
Use these choices in \(F\) from \eqref{eq:sec4-extension} and
in the quantities defined after
\eqref{eq:sec4-stein-expectation-split} and
\eqref{eq:sec4-three-terms}. These two decompositions give
\[
\left|\mathbb E^*H_{u,\eta}(R_n)-\mathbb EH_{u,\eta}(\chi_d)\right|
\le I_{11}+I_{12}+I_{13}+|I_2|+|I_3|.
\]
To apply the conditional moment bounds to these terms, note that
\endgroup
\[
b^{-p}\sum_i\mathbb E|\boldsymbol\xi_i|^p
\le\varepsilon^{-3}\lambda\le1/4.
\]
\begingroup
The conditioning estimate \eqref{eq:sec4-conditioning-bound}
and Lemma~\ref{lem:sec4-conditional-moments} therefore apply. Substituting
\(b=\varepsilon/x\) and \eqref{eq:sec4-budget} into the moment
bounds gives
\endgroup
\begingroup
the following estimates. Their constants depend only on \(d\), since
\(\varepsilon^{1-p}\le\varepsilon^{-2}\),
\(\varepsilon^{2-p}\le\varepsilon^{-1}\), and
\(\varepsilon^{j-p}\le1\) for \(j=3,4\) and \(2<p\le3\):
\endgroup
\begin{equation}\label{eq:sec4-prob-moments}
\begin{aligned}
\sum_i|\mathbb E^*\boldsymbol\xi_i|
 &\le\frac{C_d\lambda}{x},\qquad
 \|\mathbb E^*V_n-I_d\|_{\mathrm F}\le\frac{C_d\lambda}{x^2},\\
\sum_i\mathbb E^*|\boldsymbol\xi_i|^j
 &\le\frac{C_d\lambda}{x^j},\qquad j=3,4.
\end{aligned}
\end{equation}
\begingroup
We first estimate the two matrix correction terms.
For the event \(\mathcal G\) in \eqref{eq:sec4-good-matrix},
the bound \eqref{eq:sec4-matrix-corrections} gives
\endgroup
\[
|I_2|+|I_3|
\le C_d\mathbb E^*[(1+\|V_n\|_{\mathrm F})\mathbf1_{\mathcal G^c}].
\]
\begingroup
The following lemma bounds this expectation on the required
tail scale.
\endgroup
\begin{lemma}\label{lem:sec4-truncated}
Under the assumptions of Lemma~\ref{lem:sec4-conditional-moments}
and \eqref{eq:sec4-budget}, let \(x\ge4\) and use the
conditional law in \eqref{eq:sec4-conditioning} with \(b=\varepsilon/x\). Then
\begin{equation}\label{eq:sec4-full-matrix}
\mathbb E^*[(1+\|V_n\|_{\mathrm F})\mathbf1_{\mathcal G^c}]
\le C_d\lambda e^{-12x^2}
\le C_d\lambda\mathbb P(\chi_d^2\ge x^2).
\end{equation}
\end{lemma}

Substituting \eqref{eq:sec4-full-matrix} into
\eqref{eq:sec4-matrix-corrections} gives
\begin{equation}\label{eq:sec4-large-matrix-corrections}
|I_2|+|I_3|\le C_d\lambda\mathbb P(\chi_d^2\ge x^2).
\end{equation}

\begingroup
We next estimate \(I_{11}\). Splitting its expectation over
\(\mathcal G_i\) and \(\mathcal G_i^c\), as defined in
\eqref{eq:sec4-good-deleted-matrix}, gives
\endgroup
\[
\begin{aligned}
I_{11}
&=\sum_i|\mathbb E^*\boldsymbol\xi_i|\,
  \mathbb E^*[|F(S_{-i},V_{-i})|;\mathcal G_i]
+\sum_i|\mathbb E^*\boldsymbol\xi_i|\,
  \mathbb E^*[|F(S_{-i},V_{-i})|;\mathcal G_i^c].
\end{aligned}
\]
\begingroup
The uniform bound \(|F|\le C_d\) alone does not supply the tail
factor needed here. Part~\textup{(i)} of the following lemma
gives a more precise bound for \(F\), which we use immediately
in \(I_{11}\). Part~\textup{(ii)} gives the corresponding derivative
bounds for \(I_{12}\) and \(I_{13}\).
\endgroup

\begin{samepage}
\begin{lemma}\label{lem:sec4-solution}
Let $x\ge4$, $0<\eta\le1/(10x)$, and $H=H_{u,\eta}$ with
$u\in\{x-\eta,x\}$.
\begin{enumerate}
\item[(i)]
For $k=1,3$, define
\begin{equation}\label{eq:sec4-stein-envelopes}
A_{k,x}(y)=
\begin{cases}
\mathbb P(\chi_d^2\ge x^2)e^{y^2/2}(1+y)^{k-d},
    &0\le y\le x-x^{-1},\\[3pt]
(1+x)^{k-2},&y>x-x^{-1}.
\end{cases}
\end{equation}
For $s\in\mathbb R^d$ and $v\succeq0$, put
$y=(s^\top K(v)^{-1}s)^{1/2}$. Then
\begin{equation}\label{eq:sec4-stein-field-envelope}
\begin{aligned}
|F(s,v)|&\le C_dA_{1,x}(y),\\
\|\partial_vF(s,v)\|&\le C_dA_{3,x}(y)\le C_d(1+x).
\end{aligned}
\end{equation}
\item[(ii)]
Let $s,z\in\mathbb R^d$ and
$v=v^\top$ satisfy
\[
\|v-I_d\|_{\mathrm F}\le1/8,\qquad |z|\le\varepsilon/x,
\]
and put $y=(s^\top v^{-1}s)^{1/2}$. Then
\begin{equation}\label{eq:sec4-stein-path-matrix}
\sup_{0\le t\le1}
\left\|\partial_vF(s+tz,v+tzz^\top)\right\|
\le C_dA_{3,x}(y),
\end{equation}
and
\begin{equation}\label{eq:sec4-stein-path-jacobian}
\begin{aligned}
&\sup_{0\le t\le1}
\left\|\partial_sF(s+z,v+zz^\top)
       -\partial_sF(s+tz,v+tzz^\top)\right\|_{\mathrm F}\\
&\quad\le C_d|z|A_{3,x}(y)
 +C_d\min\left\{1,\frac{C_d|z|}{\eta}\right\}
       \mathbf1_{\{|y-x|<C_d|z|+\eta\}}.
\end{aligned}
\end{equation}
In the last term, \(C_d|z|\) accounts for the change in radius along
the path near \(x\), while \(\eta\) is the smoothing width.

\end{enumerate}
\end{lemma}
\end{samepage}

\begingroup
On \(\mathcal G_i\), the matrix-valued function \(K(\cdot)\)
in \eqref{eq:sec4-extension} satisfies \(K(V_{-i})=V_{-i}\).
Using the deleted radius \(R_{-i}\) introduced before
Lemma~\ref{lem:sec4-small-increments},
\eqref{eq:sec4-stein-field-envelope} bounds
\(|F(S_{-i},V_{-i})|\) by \(C_dA_{1,x}(R_{-i})\) on this event.
On \(\mathcal G_i^c\), we use \(|F|\le C_d\) from
\eqref{eq:sec4-truncated-field}. Substitution into \(I_{11}\) gives
\endgroup
\[
I_{11}\le C_d\sum_i|\mathbb E^*\boldsymbol\xi_i|
 \left\{\mathbb E^*[A_{1,x}(R_{-i});\mathcal G_i]
              +\mathbb P^*(\mathcal G_i^c)\right\}.
\]
\begingroup
It remains to bound the expectation of \(A_{1,x}(R_{-i})\)
on \(\mathcal G_i\) and the probability of \(\mathcal G_i^c\).
The following lemma provides both estimates, together with the
corresponding expectation for \(A_{3,x}\).
\endgroup

\begin{lemma}
\label{lem:sec4-conditional}
Under the assumptions of Lemma~\ref{lem:sec4-truncated}, use the
conditional law in \eqref{eq:sec4-conditioning} with \(b=\varepsilon/x\).
The following bounds hold uniformly in $i$.
The deleted matrix satisfies
\begin{equation}\label{eq:sec4-prob-bad-matrices}
\sup_{1\le i\le n}
\mathbb E^*[(1+\|V_{-i}\|_{\mathrm F})\mathbf1_{\mathcal G_i^c}]
 \le C_d\lambda e^{-2x^2}.
\end{equation}
For the functions $A_{k,x}$ in Lemma~\ref{lem:sec4-solution}(i),
\begin{equation}\label{eq:sec4-prob-envelope-integrals}
\mathbb E^*[A_{k,x}(R_{-i})\mathbf1_{\mathcal G_i}]
\le C_d\mathbb P(\chi_d^2\ge x^2)(1+x)^k,
\qquad k=1,3.
\end{equation}
\end{lemma}

\begingroup
Apply \eqref{eq:sec4-prob-envelope-integrals} with \(k=1\) and
use \eqref{eq:sec4-prob-bad-matrices} to bound
\(\mathbb P^*(\mathcal G_i^c)\).
Then \eqref{eq:sec4-prob-moments} gives
\endgroup
\begin{equation}\label{eq:sec4-mean-term}
\begin{aligned}
I_{11}
&\le C_d\sum_i|\mathbb E^*\boldsymbol\xi_i|
 \left\{(1+x)\mathbb P(\chi_d^2\ge x^2)+\lambda e^{-2x^2}\right\}
\le C_d\lambda\mathbb P(\chi_d^2\ge x^2).
\end{aligned}
\end{equation}

\begingroup
Here the factor \(1+x\) is absorbed by the bound
\(\sum_i|\mathbb E^*\boldsymbol\xi_i|\le C_d\lambda/x\), since
\(x>4\). This completes the estimate of \(I_{11}\).

We turn to \(I_{12}\), defined after \eqref{eq:sec4-three-terms}.
Its integrand contains the matrix derivative
\(\partial_vF\) along the path obtained by adding back
\(\boldsymbol\xi_i\). Split its expectation over
\(\mathcal G_i\) and \(\mathcal G_i^c\).
On \(\mathcal G_i\), apply \eqref{eq:sec4-stein-path-matrix} with
\endgroup
\begin{equation}\label{eq:sec4-large-path-arguments}
(s,v,z)=(S_{-i},V_{-i},\boldsymbol\xi_i),\qquad y=R_{-i}.
\end{equation}
\begingroup
This bounds the derivative along the path by
\(C_dA_{3,x}(R_{-i})\). On \(\mathcal G_i^c\), use the global
bound \(\|\partial_vF\|\le C_d(1+x)\) from
\eqref{eq:sec4-stein-field-envelope}.
Since \(\boldsymbol\xi_i\) is independent of
\((S_{-i},V_{-i})\) under \(\mathbb P^*\), we obtain
\endgroup
\begin{equation}\label{eq:sec4-matrix-term-target}
I_{12}\le C_d\sum_i\mathbb E^*|\boldsymbol\xi_i|^3
 \left\{\mathbb E^*[A_{3,x}(R_{-i});\mathcal G_i]
             +(1+x)\mathbb P^*(\mathcal G_i^c)\right\}.
\end{equation}
For the first term in \eqref{eq:sec4-matrix-term-target},
\eqref{eq:sec4-prob-envelope-integrals} and
\eqref{eq:sec4-prob-moments} give
\begin{equation}\label{eq:sec4-matrix-envelope-product}
\begin{aligned}
&\sum_i\mathbb E^*|\boldsymbol\xi_i|^3\,
 \mathbb E^*[A_{3,x}(R_{-i});\mathcal G_i]
\le C_d(1+x)^3\mathbb P(\chi_d^2\ge x^2)
             \sum_i\mathbb E^*|\boldsymbol\xi_i|^3
\le C_d\lambda\mathbb P(\chi_d^2\ge x^2).
\end{aligned}
\end{equation}
The second term in \eqref{eq:sec4-matrix-term-target} is at most
\begin{equation}\label{eq:sec4-matrix-bad-product}
C_d(1+x)\lambda e^{-2x^2}
 \sum_i\mathbb E^*|\boldsymbol\xi_i|^3
\le C_d\lambda\mathbb P(\chi_d^2\ge x^2),
\end{equation}
by \eqref{eq:sec4-prob-bad-matrices} and
\eqref{eq:sec4-prob-moments}. Substituting
\eqref{eq:sec4-matrix-envelope-product} and
\eqref{eq:sec4-matrix-bad-product} into
\eqref{eq:sec4-matrix-term-target} gives
\begin{equation}\label{eq:sec4-matrix-good-part}
I_{12}\le C_d\lambda\mathbb P(\chi_d^2\ge x^2).
\end{equation}
Combining
\eqref{eq:sec4-mean-term} and \eqref{eq:sec4-matrix-good-part} gives
\begin{equation}\label{eq:sec4-first-two}
I_{11}+I_{12}\le C_d\lambda\mathbb P(\chi_d^2\ge x^2).
\end{equation}

\begingroup
The estimates of \(I_{11},I_{12},I_2,I_3\) now have the required
factor \(\lambda\mathbb P(\chi_d^2\ge x^2)\).
It remains to estimate \(I_{13}\), which contains the difference
of Jacobians along the path.
Use the split in \eqref{eq:sec4-small-third-split}, with the
present choices of \(b\) and \(H_{u,\eta}\), to estimate
\(I_{131}\) and \(I_{132}\) separately.
For \(I_{131}\), the global Jacobian bound gives
\eqref{eq:sec4-small-bad-matrix-target}.
Using \eqref{eq:sec4-prob-bad-matrices} for the exceptional
probability and \eqref{eq:sec4-conditional-moments} for the
second moments, we obtain
\endgroup
\begin{equation}\label{eq:sec4-large-third-bad}
\begin{aligned}
I_{131}
&\le C_d\sum_i\mathbb E^*|\boldsymbol\xi_i|^2
                   \mathbb P^*(\mathcal G_i^c)
\le C_d\lambda e^{-2x^2}\sum_i\mathbb E^*|\boldsymbol\xi_i|^2
\le C_d\lambda\mathbb P(\chi_d^2\ge x^2).
\end{aligned}
\end{equation}
\begingroup
For \(I_{132}\), apply the Jacobian increment bound
\eqref{eq:sec4-stein-path-jacobian} with \((s,v,z)\) as in
\eqref{eq:sec4-large-path-arguments}. Multiplying by
\(|\boldsymbol\xi_i|^2\) and taking expectations gives
\[
\begin{aligned}
I_{132}
&\le C_d\sum_i\mathbb E^*\left[
 |\boldsymbol\xi_i|^3A_{3,x}(R_{-i});\mathcal G_i\right]\\
&\quad+C_d\sum_i\mathbb E^*\left[
 |\boldsymbol\xi_i|^2
 \min\left\{1,\frac{C_d|\boldsymbol\xi_i|}{\eta}\right\}
 \mathbf1_{\{|R_{-i}-x|<C_d|\boldsymbol\xi_i|+\eta\}}\right].
\end{aligned}
\]
In the second term we have dropped the indicator of
\(\mathcal G_i\), since the integrand is nonnegative.
By independence, the first term is the product already bounded
in \eqref{eq:sec4-matrix-envelope-product}. Hence
\endgroup
\begin{equation}\label{eq:sec4-shell-target}
\begin{aligned}
I_{132}
&\le C_d\lambda\mathbb P(\chi_d^2\ge x^2)
+C_d\sum_i\mathbb E^*\left[
 |\boldsymbol\xi_i|^2
 \min\left\{1,\frac{C_d|\boldsymbol\xi_i|}{\eta}\right\}
 \mathbf1_{\{|R_{-i}-x|<C_d|\boldsymbol\xi_i|+\eta\}}\right].
\end{aligned}
\end{equation}
\begingroup
We must therefore estimate the weighted thin-shell probability
remaining in \eqref{eq:sec4-shell-target}.
The quantity \(\mathcal B(\lambda)\) in
\eqref{eq:sec4-small-error} controls distribution-function
differences. Here we need to compare the tails of \(R_{-i}\)
at the two ends of the shell relative to the corresponding
chi tail probabilities. We therefore introduce a separate
quantity for the error in the tail ratio. For \(\lambda>0\), define
\endgroup
\begin{equation}\label{eq:sec4-error}
\mathcal E(\lambda)=
\sup\left|
\frac{\mathbb P(R_n^2\ge x^2)}{\mathbb P(\chi_d^2\ge x^2)}-1
\right|,
\end{equation}
where the supremum is over \(n\ge1\), \(x\ge0\), and centered
independent arrays with total covariance \(I_d\) and
\[
(1+x)^p\sum_i\mathbb E|\boldsymbol\xi_i|^p\le\lambda.
\]
\begingroup
Our aim is to prove \(\mathcal E(\lambda)\le C_d\lambda\).
At this stage we do not assume this estimate. The following
corollary bounds the deleted tail in terms of
\(\mathcal E(2\lambda)\); the enlarged argument allows for
deleting an observation and normalizing the remaining covariance.
\endgroup

\begin{corollary}\label{lem:sec4-deleted-comparison}
Under the assumptions of Lemma~\ref{lem:sec4-conditional}, for
\(3\le t\le x+1/x\), uniformly in \(i\),
\begin{equation}\label{eq:sec4-prob-deleted-comparison}
\left|\mathbb P^*(R_{-i}\ge t)-\mathbb P(\chi_d^2\ge t^2)\right|
\le C_d\mathbb P(\chi_d^2\ge t^2)
       \{\mathcal E(2\lambda)+\lambda\}.
\end{equation}
\end{corollary}

\begingroup
First take a deterministic width \(0<h\le1/x\). Then
\(3\le x-h<x+h\le x+1/x\), so
\eqref{eq:sec4-prob-deleted-comparison} applies at both endpoints.
The chi density gives
\[
\begin{aligned}
\mathbb P(x-h\le\chi_d<x+h)
 &\le C_dxh\,\mathbb P(\chi_d^2\ge x^2),\\
\mathbb P(\chi_d\ge x-h)
 &\le C_d\mathbb P(\chi_d^2\ge x^2).
\end{aligned}
\]
Subtracting the two deleted tail probabilities and using these
estimates, we obtain
\endgroup
\begin{equation}\label{eq:sec4-shell-probability}
\begin{aligned}
\mathbb P^*(|R_{-i}-x|<h)
&\le\mathbb P^*(R_{-i}\ge x-h)
       -\mathbb P^*(R_{-i}\ge x+h)\\
&\le C_d\mathbb P(\chi_d^2\ge x^2)
 \{xh+\mathcal E(2\lambda)+\lambda\}.
\end{aligned}
\end{equation}
\begingroup
In \eqref{eq:sec4-shell-target}, the width is
\(C_d|\boldsymbol\xi_i|+\eta\), which is random.
Under \(\mathbb P^*\), the truncation and the fixed sufficiently
small \(\varepsilon\) ensure that
\[
0<C_d|\boldsymbol\xi_i|+\eta
\le C_db+\eta\le\frac{C_d\varepsilon+1/10}{x}\le\frac1x
\quad\mathbb P^*\text{-a.s.}
\]
We first condition on \(\boldsymbol\xi_i\). The width is then
fixed, while the law of \(R_{-i}\) is unchanged by independence.
Thus \eqref{eq:sec4-shell-probability} gives
\begin{equation}\label{eq:sec4-large-conditional-shell}
\begin{aligned}
&\mathbb P^*\left(|R_{-i}-x|<C_d|\boldsymbol\xi_i|+\eta
                \mid\boldsymbol\xi_i\right)
\le C_d\mathbb P(\chi_d^2\ge x^2)
 \left\{x(C_d|\boldsymbol\xi_i|+\eta)
              +\mathcal E(2\lambda)+\lambda\right\}.
\end{aligned}
\end{equation}
To insert this estimate into the weighted expectation, use
the two inequalities in \eqref{eq:sec4-shell-weight-inequalities}
with \(w=C_d|\boldsymbol\xi_i|\). They imply
\[
\begin{aligned}
|\boldsymbol\xi_i|^2
 \min\left\{1,\frac{C_d|\boldsymbol\xi_i|}{\eta}\right\}
 (C_d|\boldsymbol\xi_i|+\eta)
 &\le C_d|\boldsymbol\xi_i|^3,\\
|\boldsymbol\xi_i|^2
 \min\left\{1,\frac{C_d|\boldsymbol\xi_i|}{\eta}\right\}
 &\le \frac{C_d|\boldsymbol\xi_i|^3}{\eta}.
\end{aligned}
\]
The first inequality handles the term proportional to the shell
width; the second handles \(\mathcal E(2\lambda)+\lambda\).
Taking expectations and using the third-moment bound in
\eqref{eq:sec4-prob-moments}, we obtain
\endgroup
\begin{equation}\label{eq:sec4-large-shell-contribution}
\begin{aligned}
&\sum_i\mathbb E^*\left[
 |\boldsymbol\xi_i|^2
 \min\left\{1,\frac{C_d|\boldsymbol\xi_i|}{\eta}\right\}
 \mathbf1_{\{|R_{-i}-x|<C_d|\boldsymbol\xi_i|+\eta\}}\right]\\
&\quad\le C_d\mathbb P(\chi_d^2\ge x^2)
 \left\{x\sum_i\mathbb E^*|\boldsymbol\xi_i|^3
 +\frac{\mathcal E(2\lambda)+\lambda}{\eta}
       \sum_i\mathbb E^*|\boldsymbol\xi_i|^3\right\}\\
&\quad\le C_d\mathbb P(\chi_d^2\ge x^2)
 \left\{\frac{\lambda}{x^2}
 +\frac{\lambda}{x^3\eta}
       [\mathcal E(2\lambda)+\lambda]\right\}.
\end{aligned}
\end{equation}
\begingroup
In particular, the coefficient \(\lambda/(x^3\eta)\) comes from
the third-moment sum \(C_d\lambda/x^3\) and the factor
\(1/\eta\) in the second weight inequality.
Substitute \eqref{eq:sec4-large-shell-contribution} into
\eqref{eq:sec4-shell-target}. Returning to
\eqref{eq:sec4-small-third-split} and using
\eqref{eq:sec4-large-third-bad} gives
\endgroup
\begin{equation}\label{eq:sec4-third}
\begin{aligned}
I_{13}=I_{131}+I_{132}
&\le C_d\lambda\mathbb P(\chi_d^2\ge x^2)+I_{132}\\
&\le C_d\mathbb P(\chi_d^2\ge x^2)
 \left\{\lambda+
 \frac{\lambda}{x^3\eta}[\mathcal E(2\lambda)+\lambda]\right\}.
\end{aligned}
\end{equation}

\begingroup
We have now estimated all five terms in the conditional Stein
decomposition. Combining \eqref{eq:sec4-large-matrix-corrections},
\eqref{eq:sec4-first-two}, and \eqref{eq:sec4-third} with
\eqref{eq:sec4-stein-expectation-split} and
\eqref{eq:sec4-three-terms} gives, for \(u\in\{x-\eta,x\}\),
\endgroup
\begin{equation}\label{eq:sec4-stein-expectation}
\begin{aligned}
&\left|\mathbb E^*H_{u,\eta}(R_n)-\mathbb EH_{u,\eta}(\chi_d)\right|
\le C_d\mathbb P(\chi_d^2\ge x^2)
 \left\{\lambda+\frac{\lambda}{x^3\eta}
             [\mathcal E(2\lambda)+\lambda]\right\}.
\end{aligned}
\end{equation}

\begingroup
The remaining task is to pass from these smooth conditional
expectations to the original tail probability. We split the
target difference as
\begin{equation}\label{eq:sec4-large-tail-decomposition}
\begin{aligned}
\left|\mathbb P(R_n^2\ge x^2)-\mathbb P(\chi_d^2\ge x^2)\right|
&\le\left|\mathbb P(R_n^2\ge x^2)
                       -\mathbb P^*(R_n^2\ge x^2)\right|\\
&\quad+\left|\mathbb P^*(R_n^2\ge x^2)
                       -\mathbb P(\chi_d^2\ge x^2)\right|.
\end{aligned}
\end{equation}
We first bound the second term by smoothing. The two cutoffs satisfy
\[
H_{x-\eta,\eta}(y)\le\mathbf1_{\{y<x\}}\le H_{x,\eta}(y),
\qquad y\ge0.
\]
Their expectations therefore bound
\(\mathbb P^*(R_n<x)\) from below and above. Taking complements
gives the desired tail with \(R_n\ge x\). The discrepancy for
\(\chi_d\) is confined to \([x-\eta,x+\eta]\), whose probability
is at most \(C_dx\eta\,\mathbb P(\chi_d^2\ge x^2)\).
Consequently,
\endgroup
\begin{equation}\label{eq:sec4-smoothing}
\begin{aligned}
&\left|\mathbb P^*(R_n^2\ge x^2)-\mathbb P(\chi_d^2\ge x^2)\right|\\
&\quad\le\mathbb P(x-\eta\le\chi_d\le x+\eta)
 +\max_{u\in\{x-\eta,x\}}
 \left|\mathbb E^*H_{u,\eta}(R_n)
             -\mathbb EH_{u,\eta}(\chi_d)\right|\\
&\quad\le C_dx\eta\,\mathbb P(\chi_d^2\ge x^2)
 +\max_{u\in\{x-\eta,x\}}
 \left|\mathbb E^*H_{u,\eta}(R_n)
             -\mathbb EH_{u,\eta}(\chi_d)\right|.
\end{aligned}
\end{equation}

\begingroup
We next bound the first term in
\eqref{eq:sec4-large-tail-decomposition}, which removes the
conditioning. The estimate
\(\mathbb P(\mathcal N^c)\le C_d\lambda\) alone would not suffice:
it lacks the tail factor \(\mathbb P(\chi_d^2\ge x^2)\).
The transfer estimate \eqref{eq:sec4-prelim-transfer} in the
following lemma retains this factor.
\endgroup

\begin{lemma}\label{lem:sec4-preliminary}
Let \(d\ge2\), \(2<p\le3\), and \(x\ge4\).
Use the fixed \(\varepsilon\) and \(\lambda_*\). Let
\(\boldsymbol\xi_1,\ldots,\boldsymbol\xi_n\) be independent centered
random vectors in \(\mathbb R^d\) such that
\begin{equation}\label{eq:sec4-prelim-budget}
\sum_{i=1}^n\mathbb E(\boldsymbol\xi_i\boldsymbol\xi_i^\top)=I_d,
\qquad
(1+x)^p\sum_{i=1}^n\mathbb E|\boldsymbol\xi_i|^p
\le\lambda\le\lambda_*.
\end{equation}
For the statistic in \eqref{def:Rn}, we have
\begin{equation}\label{eq:sec4-prelim-coarse}
\mathbb P(R_n^2\ge x^2)
\le C_d\mathbb P(\chi_d^2\ge x^2).
\end{equation}
Under the conditional law in \eqref{eq:sec4-conditioning} with \(b=\varepsilon/x\),
\(\mathbb P(\mathcal N)\ge3/4\), and
\begin{equation}
\left|\mathbb P^*(R_n^2\ge x^2)-\mathbb P(R_n^2\ge x^2)\right|
\le C_d\lambda\mathbb P(\chi_d^2\ge x^2).
\label{eq:sec4-prelim-transfer}
\end{equation}
\end{lemma}

\begingroup
Substitute \eqref{eq:sec4-stein-expectation} into
\eqref{eq:sec4-smoothing}, and use
\eqref{eq:sec4-prelim-transfer} for the conditioning error in
\eqref{eq:sec4-large-tail-decomposition}.
Dividing by \(\mathbb P(\chi_d^2\ge x^2)\), we obtain
\endgroup
\begin{equation}\label{eq:sec4-one-step}
\begin{aligned}
\left|\frac{\mathbb P(R_n^2\ge x^2)}
                 {\mathbb P(\chi_d^2\ge x^2)}-1\right|
&\le
\frac{|\mathbb P(R_n^2\ge x^2)-\mathbb P^*(R_n^2\ge x^2)|
     +|\mathbb P^*(R_n^2\ge x^2)-\mathbb P(\chi_d^2\ge x^2)|}
     {\mathbb P(\chi_d^2\ge x^2)}\\
&\le C_d\left\{\lambda+x\eta+
 \frac{\lambda}{x^3\eta}[\mathcal E(2\lambda)+\lambda]\right\}.
\end{aligned}
\end{equation}
\begingroup
This estimate has two contributions depending on the smoothing
width: \(x\eta\) from the chi shell and
\(\lambda[\mathcal E(2\lambda)+\lambda]/(x^3\eta)\)
from the deleted shell.
We choose \(\eta\) so that the former is of order \(\lambda\)
and the coefficient of \(\mathcal E(2\lambda)\) is small. Set
\endgroup
\begin{equation}\label{eq:sec4-large-smoothing-width}
\eta=\frac{m\lambda}{x^3}\le\frac1{10x},
\end{equation}
\begingroup
This choice is admissible because
\(m\lambda\le m\lambda_*/2\le1/2\) and \(x>4\) imply
\(x\eta=m\lambda/x^2\le1/(2x^2)<1/10\).
Moreover,
\[
\frac{\lambda}{x^3\eta}=\frac1m,\qquad
x\eta=\frac{m\lambda}{x^2}\le m\lambda.
\]
The remaining term \(\lambda^2/(x^3\eta)=\lambda/m\) is also
at most \(m\lambda\).
Thus \eqref{eq:sec4-one-step} is bounded by
\(C_dm\lambda+(C_d/m)\mathcal E(2\lambda)\) for \(x>4\).
For \(0\le x\le4\), the already established
\eqref{eq:sec4-small-cmd} gives a bound of \(C_d\lambda\).
Taking the supremum over all arrays and thresholds in
\eqref{eq:sec4-error} therefore yields
\endgroup
\begin{equation}\label{eq:sec4-recursion}
\mathcal E(\lambda)\le C_dm\lambda+\frac{C_d}{m}\mathcal E(2\lambda),
\qquad0<\lambda\le\lambda_*/2,
\end{equation}
\begingroup
By the choice of the fixed constant \(m\) in
\eqref{eq:sec4-recursion-choice}, \(C_d/m\le1/4\).
We can therefore iterate this inequality. To control the final
term of the iteration, we first need a bound for
\(\mathcal E\) when its argument reaches a fixed positive size.
The coarse tail estimate \eqref{eq:sec4-prelim-coarse}, together
with \eqref{eq:sec4-small-cmd}, gives
\begin{equation}\label{eq:sec4-large-coarse-error}
\mathcal E(\lambda)\le C_d,\qquad0<\lambda\le\lambda_*.
\end{equation}
This bound follows without using the estimate
\(\mathcal E(\lambda)\le C_d\lambda\) that we are proving.

For \(0<\lambda\le\lambda_*/2\), choose \(k\ge0\) with
\(\lambda_*/2\le2^k\lambda<\lambda_*\).
Apply \eqref{eq:sec4-recursion} successively at
\(\lambda,2\lambda,\ldots,2^{k-1}\lambda\) when \(k\ge1\),
and use \eqref{eq:sec4-large-coarse-error} at the final argument
\(2^k\lambda\). This gives
\endgroup
\[
\begin{aligned}
\mathcal E(\lambda)
&\le C_dm\lambda\sum_{j=0}^{k-1}(2C_d/m)^j
 +(C_d/m)^k\mathcal E(2^k\lambda)\\
&\le2C_dm\lambda+C_d4^{-k}\\
&\le2C_dm\lambda+C_d(2\lambda/\lambda_*)^2
\le C_d\lambda.
\end{aligned}
\]
\begingroup
Here the geometric sum is bounded by \(2\), and
\(4^{-k}\le(2\lambda/\lambda_*)^2\) follows from the choice of
\(k\). Since \(\lambda_*\) and \(m\) are fixed, the last bound
is \(C_d\lambda\).
Finally, choose \(c_d\le\lambda_*/2\) in \eqref{eq:CMD-x-range}
and substitute the moment expression from \eqref{eq:sec4-budget}
into the bound for \(\mathcal E\). By \eqref{eq:sec4-error},
\endgroup
\[
\begin{aligned}
\left|\frac{\mathbb P(R_n^2\ge x^2)}
                 {\mathbb P(\chi_d^2\ge x^2)}-1\right|
&\le\mathcal E\left((1+x)^p\sum_i\mathbb E|\boldsymbol\xi_i|^p\right)
\le C_d(1+x)^p\sum_i\mathbb E|\boldsymbol\xi_i|^p.
\end{aligned}
\]
\begingroup
Combining this estimate with \eqref{eq:sec4-small-cmd} for
\(0\le x\le4\), we obtain \eqref{eq:CMD-Rn-1} for all
\(x\ge0\) satisfying \eqref{eq:CMD-x-range}.
This completes the proof of Theorem~\ref{thm:main-CMD-T}.
\endgroup
\end{proof}

\begingroup
\begin{proof}[Proof of Theorem~\ref{thm:hotelling}]
We derive the Cram\'er-type moderate deviation results
\eqref{eq:CMD-T-1} and \eqref{eq:CMD-two-1}
from Theorem~\ref{thm:main-CMD-T}.
\begingroup
Fix the constants supplied by Theorem~\ref{thm:main-CMD-T}.
We choose the range constant \(c_d\) in
Theorem~\ref{thm:hotelling} sufficiently small relative to these
fixed constants. All choices below depend only on \(d\) and are
uniform for \(2<p\le3\).
\endgroup
The case \(x=0\) is immediate, so fix \(x>0\).

For part~\textup{(i)}, use the standardized vectors
\(\boldsymbol\xi_i\) in \eqref{def:xi}, and write
\[
 \rho=\sum_{i=1}^n\mathbb E|\boldsymbol\xi_i|^p
 =n^{-p/2}\sum_{i=1}^n\mathbb E
   |\Sigma_n^{-1/2}(\mathbf X_i-\boldsymbol\mu_n)|^p,
 \qquad L=(1+x)^p\rho\le c_d.
\]
Their covariance matrices sum to \(I_d\).
Our goal is to show that
\[
 \left|\frac{\mathbb P(T^2_{H,n}\ge x^2)}{Q_d(x^2)}-1\right|
 \le C_dL.
\]
To use Theorem~\ref{thm:main-CMD-T}, we first express this
quantity in terms of \(R_n^2\).
An invertible
affine transformation preserves the Hotelling statistic: under
\(\mathbf X_i\mapsto C\mathbf X_i+\mathbf b\) and
\(\boldsymbol\mu_0\mapsto C\boldsymbol\mu_0+\mathbf b\),
\[
 \bar{\mathbf X}_n-\boldsymbol\mu_0
 \mapsto C(\bar{\mathbf X}_n-\boldsymbol\mu_0),\qquad
 \widehat\Sigma_n\mapsto C\widehat\Sigma_n C^\top.
\]
The inverse transformation cancels in the quadratic form on the
positive-definite event, while singularity is preserved and both
statistics are zero on its complement.

Let \(A=\{\widehat\Sigma_n\succ0\}\). Under the null hypothesis, with
\(\bar{\boldsymbol\xi}=n^{-1}\sum_i\boldsymbol\xi_i\),
\[
 M:=\sum_i(\boldsymbol\xi_i-\bar{\boldsymbol\xi})
              (\boldsymbol\xi_i-\bar{\boldsymbol\xi})^\top
   =V_n-\frac{S_nS_n^\top}{n}
   =\Sigma_n^{-1/2}\widehat\Sigma_n\Sigma_n^{-1/2}.
\]
On \(A\), the matrices \(M\) and \(V_n\) are positive definite, and
\[
 T^2_{H,n}=S_n^\top M^{-1}S_n,\qquad
 R_n^2=S_n^\top V_n^{-1}S_n<n.
\]
The strict inequality follows from
\(V_n^{-1/2}MV_n^{-1/2}
=I_d-V_n^{-1/2}S_nS_n^\top V_n^{-1/2}/n\succ0\).
The Sherman--Morrison formula therefore yields
\[
 M^{-1}=V_n^{-1}
 +\frac{V_n^{-1}S_nS_n^\top V_n^{-1}}{n-R_n^2},
 \qquad
 T^2_{H,n}=\frac{R_n^2}{1-R_n^2/n},\qquad
 R_n^2=\frac{nT^2_{H,n}}{n+T^2_{H,n}}.
\]
Put \(y^2=nx^2/(n+x^2)\). Monotonicity of
\(r\mapsto r/(1-r/n)\) on \([0,n)\) gives
\begin{equation}\label{eq:prob-equality-1}
 \mathbb P(T^2_{H,n}\ge x^2,A)
 =\mathbb P\!\left(R_n^2\ge\frac{nx^2}{n+x^2},A\right).
\end{equation}
Since \(T^2_{H,n}=0\) on \(A^c\) and \(x>0\),
\begin{equation}\label{eq:prob-comparison-1}
 \begin{aligned}
 0&\le\mathbb P(R_n^2\ge y^2)-\mathbb P(T^2_{H,n}\ge x^2)
  =\mathbb P(R_n^2\ge y^2,A^c)\le\mathbb P(A^c).
 \end{aligned}
\end{equation}
Dividing \eqref{eq:prob-comparison-1} by \(Q_d(x^2)\) and
adding and subtracting \(Q_d(y^2)/Q_d(x^2)\), we obtain
\begin{equation}\label{eq:proof-hotelling-product}
 \begin{aligned}
 \left|\frac{\mathbb P(T^2_{H,n}\ge x^2)}{Q_d(x^2)}-1\right|
 &\le\frac{Q_d(y^2)}{Q_d(x^2)}
 \left|\frac{\mathbb P(R_n^2\ge y^2)}{Q_d(y^2)}-1\right|
 +\left(\frac{Q_d(y^2)}{Q_d(x^2)}-1\right)
 +\frac{\mathbb P(A^c)}{Q_d(x^2)}.
 \end{aligned}
\end{equation}
The three terms come from the approximation for \(R_n^2\),
the change of threshold from \(y\) to \(x\), and the singular
event \(A^c\), respectively. We estimate the first two terms
together, then the third.

For the absolute-value factor in the first term,
\(y\le x\) and Theorem~\ref{thm:main-CMD-T} give
\[
 \left|\frac{\mathbb P(R_n^2\ge y^2)}{Q_d(y^2)}-1\right|
 \le C_d(1+y)^p\rho\le C_dL.
\]
To finish estimating the first term, we must also bound
\(Q_d(y^2)/Q_d(x^2)\). We obtain this bound together with
the second term by comparing the two chi-squared tails.
For \(d\ge2\), the representation
\[
 e^{t/2}Q_d(t)=\frac1{2^{d/2}\Gamma(d/2)}
 \int_0^\infty(t+s)^{d/2-1}e^{-s/2}\,ds
\]
is nondecreasing in \(t\). Consequently,
\begin{equation}\label{eq:hotelling-chi-tail-shift}
 0\le\log\frac{Q_d(s)}{Q_d(t)}\le\frac{t-s}{2},
 \qquad 0\le s\le t.
\end{equation}
In particular,
\begin{equation}\label{eq:log-prob-bound-high-low}
 0\le\log\frac{Q_d(nx^2/(n+x^2))}{Q_d(x^2)}
 \le\frac{x^4}{2(n+x^2)}.
\end{equation}
To bound \(x^4/n\) by \(L\), we use the moment assumption.
Jensen's inequality gives
\begin{equation}\label{eq:hotelling-one-moment-lower}
 \rho\ge\sum_{i=1}^n
 (\mathbb E|\boldsymbol\xi_i|^2)^{p/2}
 \ge d^{p/2}n^{-(p-2)/2}.
\end{equation}
Taking \(c_d<1\) ensures \(n\ge2\).

By \eqref{eq:hotelling-one-moment-lower},
\(2p/(p-2)\ge4\), and \(L\le1\),
\begin{equation}\label{eq:proof-CMD-absorption}
 \frac{x^4}{n}
 \le d^{-p/(p-2)}\rho^{2/(p-2)}(1+x)^{2p/(p-2)}
 =d^{-p/(p-2)}L^{2/(p-2)}\le L.
\end{equation}
Thus \(1\le Q_d(y^2)/Q_d(x^2)\le e^L\le e\), and
\begin{equation}\label{eq:proof-chi-ratio}
 0\le\frac{Q_d(y^2)}{Q_d(x^2)}-1
 \le C\frac{x^4}{n}\le C_dL.
\end{equation}
Equation~\eqref{eq:proof-chi-ratio} bounds the second term.
Combining \(Q_d(y^2)/Q_d(x^2)\le e\) with the preceding
application of Theorem~\ref{thm:main-CMD-T} now gives the bound
for the entire first term:
\[
 \frac{Q_d(y^2)}{Q_d(x^2)}
 \left|\frac{\mathbb P(R_n^2\ge y^2)}{Q_d(y^2)}-1\right|
 \le C_dL.
\]

It remains to bound the third term,
\(\mathbb P(A^c)/Q_d(x^2)\). Since
\(M=\Sigma_n^{-1/2}\widehat\Sigma_n\Sigma_n^{-1/2}\),
we have \(A^c=\{\lambda_{\min}(M)=0\}\). We claim that
\begin{equation}\label{eq:hotelling-singularity}
 \mathbb P\{\lambda_{\min}(M)=0\}
 \le d\exp\{-\kappa_p\rho^{-2/(p-2)}\},
 \qquad \kappa_p=\frac1{12}\,4^{-2/(p-2)}.
\end{equation}
To prove this, we center within disjoint blocks, so that their
scatter matrices remain independent. Partition \(\{1,\ldots,n\}\)
into pairs, with one triple when \(n\) is odd. For the resulting
blocks \(I_k\), put
\[
 m_k=|I_k|,\qquad
 \bar{\boldsymbol\xi}_k=\frac1{m_k}\sum_{i\in I_k}\boldsymbol\xi_i,
 \qquad
 B_k=\sum_{i\in I_k}
 (\boldsymbol\xi_i-\bar{\boldsymbol\xi}_k)
 (\boldsymbol\xi_i-\bar{\boldsymbol\xi}_k)^\top.
\]
Then the \(B_k\) are independent and positive semidefinite, and
\[
 M=\sum_k B_k+\sum_k m_k
 (\bar{\boldsymbol\xi}_k-\bar{\boldsymbol\xi})
 (\bar{\boldsymbol\xi}_k-\bar{\boldsymbol\xi})^\top
 \succeq\sum_k B_k.
\]
Independence and centering give
\[
 \mathbb E B_k=\left(1-\frac1{m_k}\right)
 \sum_{i\in I_k}\mathbb E
 (\boldsymbol\xi_i\boldsymbol\xi_i^\top),
 \qquad \sum_k\mathbb E B_k\succeq\frac12 I_d.
\]
With \(\|\cdot\|\) denoting the operator norm,
\(\|B_k\|\le\sum_{i\in I_k}|\boldsymbol\xi_i|^2\) and
\(m_k\le3\) imply
\[
 \sum_k\mathbb E\|B_k\|^{p/2}\le3^{p/2-1}\rho.
\]
Set \(\tau=3(4\rho)^{2/(p-2)}\) and
\(Y_k=B_k\mathbf1_{\{\|B_k\|\le\tau\}}\). These matrices are
independent, \(0\preceq Y_k\preceq\tau I_d\), and
\[
 \left\|\sum_k\mathbb E(B_k-Y_k)\right\|
 \le\tau^{1-p/2}\sum_k\mathbb E\|B_k\|^{p/2}
 \le\frac14.
\]
Thus \(\sum_k\mathbb E Y_k\succeq I_d/4\). The matrix Chernoff
bound of \citet[Corollary~5.2]{Tropp2012}, with \(\delta=1\), yields
\[
 \mathbb P\left\{\lambda_{\min}\left(\sum_kY_k\right)=0\right\}
 \le d\exp\left\{-\frac1{4\tau}\right\}.
\]
Since \(M\succeq\sum_kY_k\succeq0\), singularity of \(M\)
implies singularity of \(\sum_kY_k\). Substitution of \(\tau\)
proves \eqref{eq:hotelling-singularity}.

\begingroup
Returning to the third term in \eqref{eq:proof-hotelling-product},
we make the singularity estimate uniform in \(p\). For
\(\rho\le1/4\) and \(2<p\le3\),
\[
 \kappa_p\rho^{-2/(p-2)}
 =\frac1{12}(4\rho)^{-2/(p-2)}
 \ge\frac1{192\rho^2}.
\]
Also, \(Q_d(x^2)\ge e^{-x^2/2}\) for \(d\ge2\), and
\(L\le c_d\le1\) implies \(x^2\le\rho^{-2/p}\le\rho^{-1}\).
Taking \(c_d\le1/192\), we obtain
\begin{equation}\label{eq:hotelling-one-singular-relative}
 \frac{\mathbb P(A^c)}{Q_d(x^2)}
 \le d\exp\left\{-\frac1{192\rho^2}+\frac{x^2}{2}\right\}
 \le d\exp\left\{-\frac1{384\rho^2}\right\}
 \le C_d\rho\le C_dL.
\end{equation}
\endgroup
All three terms in \eqref{eq:proof-hotelling-product} are
therefore bounded by \(C_dL\). This proves \eqref{eq:CMD-T-1}.

For part~\textup{(ii)}, write
\[
 \begin{aligned}
 \rho={}&
 \frac{\mathbb E|\Sigma_X^{-1/2}
              (\mathbf X_1-\boldsymbol\mu_X)|^p}{n_1^{(p-2)/2}}
 +\frac{\mathbb E|\Sigma_Y^{-1/2}
              (\mathbf Y_1-\boldsymbol\mu_Y)|^p}{n_2^{(p-2)/2}},\\
 L={}&(1+x)^p\rho\le c_d,\qquad m=n_1\wedge n_2.
 \end{aligned}
\]
We need to prove
\[
 \left|\frac{\mathbb P(T_{n_1,n_2}^2\ge x^2)}{Q_d(x^2)}-1\right|
 \le C_dL.
\]
To compare this statistic with the self-normalized form in
Theorem~\ref{thm:main-CMD-T}, we standardize the two samples jointly.
Set \(D=\Sigma_X/n_1+\Sigma_Y/n_2\) and define
\[
 \begin{aligned}
 \boldsymbol\zeta_i&=\frac1{n_1}D^{-1/2}
     (\mathbf X_i-\boldsymbol\mu_X),&&1\le i\le n_1,\\
 \boldsymbol\zeta_{n_1+j}&=-\frac1{n_2}D^{-1/2}
     (\mathbf Y_j-\boldsymbol\mu_Y),&&1\le j\le n_2.
 \end{aligned}
\]
These vectors are independent and centered, with total covariance
\(I_d\). The matrices
\(n_1^{-1/2}D^{-1/2}\Sigma_X^{1/2}\) and
\(n_2^{-1/2}D^{-1/2}\Sigma_Y^{1/2}\) have operator norm at most
one, because their products with their transposes sum to \(I_d\).
It follows that
\[
 \sum_{i=1}^{n_1+n_2}\mathbb E|\boldsymbol\zeta_i|^p\le\rho.
\]
Put
\[
 S_X=\sum_{i=1}^{n_1}\boldsymbol\zeta_i,\qquad
 S_Y=\sum_{j=1}^{n_2}\boldsymbol\zeta_{n_1+j},\qquad
 V=\sum_{i=1}^{n_1+n_2}\boldsymbol\zeta_i\boldsymbol\zeta_i^\top.
\]
Then the standardized sample denominator is
\[
 W=D^{-1/2}\left(\frac{\widehat\Sigma_X}{n_1}
                 +\frac{\widehat\Sigma_Y}{n_2}\right)D^{-1/2}
  =V-\frac{S_XS_X^\top}{n_1}-\frac{S_YS_Y^\top}{n_2}.
\]
Let \(A=\{W\succ0\}\). Under the null hypothesis, on \(A\),
\[
 T_{n_1,n_2}^2=(S_X+S_Y)^\top W^{-1}(S_X+S_Y).
\]
Replacing \(W\) by \(V\) gives the self-normalized statistic
\[
 R_+^2=(S_X+S_Y)^\top V^+(S_X+S_Y),
\]
defined on the whole probability space. Here \(V^+\) is the
Moore--Penrose inverse. Since \(S_X+S_Y\) lies in the range of
\(V\), this agrees with \eqref{def:Rn} for the combined array.
Thus Theorem~\ref{thm:main-CMD-T} applies to \(R_+^2\), with
moment sum at most \(\rho\).

We first prove the lower bound for the tail ratio. On \(A\),
\(V\succ0\) and \(W\preceq V\), so
\(T_{n_1,n_2}^2\ge R_+^2\). Consequently,
\[
 \mathbb P(T_{n_1,n_2}^2\ge x^2)
 \ge\mathbb P(R_+^2\ge x^2,A)
 \ge\mathbb P(R_+^2\ge x^2)-\mathbb P(A^c).
\]
Dividing by \(Q_d(x^2)\) gives
\begin{equation}\label{eq:hotelling-two-lower-error}
 1-\frac{\mathbb P(T_{n_1,n_2}^2\ge x^2)}{Q_d(x^2)}
 \le\left|\frac{\mathbb P(R_+^2\ge x^2)}{Q_d(x^2)}-1\right|
     +\frac{\mathbb P(A^c)}{Q_d(x^2)}.
\end{equation}
The first term is bounded directly by
Theorem~\ref{thm:main-CMD-T}:
\[
 \left|\frac{\mathbb P(R_+^2\ge x^2)}{Q_d(x^2)}-1\right|
 \le C_d(1+x)^p\sum_i\mathbb E|\boldsymbol\zeta_i|^p
 \le C_dL.
\]
For the second term, we use the singularity estimate
\eqref{eq:hotelling-singularity}. Each standardized \(p\)th
moment is at least \(d^{p/2}\), so
\begin{equation}\label{eq:hotelling-two-moment-lower}
 \rho\ge d^{p/2}m^{-(p-2)/2},\qquad
 \frac1m\le d^{-p/(p-2)}\rho^{2/(p-2)}.
\end{equation}
In particular, \(n_1,n_2\ge2\) when \(c_d<1\).
If the sum of two positive semidefinite matrices is singular,
each summand is singular. Apply \eqref{eq:hotelling-singularity}
to the independent centered vectors
\(n_1^{-1/2}\Sigma_X^{-1/2}(\mathbf X_i-\boldsymbol\mu_X)\),
whose covariance matrices sum to \(I_d\) and whose \(p\)th-moment
sum is at most \(\rho\). This gives
\[
 \mathbb P(A^c)\le\mathbb P\{\widehat\Sigma_X\text{ is singular}\}
 \le d\exp\{-\kappa_p\rho^{-2/(p-2)}\}.
\]
Since \(x^2\le\rho^{-2/p}\) and \(\rho\le c_d\), the
calculation in \eqref{eq:hotelling-one-singular-relative} yields
\begin{equation}\label{eq:hotelling-two-singular-relative}
 \frac{\mathbb P(A^c)}{Q_d(x^2)}\le C_d\rho\le C_dL.
\end{equation}
Both terms in \eqref{eq:hotelling-two-lower-error} are now
controlled, and hence
\[
 \frac{\mathbb P(T_{n_1,n_2}^2\ge x^2)}{Q_d(x^2)}
 \ge1-C_dL.
\]

It remains to prove the upper bound. For this, we need to bound
\(W\) from below in terms of \(V\), which requires control of
both centering terms in \(V-W\). Along with \(R_+^2\), introduce
\[
 R_-^2=(S_X-S_Y)^\top V^+(S_X-S_Y).
\]
The vector \(S_X-S_Y\) also lies in the range of \(V\). Thus
\(R_-^2\) is the statistic in \eqref{def:Rn} for the combined
array with the signs in the second group reversed. This change
preserves independence, centering, total covariance, and the
moment sum, so Theorem~\ref{thm:main-CMD-T} also applies to
\(R_-^2\).
On \(A\), the identity
\[
 S_X^\top V^{-1}S_X+S_Y^\top V^{-1}S_Y
 =\frac{R_+^2+R_-^2}{2}
\]
gives
\[
 \begin{aligned}
 \|V^{-1/2}(V-W)V^{-1/2}\|_{\mathrm{op}}
 &\le\frac{S_X^\top V^{-1}S_X}{n_1}
      +\frac{S_Y^\top V^{-1}S_Y}{n_2}
 \le\frac{R_+^2+R_-^2}{2m}.
 \end{aligned}
\]
We therefore separate off the event that \(R_-^2\) is large.
For now, take \(v>0\) such that
\(\delta=(x^2+v^2)/(2m)<1\), and put
\(y=x\sqrt{1-\delta}\). On
\(A\cap\{R_+^2<y^2,\ R_-^2<v^2\}\), the preceding bound
is at most \(\delta\), so \(W\succeq(1-\delta)V\) and
\[
 T_{n_1,n_2}^2\le\frac{R_+^2}{1-\delta}<x^2.
\]
Thus, on \(A\), the event \(\{T_{n_1,n_2}^2\ge x^2\}\)
requires either \(R_+^2\ge y^2\) or \(R_-^2\ge v^2\).
Since the statistic is zero on \(A^c\) and \(x>0\),
\begin{equation}\label{eq:hotelling-two-comparison}
 \mathbb P(T_{n_1,n_2}^2\ge x^2)
 \le\mathbb P(R_+^2\ge y^2)+\mathbb P(R_-^2\ge v^2).
\end{equation}
Dividing by \(Q_d(x^2)\) and adding and subtracting
\(Q_d(y^2)/Q_d(x^2)\), we obtain
\begin{equation}\label{eq:hotelling-two-upper-error}
 \begin{aligned}
 \frac{\mathbb P(T_{n_1,n_2}^2\ge x^2)}{Q_d(x^2)}-1
 &\le\frac{Q_d(y^2)}{Q_d(x^2)}
 \left|\frac{\mathbb P(R_+^2\ge y^2)}{Q_d(y^2)}-1\right|
 +\left(\frac{Q_d(y^2)}{Q_d(x^2)}-1\right)
 +\frac{\mathbb P(R_-^2\ge v^2)}{Q_d(x^2)}.
 \end{aligned}
\end{equation}
The first two terms involve the same chi-squared tail ratio,
as in part~\textup{(i)}. The third term comes from the event
separated off above. We choose \(v\) to make this extra tail
small while keeping the change from \(x\) to \(y\) small.

Take \(c_d<e^{-1}\) and set
\[
 v^2=2x^2+4\log(1/\rho),\qquad
 \delta=\frac{3x^2+4\log(1/\rho)}{2m}.
\]
We first check that this choice satisfies \(\delta<1\), as
required for \eqref{eq:hotelling-two-comparison}.
By \eqref{eq:hotelling-two-moment-lower} and
\(x^2\le\rho^{-2/p}\),
\begingroup
\[
 \begin{aligned}
 \delta&\le C_d\left\{
 \rho^{2/(p-2)-2/p}+\rho^{2/(p-2)}\log(1/\rho)\right\}
 \le C_d\{\rho+\rho^2\log(1/\rho)\}\le\frac12.
 \end{aligned}
\]
Here we used \(2/(p-2)-2/p\ge1\) and \(2/(p-2)\ge2\),
and chose \(c_d\) sufficiently small independently of \(p\).
\endgroup Hence \(y=x\sqrt{1-\delta}\)
is well defined with \(0<y\le x\), and
\eqref{eq:hotelling-two-upper-error} applies.

For the absolute-value factor in the first term,
Theorem~\ref{thm:main-CMD-T} gives
\[
 \left|\frac{\mathbb P(R_+^2\ge y^2)}{Q_d(y^2)}-1\right|
 \le C_d(1+y)^p\sum_i\mathbb E|\boldsymbol\zeta_i|^p
 \le C_dL.
\]
To complete the first term, we must also bound
\(Q_d(y^2)/Q_d(x^2)\). We do this together with the second term.
Since \(x^2-y^2=x^2\delta\),
\eqref{eq:hotelling-chi-tail-shift} reduces the required
chi-squared tail comparison to a bound on \(x^2\delta\).
By \eqref{eq:hotelling-two-moment-lower},
\(2p/(p-2)\ge4\), and \(p\ge2\),
\[
 \begin{aligned}
 x^2\delta
 &\le C_d\rho^{2/(p-2)}
       \{(1+x)^4+(1+x)^2\log(1/\rho)\}\\
 &\le C_d\{L^{2/(p-2)}
       +L\rho^{2/(p-2)-1}\log(1/\rho)\}
 \le C_dL.
 \end{aligned}
\]
\begingroup
The last inequality is uniform in \(p\), because \(L\le1\) gives
\(L^{2/(p-2)}\le L\), and
\[
 \rho^{2/(p-2)-1}\log(1/\rho)
 \le\rho\log(1/\rho)\le1/e.
\]
\endgroup
Consequently,
\[
 1\le\frac{Q_d(y^2)}{Q_d(x^2)}
 \le e^{x^2\delta/2}\le1+C_dL.
\]
This bounds the second term in \eqref{eq:hotelling-two-upper-error}.
Together with the preceding application of
Theorem~\ref{thm:main-CMD-T}, it also bounds the entire first term:
\[
 \frac{Q_d(y^2)}{Q_d(x^2)}
 \left|\frac{\mathbb P(R_+^2\ge y^2)}{Q_d(y^2)}-1\right|
 \le(1+C_dL)C_dL\le C_dL.
\]

It remains to bound the third term,
\(\mathbb P(R_-^2\ge v^2)/Q_d(x^2)\).
\begingroup
To apply Theorem~\ref{thm:main-CMD-T} at \(v\), use \(p\le3\)
and \(\log(1/\rho)\ge1\) to obtain
\[
 (1+v)^p\sum_i\mathbb E|\boldsymbol\zeta_i|^p
 \le C\{L+\rho[\log(1/\rho)]^{3/2}\}
 \le C\sqrt{c_d},
\]
where \(C\) is absolute. The last inequality uses
\(\rho[\log(1/\rho)]^{3/2}\le2\sqrt\rho\) and
\(\rho\le L\le c_d\le1\). Thus we can choose \(c_d\),
independently of \(p\), so that this expression lies within the
fixed range of Theorem~\ref{thm:main-CMD-T}.
\endgroup
The theorem therefore gives
\(\mathbb P(R_-^2\ge v^2)\le C_dQ_d(v^2)\).
Using \(e^{-t/2}\le Q_d(t)\le2^{d/2}e^{-t/4}\), valid for
\(t\ge0\) and \(d\ge2\), and the chosen value of \(v\), we obtain
\begin{equation}\label{eq:hotelling-two-auxiliary-tail}
 \begin{aligned}
 \frac{\mathbb P(R_-^2\ge v^2)}{Q_d(x^2)}
 &\le C_d\frac{Q_d(v^2)}{Q_d(x^2)}
 \le C_d\exp\{-v^2/4+x^2/2\}
 =C_d\rho\le C_dL.
 \end{aligned}
\end{equation}
All three terms in \eqref{eq:hotelling-two-upper-error} are now
bounded by \(C_dL\), giving
\[
 \frac{\mathbb P(T_{n_1,n_2}^2\ge x^2)}{Q_d(x^2)}
 \le1+C_dL.
\]
Together with the lower bound proved above, this establishes
\eqref{eq:CMD-two-1}.
\end{proof}
\endgroup

\begingroup
\begin{proof}[Proof of Corollary~\ref{cor:unif}]
We first prove \eqref{eq:BE-all}. Fix the range constant \(c_d\)
from Theorem~\ref{thm:main-CMD-T}. If
\(\sum_{i=1}^n\mathbb E|\boldsymbol\xi_i|^p\ge c_d\), then
\[
 \sup_{t\ge0}\left|\mathbb P(R_n^2\ge t)-Q_d(t)\right|
 \le1\le\frac1{c_d}\sum_{i=1}^n\mathbb E|\boldsymbol\xi_i|^p.
\]
Otherwise, let \(a>0\) be the endpoint of the moderate deviation
range, determined by
\[
 (1+a)^p\sum_{i=1}^n\mathbb E|\boldsymbol\xi_i|^p=c_d.
\]
For \(0\le x\le a\), Theorem~\ref{thm:main-CMD-T} gives
\[
 \begin{aligned}
 \left|\mathbb P(R_n^2\ge x^2)-Q_d(x^2)\right|
 &\le C_d(1+x)^pQ_d(x^2)
       \sum_{i=1}^n\mathbb E|\boldsymbol\xi_i|^p
 \le C_d\sum_{i=1}^n\mathbb E|\boldsymbol\xi_i|^p.
 \end{aligned}
\]
Here we used \(p\le3\) and the exponential decay of the chi-squared
tail, which imply
\[
 (1+x)^pQ_d(x^2)\le(1+x)^3Q_d(x^2)\le C_d,
 \qquad x\ge0.
\]

To handle \(x>a\), apply Theorem~\ref{thm:main-CMD-T} at \(a\):
\[
 \mathbb P(R_n^2\ge a^2)\le(1+C_dc_d)Q_d(a^2).
\]
The preceding chi-squared tail bound and the defining equality for
\(a\) also give
\[
 Q_d(a^2)\le\frac{C_d}{(1+a)^p}
 =\frac{C_d}{c_d}\sum_{i=1}^n\mathbb E|\boldsymbol\xi_i|^p.
\]
Thus, by monotonicity of both tails, for \(x>a\),
\[
 \begin{aligned}
 \left|\mathbb P(R_n^2\ge x^2)-Q_d(x^2)\right|
 &\le\mathbb P(R_n^2\ge a^2)+Q_d(a^2)\\
 &\le(2+C_d c_d)Q_d(a^2)
 \le C_d\sum_{i=1}^n\mathbb E|\boldsymbol\xi_i|^p.
 \end{aligned}
\]
Combining the two ranges proves \eqref{eq:BE-all}, since
\(t=x^2\) covers all nonnegative thresholds.

Applying the same argument to \eqref{eq:CMD-T-1} and
\eqref{eq:CMD-two-1}, with their respective moment expressions,
proves \eqref{eq:BE-one} and \eqref{eq:BE-two}, respectively.
\end{proof}
\endgroup

\begingroup
\begin{proof}[Proof of Corollary~\ref{cor:iid}]
For one sample, the i.i.d. assumption and \eqref{def:xi} give
\begin{equation}\label{eq:iid-moment-sum}
 \sum_{i=1}^n\mathbb E|\boldsymbol\xi_i|^p
 =n^{-(p-2)/2}\mathbb E
 \left|\Sigma_n^{-1/2}(\mathbf X_1-\boldsymbol\mu_n)\right|^p
 =O\!\left(n^{-(p-2)/2}\right).
\end{equation}
Substitution into \eqref{eq:BE-one} gives \eqref{eq:iid-BE-rate}.
If \(b_n=o(n^{(p-2)/(2p)})\), then
\((1+b_n)^p n^{-(p-2)/2}\to0\).
Together with \eqref{eq:iid-moment-sum}, this gives
\[
 \sup_{0\le x\le b_n}(1+x)^p
 \sum_{i=1}^n\mathbb E|\boldsymbol\xi_i|^p\longrightarrow0.
\]
Thus, for all sufficiently large \(n\), the range condition in
Theorem~\ref{thm:hotelling}\textup{(i)} holds throughout
\(0\le x\le b_n\). Applying \eqref{eq:CMD-T-1} on this interval
yields \eqref{eq:iid-CMD-rate}.

For two samples, the uniform moment bound gives
\[
 \frac{\mathbb E|\Sigma_X^{-1/2}
              (\mathbf X_1-\boldsymbol\mu_X)|^p}{n_1^{(p-2)/2}}
 +\frac{\mathbb E|\Sigma_Y^{-1/2}
              (\mathbf Y_1-\boldsymbol\mu_Y)|^p}{n_2^{(p-2)/2}}
 =O\!\left((n_1\wedge n_2)^{-(p-2)/2}\right).
\]
Substitution into \eqref{eq:BE-two} gives \eqref{eq:iid-BE-two-rate}.
If \(b_{n_1,n_2}=o((n_1\wedge n_2)^{(p-2)/(2p)})\), then
\[
 (1+b_{n_1,n_2})^p(n_1\wedge n_2)^{-(p-2)/2}\longrightarrow0.
\]
Together with the preceding two-sample moment estimate, this shows
that the range condition in Theorem~\ref{thm:hotelling}\textup{(ii)}
holds throughout \(0\le x\le b_{n_1,n_2}\) for all sufficiently
large \(n_1\wedge n_2\). Applying \eqref{eq:CMD-two-1} on this
interval yields \eqref{eq:iid-CMD-two-rate}.
\end{proof}
\endgroup

\begingroup
\begin{proof}[Proof of Example~\ref{thm:sharpness}]
Choose an integer \(M>\max\{10,2(p+1)/(p-2)\}\). For all
sufficiently large integers \(k\), set
\[
 m_k=\lceil k^M\rceil,\qquad n_k=m_k+1,\qquad x_k=k,
 \qquad C_k=I_d-k^{-3}e_1e_1^\top,
\]
where \(e_1=(1,0,\ldots,0)^\top\). Let
\(G_{k1},\ldots,G_{k m_k}\) be independent \(N_d(0,I_d)\)
vectors, independent of \(H_k\), with
\[
 \mathbb P(H_k=0)=1-1/k,\qquad
 \mathbb P(H_k=k^{-1}e_1)=\mathbb P(H_k=-k^{-1}e_1)=1/(2k).
\]
Define
\[
 \boldsymbol\xi_{ki}=m_k^{-1/2}C_k^{1/2}G_{ki}\quad(i\le m_k),
 \qquad \boldsymbol\xi_{k n_k}=H_k,
 \qquad \mathbf X_{ki}=\sqrt{n_k}\boldsymbol\xi_{ki}.
\]
The observations are independent and centered, with
\(n_k^{-1}\sum_i\mathbb E(\mathbf X_{ki}\mathbf X_{ki}^\top)
=C_k+k^{-3}e_1e_1^\top=I_d\). Their centered sample covariance
is positive definite almost surely because the first \(d+1\)
Gaussian observations are affinely independent. Moreover,
\[
 n_k^{-p/2}\sum_i\mathbb E|\mathbf X_{ki}|^p
 =\sum_i\mathbb E|\boldsymbol\xi_{ki}|^p
 =m_k^{1-p/2}\mathbb E|C_k^{1/2}G_{k1}|^p+k^{-(p+1)}
 \sim k^{-(p+1)},
\]
since \(M(p-2)/2>p+1\). Hence
\begin{equation}\label{eq:sharpness-moment-scale}
 (1+k)^p n_k^{-p/2}\sum_i\mathbb E|\mathbf X_{ki}|^p\sim k^{-1}.
\end{equation}

To compute the tail, put
\[
 S_{0k}=\sum_{i\le m_k}\boldsymbol\xi_{ki},\quad
 V_{0k}=\sum_{i\le m_k}\boldsymbol\xi_{ki}\boldsymbol\xi_{ki}^\top,
 \quad
 D_k=(S_{0k}+H_k)^\top(C_k+H_kH_k^\top)^{-1}(S_{0k}+H_k).
\]
Here \(S_{0k}\sim N_d(0,C_k)\) and \(\mathbb EV_{0k}=C_k\).
Conditional on \(H_k=0\), \(D_k\sim\chi_d^2\); conditional on
either nonzero value of \(H_k\), Gaussian symmetry gives
\[
 D_k\ \stackrel{\mathrm d}=
 \frac{(Z_1+a_k)^2}{1+a_k^2}+\sum_{j=2}^dZ_j^2,
 \qquad a_k=\frac1{k\sqrt{1-k^{-3}}},\qquad Z\sim N_d(0,I_d).
\]
The density of the transformed Gaussian vector, relative to the
standard Gaussian density, is
\begin{equation}\label{eq:sharpness-density-ratio}
 L_k(z)=\sqrt{1+a_k^2}
 \exp\left\{-\frac{a_k^2z_1^2}{2}
       +a_k\sqrt{1+a_k^2}\,z_1-\frac{a_k^2}{2}\right\}.
\end{equation}
Completing the square gives
\(0<L_k(z)\le\sqrt{1+a_k^2}e^{1/2}\le C\).
Write \(Z=|Z|U\), with \(U\) uniform on \(\mathbb S^{d-1}\)
and independent of \(|Z|\). For every \(\delta>0\),
\begin{equation}\label{eq:sharpness-radial-concentration}
 \mathbb P\{|Z|\ge(1+\delta)k\mid |Z|\ge k\}
 =\frac{Q_d((1+\delta)^2k^2)}{Q_d(k^2)}\longrightarrow0.
\end{equation}
Thus \(|Z|/k\to1\) under this conditioning. Since \(ka_k\to1\),
bounded convergence and symmetry give
\[
 \frac{\mathbb P(D_k\ge k^2\mid H_k\ne0)}{Q_d(k^2)}
 =\mathbb E\{L_k(Z)\mid |Z|\ge k\}
 \longrightarrow\mathbb E[e^{-U_1^2/2}\cosh(U_1)].
\]
Averaging over \(H_k\) therefore yields
\begin{equation}\label{eq:sharpness-comparison-tail}
 \frac{\mathbb P(D_k\ge k^2)}{Q_d(k^2)}
 =1-\frac{1-\mathbb E[e^{-U_1^2/2}\cosh(U_1)]}{k}+o(k^{-1}).
\end{equation}

We finally verify that replacing \(D_k\) by \(T^2_{H,n_k}\)
does not affect this expansion. Gaussian covariance concentration
gives
\begin{equation}\label{eq:sharpness-matrix-error}
 \mathbb P\{\|V_{0k}-C_k\|_{\mathrm{op}}>k^{-4}\}
 \le C_de^{-c_dm_kk^{-8}}=o\{k^{-1}Q_d(k^2)\},
\end{equation}
because \(M>10\). On the complementary event,
\(C_k\succeq I_d/2\) implies, for \(R_{n_k}^2\) formed from
\(\boldsymbol\xi_{k1},\ldots,\boldsymbol\xi_{k n_k}\),
\begin{equation}\label{eq:sharpness-quadratic-comparison}
 \frac{D_k}{1+2k^{-4}}\le R_{n_k}^2\le\frac{D_k}{1-2k^{-4}}.
\end{equation}
By \eqref{eq:prob-equality-1}, the two corresponding thresholds
for \(D_k\) are
\[
 (1\pm2k^{-4})\frac{n_kk^2}{n_k+k^2}=k^2+O(k^{-2}),
\]
since \(k^4/n_k=O(k^{4-M})=o(k^{-2})\).
The bounded density ratio above also bounds the density of
\(D_k\) by a constant times the chi-squared density. Consequently,
\[
\begin{aligned}
 &\left|\mathbb P(T^2_{H,n_k}\ge k^2)-\mathbb P(D_k\ge k^2)\right|\\
 &\quad\le\mathbb P\{\|V_{0k}-C_k\|_{\mathrm{op}}>k^{-4}\}
       +C_dk^{-2}Q_d(k^2)
 =o\{k^{-1}Q_d(k^2)\}.
\end{aligned}
\]
Thus \eqref{eq:sharpness-comparison-tail} holds for \(T^2_{H,n_k}\).
Since \(\log\cosh t<t^2/2\) for \(t\ne0\),
\[
 c_{d,p}:=1-\mathbb E[e^{-U_1^2/2}\cosh(U_1)]\in(0,1).
\]
Together with \eqref{eq:sharpness-moment-scale}, this proves
\eqref{eq:one-sample-sharp-rate}.
\end{proof}
\endgroup

\begingroup
\begin{proof}[Proof of Example~\ref{thm:two-sample-sharpness}]
Choose an integer \(M>\max\{12,2p/(p-2)\}\). For all sufficiently
large integers \(k\), let \(n_k=\lceil k^M\rceil\), \(x_k=k\), and
\[
 \lambda_k=k^pn_k^{-(p-2)/2},\qquad
 \pi_k=\frac{\lambda_k}{2n_k},\qquad a_k=\frac{\sqrt{2n_k}}k.
\]
Take all \(2n_k\) observations to be independent copies of
\(G+a_kJ_ke_1\), where \(G\sim N_d(0,I_d)\),
\(e_1=(1,0,\ldots,0)^\top\), and \(J_k\) is independent of \(G\), with
\[
 \mathbb P(J_k=0)=1-\pi_k,\qquad
 \mathbb P(J_k=1)=\mathbb P(J_k=-1)=\pi_k/2.
\]
These observations are centered, with
\(\Sigma_k=I_d+(\lambda_k/k^2)e_1e_1^\top\).
Since \(\Sigma_k\succeq I_d\) and \(\pi_ka_k^p=2^{p/2-1}\),
\[
 \mathbb E|\Sigma_k^{-1/2}\mathbf X_{k1}|^p
 \le2^{p-1}\{\mathbb E|G|^p+2^{p/2-1}\},
\]
uniformly in \(k\). The positive Gaussian mixture density also
ensures that both centered sample covariances are positive
definite almost surely.

Write \(n=n_k\), and let \(K\) count the nonzero shifts.
Then \(K\sim\operatorname{Binomial}(2n,\pi_k)\),
with mean \(\lambda_k\to0\), and
\[
 \mathbb P(K=0)=1-\lambda_k+O(\lambda_k^2),\qquad
 \mathbb P(K=1)=\lambda_k+O(\lambda_k^2).
\]
We evaluate the tail through the decomposition
\begin{equation}\label{eq:two-sharp-tail-decomposition}
\begin{aligned}
 \mathbb P(T_{n,n}^2\ge k^2)
 &=\mathbb P(K=0)\mathbb P(T_{n,n}^2\ge k^2\mid K=0)\\
 &\quad+\mathbb P(K=1)\mathbb P(T_{n,n}^2\ge k^2\mid K=1)\\
 &\quad+\mathbb P(T_{n,n}^2\ge k^2,\ K\ge2).
\end{aligned}
\end{equation}

For the first term, the Gaussian two-sample law gives
\[
 T_{n,n}^2\mid\{K=0\}\ \stackrel{\mathrm d}=\frac{2nA}{B},
 \qquad A\sim\chi_d^2,\quad B\sim\chi_{2n-d-1}^2,\quad A\perp B.
\]
The density ratio to \(\chi_d^2\), uniformly on \([k^2,2k^2]\), is
\[
 \{1+O_d(n^{-1})\}e^{t/2}
       \left(1+\frac{t}{2n}\right)^{-n+1/2}
 =1+O_d(k^4/n),
\]
by Stirling's formula and a Taylor expansion. Beyond \(2k^2\),
\(\mathbb P(2nA/B\ge2k^2)\le Q_d(3k^2/2)+C_de^{-c_dn}\).
Integrating the density estimate therefore gives
\begin{equation}\label{eq:two-sharp-zero-shifts}
 \frac{\mathbb P(T_{n,n}^2\ge k^2\mid K=0)}{Q_d(k^2)}
 =1+O_d(k^4/n)+O_d(e^{-c_dk^2})=1+o(\lambda_k),
\end{equation}
since \((k^4/n)/\lambda_k=(k/\sqrt n)^{4-p}\to0\).

For the remaining terms we retain the shifts in both the mean
and the centered covariance. Stack the Gaussian parts into a
\(2n\times d\) matrix \(\mathsf G\), and put
\(\boldsymbol j=(J_{k1},\ldots,J_{k,2n})^\top\). Let
\[
 P=\operatorname{diag}(P_n,P_n),\qquad
 P_n=I_n-\boldsymbol1_n\boldsymbol1_n^\top/n,\qquad
 \boldsymbol c=(\boldsymbol1_n^\top,-\boldsymbol1_n^\top)^\top.
\]
Thus \(P\) centers each group and \(\boldsymbol c\) takes their
difference. With
\[
 Z=(2n)^{-1/2}\mathsf G^\top\boldsymbol c,\qquad
 h=\boldsymbol c^\top\boldsymbol j/k,\qquad
 W=\frac{(\mathsf G+a_k\boldsymbol j e_1^\top)^\top
             P(\mathsf G+a_k\boldsymbol j e_1^\top)}{2n},
\]
the statistic is exactly
\begin{equation}\label{eq:two-sharp-conditional-statistic}
 T_{n,n}^2=(Z+he_1)^\top W^{-1}(Z+he_1).
\end{equation}
Conditional on \(\boldsymbol j\), \(Z\sim N_d(0,I_d)\) and is
independent of \(P\mathsf G\), since \(P\boldsymbol c=0\).
Write \(L_X=\sum_{i\le n}J_{ki}\) and
\(L_Y=\sum_{i>n}J_{ki}\). Then
\begin{equation}\label{eq:two-sharp-centering}
 \boldsymbol j^\top P\boldsymbol j=K-\frac{L_X^2+L_Y^2}{n},
 \qquad L_X^2+L_Y^2\le K^2,\qquad |h|\le K/k.
\end{equation}

For \(1\le K\le k^2\), compare \(W\) with
\(A_K=I_d+(K/k^2)e_1e_1^\top\). Expanding \(W\) gives
\[
 W-A_K=\left\{\frac{\mathsf G^\top P\mathsf G}{2n}
                    -(1-1/n)I_d\right\}
       +ve_1^\top+e_1v^\top-\frac{I_d}{n}
       -\frac{L_X^2+L_Y^2}{nk^2}e_1e_1^\top,
\]
where \(v=(a_k/(2n))\mathsf G^\top P\boldsymbol j\) is
conditionally Gaussian with covariance at most \((2n)^{-1}I_d\).
The last two terms have norm at most \((1+k^2)/n=o(k^{-3})\).
Gaussian covariance concentration and Gaussian vector tails
therefore yield, uniformly in these configurations,
\begin{equation}\label{eq:two-sharp-matrix-error}
 \mathbb P\{\|W-A_K\|_{\mathrm{op}}>k^{-3}\mid\boldsymbol j\}
 \le C_de^{-c_dn/k^6}=o\{\lambda_kQ_d(k^2)\}.
\end{equation}
On the complementary event, \(A_K\succeq I_d\) gives
\begin{equation}\label{eq:two-sharp-quadratic-comparison}
 \frac{D_{\boldsymbol j}}{1+k^{-3}}\le T_{n,n}^2
 \le\frac{D_{\boldsymbol j}}{1-k^{-3}},\qquad
 D_{\boldsymbol j}=\frac{(Z_1+h)^2}{1+K/k^2}+\sum_{r=2}^dZ_r^2.
\end{equation}
Completing the square in the density ratio
\(\ell_{\boldsymbol j}\) of \(A_K^{-1/2}(Z+he_1)\) to
\(N_d(0,I_d)\) gives
\begin{equation}\label{eq:two-sharp-density-bound}
 \ell_{\boldsymbol j}(z)
 \le\sqrt{1+K/k^2}\exp\left(\frac{k^2h^2}{2K}\right)
 \le\sqrt2e^{K/2}.
\end{equation}

When \(K=1\), \(h=\pm1/k\). By symmetry, the Gaussian density
calculation in \eqref{eq:sharpness-density-ratio}, with shift
parameter \(1/k\), and \eqref{eq:sharpness-radial-concentration} give
\[
 \frac{\mathbb P(D_{\boldsymbol j}\ge k^2)}{Q_d(k^2)}
 \longrightarrow\mathbb E[e^{-U_1^2/2}\cosh(U_1)],
\]
where \(U\) is uniform on \(\mathbb S^{d-1}\).
The squared-threshold shifts in
\eqref{eq:two-sharp-quadratic-comparison} are \(\pm k^{-1}\).
By \eqref{eq:two-sharp-density-bound}, they change the tail by
\(O_d(k^{-1}Q_d(k^2))\). Together with
\eqref{eq:two-sharp-matrix-error}, this proves, uniformly in
the position and sign of the single shift,
\begin{equation}\label{eq:two-sharp-one-shift}
 \frac{\mathbb P(T_{n,n}^2\ge k^2\mid K=1)}{Q_d(k^2)}
 \longrightarrow\mathbb E[e^{-U_1^2/2}\cosh(U_1)].
\end{equation}

For \(2\le K\le k^2\), the same comparison and density bound give
\[
 \mathbb P(T_{n,n}^2\ge k^2\mid\boldsymbol j)
 \le C_de^{K/2}Q_d(k^2)+C_de^{-c_dn/k^6}.
\]
The binomial generating function and Chernoff bound imply
\[
 \mathbb E[e^{K/2}\mathbf1_{\{K\ge2\}}]
 \le\frac e2\lambda_k^2e^{(\sqrt e-1)\lambda_k}=O(\lambda_k^2),
 \qquad
 \mathbb P(K>k^2)\le(e\lambda_k/k^2)^{k^2}
 =o\{\lambda_kQ_d(k^2)\}.
\]
Consequently,
\begin{equation}\label{eq:two-sharp-multiple-shifts}
 \frac{\mathbb P(T_{n,n}^2\ge k^2,\ K\ge2)}{Q_d(k^2)}
 =o(\lambda_k).
\end{equation}
Substituting \eqref{eq:two-sharp-zero-shifts},
\eqref{eq:two-sharp-one-shift}, and
\eqref{eq:two-sharp-multiple-shifts} into
\eqref{eq:two-sharp-tail-decomposition} yields
\[
 \frac{\mathbb P(T_{n_k,n_k}^2\ge k^2)}{Q_d(k^2)}
 =1-\left\{1-\mathbb E[e^{-U_1^2/2}\cosh(U_1)]\right\}
       \lambda_k+o(\lambda_k).
\]
The coefficient in braces is positive by the preceding example.
Since \((1+k)^pn_k^{-(p-2)/2}\sim\lambda_k\to0\), this proves
\eqref{eq:two-sample-sharp-rate}.
\end{proof}
\endgroup

\section{Auxiliary proofs}
\label{sec:sec4-lemma-proofs}

This section proves the preliminary tail bound, the analytic estimates
for the functions $f_H$ and $F$, and the conditional estimates used in the proof
of \eqref{eq:CMD-Rn-1}.


\begin{proof}[Proof of Lemma~\ref{lem:sec4-field-global}]
For \eqref{eq:sec4-truncated-field}, put
$y=(s^\top K(v)^{-1}s)^{1/2}$. Direct differentiation and the
spectral bounds for $K(v)$ give
\begin{equation}\label{eq:sec4-field-radial-control}
\begin{aligned}
|F(s,v)|&\le C_dy|a_H(y)|,\\
|s^\top F(s,v)|&=y^2|a_H(y)|,\\
\|\partial_sF(s,v)\|_{\mathrm F}
&\le C_d\{|a_H(y)|+y|a_H'(y)|\}.
\end{aligned}
\end{equation}
The formulas at zero are understood by continuity. We first bound
$a_H(y)$ and $ya_H'(y)$ for every $y\ge0$.
For any smooth $H$ with values in $[0,1]$, the definition of the
chi-distribution expectation gives
\[
\int_0^\infty\{H(u)-\mathbb EH(\chi_d)\}
 u^{d-1}e^{-u^2/2}\,du=0.
\]
The integral defining $a_H(y)$ can therefore also be written as the
negative integral over $[y,\infty)$. For $y\le1$, use the original
integral and $|H-\mathbb EH(\chi_d)|\le1$ to bound $|a_H(y)|$ by
$C_d$. For $y\ge1$, use the tail representation and
\[
\int_y^\infty u^{d-1}e^{-u^2/2}\,du
\le C_dy^{d-2}e^{-y^2/2}
\]
to bound it by $C_dy^{-2}$. Hence
\begin{equation}\label{eq:sec4-field-a-bound}
|a_H(y)|\le\frac{C_d}{1+y^2}.
\end{equation}
The radial equation then gives
\begin{equation}\label{eq:sec4-field-radial-derivative}
|ya_H'(y)|
\le|H(y)-\mathbb EH(\chi_d)|+(d+y^2)|a_H(y)|
\le C_d.
\end{equation}
Substituting \eqref{eq:sec4-field-a-bound} and
\eqref{eq:sec4-field-radial-derivative} into
\eqref{eq:sec4-field-radial-control} proves
\eqref{eq:sec4-truncated-field}, uniformly over all such $H$.

For \(v\) with eigenvalues in \([1/2,3/2]\),
\[
\begin{aligned}
f_H(s,v)&=v^{-1/2}F(v^{-1/2}s,I_d),\\
s^\top f_H(s,v)&=(v^{-1/2}s)^\top F(v^{-1/2}s,I_d),\\
\partial_sf_H(s,v)
&=v^{-1/2}\partial_sF(v^{-1/2}s,I_d)v^{-1/2}.
\end{aligned}
\]
By \eqref{eq:sec4-truncated-field},
\begin{equation}\label{eq:sec4-field-spectral-bounds}
|f_H(s,v)|+|s^\top f_H(s,v)|
 +\|\partial_sf_H(s,v)\|_{\mathrm F}\le C_d.
\end{equation}

For \eqref{eq:sec4-global-matrix-derivative}, differentiate
\eqref{eq:sec4-stein-field} at a positive definite matrix \(w\), with
\(y=(s^\top w^{-1}s)^{1/2}>0\), in a symmetric direction \(A\):
\begin{equation}\label{eq:sec4-stein-matrix-derivative}
\partial_v f_H(s,w)[A]
=-w^{-1}Aw^{-1}s\,a_H(y)
-\frac{a_H'(y)}{2y}w^{-1}s
 \bigl(s^\top w^{-1}Aw^{-1}s\bigr).
\end{equation}
For eigenvalues of \(w\) in \([1/2,3/2]\),
\eqref{eq:sec4-field-a-bound} and
\eqref{eq:sec4-field-radial-derivative} imply
\[
\|\partial_v f_H(s,w)\|
\le C_d\{y|a_H(y)|+y^2|a_H'(y)|\}
\le C_d(1+|s|).
\]
At \(s=0\) the derivative vanishes. The cutoff definition
\eqref{eq:sec4-extension} gives \(\|DK(v)[A]\|_{\mathrm F}\le C_d\|A\|_{\mathrm F}\),
since the cutoff \(H_{1/16,\,3/16}\) and its derivative vanish
when \(\|v-I_d\|_{\mathrm F}>1/2\).
The chain rule proves \eqref{eq:sec4-global-matrix-derivative}.

We finally prove \eqref{eq:sec4-global-jacobian-modulus}.
Use $H=H_{u,\eta}$ from \eqref{eq:sec4-cutoff}.
Equation \eqref{eq:sec4-field-spectral-bounds} at $v=I_d$ gives
$\|\partial_sf_H(w,I_d)\|_{\mathrm F}\le C_d$.
We first refine the derivative bound near zero.
For $0<y\le1$, subtracting $H(y)-\mathbb EH(\chi_d)$ in the
integral defining $a_H$ yields
\begin{equation}\label{eq:sec4-field-origin-expansion}
\begin{aligned}
d\,a_H(y)-\{H(y)-\mathbb EH(\chi_d)\}
&=d\,y^{-d}e^{y^2/2}\int_0^y
 \{H(t)-H(y)\}t^{d-1}e^{-t^2/2}\,dt\\
&\quad+\{H(y)-\mathbb EH(\chi_d)\}
 \left[d\,y^{-d}e^{y^2/2}
             \int_0^y t^{d-1}e^{-t^2/2}\,dt-1\right].
\end{aligned}
\end{equation}
The first term in \eqref{eq:sec4-field-origin-expansion} has absolute
value at most $C_dy/\eta$, by the Lipschitz bound for $H$.
Substituting $t=yv$ rewrites the expression in square brackets in
\eqref{eq:sec4-field-origin-expansion} as
\[
d\int_0^1v^{d-1}\bigl[\exp\{y^2(1-v^2)/2\}-1\bigr]\,dv
=O(y^2).
\]
Using \eqref{eq:sec4-radial-equation} and
\eqref{eq:sec4-field-a-bound}, we obtain
\begin{equation}\label{eq:sec4-field-radial-small}
|ya_H'(y)|\le C_d\min\{1,y/\eta+y^2\},
\qquad 0<y\le1.
\end{equation}
For $y\le\eta$, divide the bound $C_d(y/\eta+y^2)$ in
\eqref{eq:sec4-field-radial-small} by $y$. For $y>\eta$, including
$y\ge1$, use \eqref{eq:sec4-field-radial-derivative}. We obtain
\begin{equation}\label{eq:sec4-field-a-prime}
|a_H'(y)|\le C_d\{1+(y+\eta)^{-1}\},\qquad y>0.
\end{equation}

For a nonzero vector $w$, write $y=|w|$ and $e=w/y$.
Differentiating $f_H(s,I_d)$ with respect to $s$ gives
\[
\partial_sf_H(w,I_d)=a_H(y)I_d+ya_H'(y)ee^\top.
\]
The radial differential equation implies
\[
\{ya_H'(y)\}'
=H'(y)+2ya_H(y)-(d-y^2)a_H'(y).
\]
Since $dy=e^\top dw$ and $de=(I_d-ee^\top)dw/y$, we deduce
\begin{equation}\label{eq:sec4-field-jacobian-differential}
\begin{aligned}
&\left\|d\{\partial_sf_H(w,I_d)\}
       -H'(y)\,dy\,ee^\top\right\|_{\mathrm F}\\
&\quad\le C_d\{y|a_H(y)|+(1+y^2)|a_H'(y)|\}|dw|\\
&\quad\le C_d\{1+y+(y+\eta)^{-1}\}|dw|.
\end{aligned}
\end{equation}
The last inequality in \eqref{eq:sec4-field-jacobian-differential}
follows from \eqref{eq:sec4-field-a-prime} when $y\le1$ and
\eqref{eq:sec4-field-radial-derivative} when $y\ge1$.

Now fix $w,w'$ and return to $y=|w'|$ and $\delta=|w-w'|$.
If $\delta\ge(y+\eta)/2$, \eqref{eq:sec4-field-spectral-bounds} at
$v=I_d$ proves \eqref{eq:sec4-global-jacobian-modulus}, since
$\min\{1,\delta/(y+\eta)\}\ge1/2$.
Otherwise, every point on the segment from $w'$ to $w$ has radius
at most $y+\delta$, and its radius plus $\eta$ is at least
$(y+\eta)/2$. Integrating
\eqref{eq:sec4-field-jacobian-differential} along the segment gives
\[
C_d\delta(1+y+\delta)+C_d\frac{\delta}{y+\eta}
\]
for the contribution excluding $H'$.
Since $\delta/(y+\eta)<1/2$, the fraction equals the minimum in
\eqref{eq:sec4-global-jacobian-modulus}.

It remains to integrate the term containing $H'$.
The radius along a straight segment is monotone on at most two
subintervals. As $H$ is non-increasing and takes values in $[0,1]$,
\begin{equation}\label{eq:sec4-field-Hprime-variation}
\int_{\mathrm{segment}}|H'(y)|\,|dy|\le2.
\end{equation}
The derivative bound and the length of the segment also give
\begin{equation}\label{eq:sec4-field-Hprime-lipschitz}
\int_{\mathrm{segment}}|H'(y)|\,|dy|
\le C_d\delta/\eta.
\end{equation}
The integrals in \eqref{eq:sec4-field-Hprime-variation} and
\eqref{eq:sec4-field-Hprime-lipschitz} vanish unless the segment has a radius in
$[u,u+\eta]$, which requires $|y-u|\le\delta+\eta$.
Together, \eqref{eq:sec4-field-Hprime-variation} and
\eqref{eq:sec4-field-Hprime-lipschitz} give the last term in
\eqref{eq:sec4-global-jacobian-modulus}.
If the segment passes through zero, integrate on each side and pass
to the endpoint by continuity of $\partial_sf_H(\cdot,I_d)$.
This proves \eqref{eq:sec4-global-jacobian-modulus}.
\end{proof}


\begin{proof}[Proof of Lemma~\ref{lem:sec4-conditional-moments}]
By \eqref{eq:sec4-conditioning-bound},
\[
\mathbb P(|\boldsymbol\xi_i|\le b)
\ge\mathbb P(\mathcal N)\ge3/4.
\]
Under the product law in \eqref{eq:sec4-conditioning}, each
conditional moment is its truncated moment divided by
$\mathbb P(|\boldsymbol\xi_i|\le b)$.

Since $\mathbb E\boldsymbol\xi_i=0$, the conditional means satisfy
\[
\begin{aligned}
\sum_i|\mathbb E^*\boldsymbol\xi_i|
&=\sum_i\frac{
 |\mathbb E[\boldsymbol\xi_i;|\boldsymbol\xi_i|>b]|}
 {\mathbb P(|\boldsymbol\xi_i|\le b)}\\
&\le\frac43 b^{1-p}\sum_i\mathbb E|\boldsymbol\xi_i|^p
\le2b^{1-p}\sum_i\mathbb E|\boldsymbol\xi_i|^p.
\end{aligned}
\]

For the matrix-bias bound in \eqref{eq:sec4-conditional-moments},
subtract the original second-moment sum from the conditional sum:
\begin{equation}\label{eq:sec4-conditional-matrix-bias-decomposition}
\begin{aligned}
\mathbb E^*V_n-I_d
&=-\sum_{i=1}^n
 \mathbb E[\boldsymbol\xi_i\boldsymbol\xi_i^\top;
                         |\boldsymbol\xi_i|>b]\\
&\quad+\sum_{i=1}^n
 \frac{\mathbb P(|\boldsymbol\xi_i|>b)}
      {\mathbb P(|\boldsymbol\xi_i|\le b)}
 \mathbb E[\boldsymbol\xi_i\boldsymbol\xi_i^\top;
                         |\boldsymbol\xi_i|\le b].
\end{aligned}
\end{equation}
The Frobenius norm of $zz^\top$ is $|z|^2$.
The first sum in \eqref{eq:sec4-conditional-matrix-bias-decomposition}
has norm at most
$b^{2-p}\sum_i\mathbb E|\boldsymbol\xi_i|^p$.
For the second sum in
\eqref{eq:sec4-conditional-matrix-bias-decomposition}, use
$\mathbb P(|\boldsymbol\xi_i|\le b)^{-1}\le2$ and
$\|\mathbb E[\boldsymbol\xi_i\boldsymbol\xi_i^\top;
 |\boldsymbol\xi_i|\le b]\|_{\mathrm F}\le b^2$.
Its norm is therefore at most
\[
2b^2\sum_i\mathbb P(|\boldsymbol\xi_i|>b)
\le2b^{2-p}\sum_i\mathbb E|\boldsymbol\xi_i|^p.
\]
Applying the triangle inequality to
\eqref{eq:sec4-conditional-matrix-bias-decomposition} gives the
coefficient $3$ in \eqref{eq:sec4-conditional-moments}.

For $j=3,4$, we have $j\ge p$ and hence
$|\boldsymbol\xi_i|^j\le b^{j-p}|\boldsymbol\xi_i|^p$
on $\{|\boldsymbol\xi_i|\le b\}$.
The conditional expectation formula gives
\[
\sum_i\mathbb E^*|\boldsymbol\xi_i|^j
\le\frac43 b^{j-p}\sum_i\mathbb E|\boldsymbol\xi_i|^p
\le2b^{j-p}\sum_i\mathbb E|\boldsymbol\xi_i|^p.
\]
Finally, taking traces in the covariance normalization gives
$\sum_i\mathbb E|\boldsymbol\xi_i|^2=d$.
Consequently,
\[
\sum_i\mathbb E^*|\boldsymbol\xi_i|^2
=\sum_i\frac{\mathbb E[|\boldsymbol\xi_i|^2;
                         |\boldsymbol\xi_i|\le b]}
                 {\mathbb P(|\boldsymbol\xi_i|\le b)}
\le\frac43\sum_i\mathbb E|\boldsymbol\xi_i|^2
=\frac{4d}{3}.
\]
This proves \eqref{eq:sec4-conditional-moments}.
\end{proof}

\begin{proof}[Proof of Lemma~\ref{lem:sec4-small-increments}]
For \eqref{eq:sec4-small-path}, we compare both path points with
the deleted point. We first give the required matrix estimates.
For $s,z\in\mathbb R^d$, let $v=v^\top$ satisfy $\|v-I_d\|_{\mathrm F}\le1/8$, and let
$|z|\le\varepsilon/4$. Put $y=|v^{-1/2}s|$.
For $0\le t\le1$,
\begin{equation}\label{eq:sec4-deletion-path-spectrum}
\|v+tzz^\top-I_d\|_{\mathrm F}
\le\frac18+\frac{\varepsilon^2}{16}\le\frac3{16}<\frac14.
\end{equation}
Thus the path matrices have eigenvalues in $[3/4,5/4]$, and
\eqref{eq:sec4-extension} gives
$K(v+tzz^\top)=v+tzz^\top$.
The inverse square-root map is Lipschitz on this spectral interval:
\begin{equation}\label{eq:sec4-deletion-inverse-root}
\left\|(v+tzz^\top)^{-1/2}-v^{-1/2}\right\|_{\mathrm F}
\le C_dt|z|^2.
\end{equation}
Using $|s|\le\sqrt{9/8}\,y$ and
\eqref{eq:sec4-deletion-inverse-root}, we obtain
\begin{equation}\label{eq:sec4-generic-path-distance}
\begin{aligned}
w_t&=(v+tzz^\top)^{-1/2}(s+tz),\qquad0\le t\le1,\\
|w_0|&=y,\qquad
\sup_{0\le t\le1}|w_t-w_0|
\le C_d\bigl(|z|+|z|^2y\bigr).
\end{aligned}
\end{equation}

Since \(b\le\varepsilon/4\), on
\(\mathcal G_i\cap\{|\boldsymbol\xi_i|\le b\}\) we may substitute
\((s,v,z)=(S_{-i},V_{-i},\boldsymbol\xi_i)\) into
\eqref{eq:sec4-deletion-path-spectrum} and
\eqref{eq:sec4-generic-path-distance} to obtain
\begin{equation}\label{eq:sec4-deletion-path-distance}
\begin{aligned}
w_{i,t}
&=(V_{-i}+t\boldsymbol\xi_i\boldsymbol\xi_i^\top)^{-1/2}
  (S_{-i}+t\boldsymbol\xi_i),\qquad0\le t\le1,\\
|w_{i,0}|&=R_{-i},\\
\sup_{0\le t\le1}|w_{i,t}-w_{i,0}|
&\le C_d\bigl(|\boldsymbol\xi_i|+|\boldsymbol\xi_i|^2R_{-i}\bigr).
\end{aligned}
\end{equation}
The path matrices have eigenvalues in $[3/4,5/4]$, and
$K(V_{-i}+t\boldsymbol\xi_i\boldsymbol\xi_i^\top)
 =V_{-i}+t\boldsymbol\xi_i\boldsymbol\xi_i^\top$.

Along the path,
\begin{equation}\label{eq:sec4-deletion-jacobian-conjugation}
\begin{aligned}
&\partial_sF\left(S_{-i}+t\boldsymbol\xi_i,
                 V_{-i}+t\boldsymbol\xi_i\boldsymbol\xi_i^\top\right)\\
&\quad=(V_{-i}+t\boldsymbol\xi_i\boldsymbol\xi_i^\top)^{-1/2}
 \partial_sf_H(w_{i,t},I_d)
 (V_{-i}+t\boldsymbol\xi_i\boldsymbol\xi_i^\top)^{-1/2}.
\end{aligned}
\end{equation}
Using \eqref{eq:sec4-deletion-inverse-root},
\eqref{eq:sec4-deletion-jacobian-conjugation}, and
\eqref{eq:sec4-truncated-field} at $v=I_d$ gives
\begin{equation}\label{eq:sec4-deletion-jacobian-reduction}
\begin{aligned}
&\left\|\partial_sF(S_n,V_n)
 -\partial_sF\left(S_{-i}+t\boldsymbol\xi_i,
                  V_{-i}+t\boldsymbol\xi_i\boldsymbol\xi_i^\top\right)
 \right\|_{\mathrm F}\\
&\quad\le C_d|\boldsymbol\xi_i|^2
 +C_d\sum_{\tau\in\{t,1\}}
 \|\partial_sf_H(w_{i,\tau},I_d)
       -\partial_sf_H(w_{i,0},I_d)\|_{\mathrm F}.
\end{aligned}
\end{equation}
Apply \eqref{eq:sec4-global-jacobian-modulus} with
$w'=w_{i,0}$ and $w=w_{i,\tau}$, for $\tau\in\{t,1\}$.
Since $|\boldsymbol\xi_i|\le b\le1$,
\eqref{eq:sec4-deletion-path-distance} gives
\begin{equation}\label{eq:sec4-deletion-polynomial-modulus}
\begin{aligned}
&|\boldsymbol\xi_i|^2+|w_{i,\tau}-w_{i,0}|
 (1+R_{-i}+|w_{i,\tau}-w_{i,0}|)\\
&\quad\le C_d|\boldsymbol\xi_i|(1+R_{-i}^2),
\end{aligned}
\end{equation}
and
\begin{equation}\label{eq:sec4-deletion-inverse-modulus}
\min\left\{1,\frac{|w_{i,\tau}-w_{i,0}|}{R_{-i}+\eta}\right\}
\le C_d|\boldsymbol\xi_i|^2
 +\min\left\{1,\frac{C_d|\boldsymbol\xi_i|}{R_{-i}+\eta}\right\}.
\end{equation}
The first term on the right of \eqref{eq:sec4-deletion-inverse-modulus}
is bounded by $C_d|\boldsymbol\xi_i|(1+R_{-i}^2)$.
For the indicator in \eqref{eq:sec4-global-jacobian-modulus},
the event $|R_{-i}-u|\le|w_{i,\tau}-w_{i,0}|+\eta$ implies
\begin{equation}\label{eq:sec4-deletion-shell-width}
\begin{aligned}
R_{-i}
&\le\frac{u+\eta+C_d|\boldsymbol\xi_i|}
              {1-C_d|\boldsymbol\xi_i|^2}
\le C_d(1+u),\\
|w_{i,\tau}-w_{i,0}|
&\le C_d|\boldsymbol\xi_i|(1+u).
\end{aligned}
\end{equation}
The denominator in \eqref{eq:sec4-deletion-shell-width} is positive
for the common sufficiently small $\varepsilon$.
Equation \eqref{eq:sec4-deletion-shell-width} bounds the indicator term
in \eqref{eq:sec4-global-jacobian-modulus} by the last term of
\eqref{eq:sec4-small-path}.
Substituting \eqref{eq:sec4-deletion-polynomial-modulus},
\eqref{eq:sec4-deletion-inverse-modulus}, and
\eqref{eq:sec4-deletion-shell-width} in
\eqref{eq:sec4-deletion-jacobian-reduction} proves
\eqref{eq:sec4-small-path}, uniformly in $t$.
\end{proof}

\begin{proof}[Proof of Lemma~\ref{lem:sec4-deleted-expectations}]
For \eqref{eq:sec4-small-inverse-contribution} and
\eqref{eq:sec4-small-shell-contribution}, we first prove the
distribution comparison \eqref{eq:sec4-small-deleted-cdf}.
The leave-one-out sample under \(\mathbb P\) has covariance
\(C_i=I_d-\mathbb E(\boldsymbol\xi_i\boldsymbol\xi_i^\top)
\succeq(1-\lambda^{2/p})I_d\).
For \(0<\lambda\le\lambda_0\),
\begin{equation}\label{eq:sec4-small-deleted-normalization}
\begin{aligned}
\|C_i^{-1/2}\|_{\mathrm{op}}^p
&\le(1-\lambda^{2/p})^{-p/2},\\
\sum_{j\ne i}\mathbb E|C_i^{-1/2}\boldsymbol\xi_j|^p
&\le(1-\lambda^{2/p})^{-p/2}
       \sum_{j\ne i}\mathbb E|\boldsymbol\xi_j|^p
\le2\lambda.
\end{aligned}
\end{equation}
The vectors \(C_i^{-1/2}\boldsymbol\xi_j\), \(j\ne i\), are centered
and independent, with total covariance \(I_d\).
The normalization leaves \(R_{-i}\) unchanged, including when
\(V_{-i}\) is singular. Indeed, for arbitrary vectors \(z_j\), put
\(s=\sum_jz_j\) and \(v=\sum_jz_jz_j^\top\).
Since \(s\) belongs to the range of \(v\), completing the square gives
\begin{equation}\label{eq:sec4-prob-variational}
\begin{aligned}
s^\top v^+s
&=\sup_{\theta\in\mathbb R^d}
 \{2\theta^\top s-\theta^\top v\theta\}\\
&=\sup_{\theta\in\mathbb R^d}
 \sum_j\{2\theta^\top z_j-(\theta^\top z_j)^2\}.
\end{aligned}
\end{equation}
An invertible linear transformation of the observations in
\eqref{eq:sec4-prob-variational} is absorbed by \(\theta\).
By \eqref{eq:sec4-small-error} and
\eqref{eq:sec4-small-deleted-normalization}, the error under
\(\mathbb P\) is at most \((1+t)\mathcal B(2\lambda)\).
Under \(\mathbb P^*\), the law of \(R_{-i}\) is its original law
conditioned on \(\bigcap_{j\ne i}\{|\boldsymbol\xi_j|\le b\}\).
Since $b=\varepsilon/4$ is fixed, this event has complement of probability at most
\(b^{-p}\sum_{j\ne i}\mathbb E|\boldsymbol\xi_j|^p\le C_d\lambda\).
Consequently,
\begin{equation}\label{eq:sec4-small-deleted-cdf}
|\mathbb P^*(R_{-i}\le t)-\mathbb P(\chi_d\le t)|
\le(1+t)\mathcal B(2\lambda)+C_d\lambda,\qquad t\ge0.
\end{equation}
Equation \eqref{eq:sec4-small-deleted-cdf} also holds with
\(R_{-i}<t\) and \(\chi_d<t\). For \(t>0\), take limits from below.
At \(t=0\), both probabilities are zero.

Under $\mathbb P^*$, $R_{-i}$ is independent of
$\boldsymbol\xi_i$. Thus \eqref{eq:sec4-small-deleted-cdf} can also
be applied conditional on $\boldsymbol\xi_i$.
For \eqref{eq:sec4-small-inverse-contribution}, fix
\(\boldsymbol\xi_i\). The function
\[
y\longmapsto
\min\{1,C_d|\boldsymbol\xi_i|/(y+\eta)\}
\mathbf1_{\{y\le1\}}
\]
is non-increasing on \([0,\infty)\) and has total decrease
\(\min\{1,C_d|\boldsymbol\xi_i|/\eta\}\).
Stieltjes integration of \eqref{eq:sec4-small-deleted-cdf}, conditional
on \(\boldsymbol\xi_i\), gives
\begin{equation}\label{eq:sec4-small-inverse-conditional}
\begin{aligned}
&\mathbb E^*\left[
 \min\left\{1,\frac{C_d|\boldsymbol\xi_i|}{R_{-i}+\eta}\right\}
 \mathbf1_{\{R_{-i}\le1\}}\,\middle|\,\boldsymbol\xi_i\right]\\
&\quad\le C_d|\boldsymbol\xi_i|
 +C_d\{\mathcal B(2\lambda)+\lambda\}
 \min\{1,C_d|\boldsymbol\xi_i|/\eta\}.
\end{aligned}
\end{equation}
The chi term in \eqref{eq:sec4-small-inverse-conditional} is at most
\(C_d|\boldsymbol\xi_i|\mathbb E\chi_d^{-1}\), with
\(\mathbb E\chi_d^{-1}<\infty\) for \(d\ge2\).
Multiplying \eqref{eq:sec4-small-inverse-conditional} by
\(|\boldsymbol\xi_i|^2\), taking expectations, and using
\eqref{eq:sec4-small-conditional-moments} proves
\eqref{eq:sec4-small-inverse-contribution}.

For \eqref{eq:sec4-small-shell-contribution}, write
\begin{equation}\label{eq:sec4-interval-cdf}
\begin{aligned}
&\mathbb P^*(|R_{-i}-u|\le C_d|\boldsymbol\xi_i|(1+u)+\eta
             \mid\boldsymbol\xi_i)\\
&\quad=\mathbb P^*(R_{-i}\le u+C_d|\boldsymbol\xi_i|(1+u)+\eta
                    \mid\boldsymbol\xi_i)\\
&\qquad-\mathbb P^*(R_{-i}<(u-C_d|\boldsymbol\xi_i|(1+u)-\eta)_+
                      \mid\boldsymbol\xi_i).
\end{aligned}
\end{equation}
Equation \eqref{eq:sec4-interval-cdf} also holds when the interval
contains zero. At both endpoints in \eqref{eq:sec4-interval-cdf},
\(1+t\le C_d(1+u)\), since \(|\boldsymbol\xi_i|\le b\) and \(\eta\le1\).
Apply \eqref{eq:sec4-small-deleted-cdf} for \(R_{-i}\le t\) and
\(R_{-i}<t\) in \eqref{eq:sec4-interval-cdf}.
The bounded chi density gives
\begin{equation}\label{eq:sec4-small-interval-probability}
\begin{aligned}
&\mathbb P^*(|R_{-i}-u|\le C_d|\boldsymbol\xi_i|(1+u)+\eta
             \mid\boldsymbol\xi_i)\\
&\quad\le C_d\{|\boldsymbol\xi_i|(1+u)+\eta\}
 +C_d(1+u)\mathcal B(2\lambda)+C_d\lambda.
\end{aligned}
\end{equation}
For \(w\ge0\),
\begin{equation}\label{eq:sec4-shell-weight-inequalities}
\min\{1,w/\eta\}(w+\eta)\le2w,\qquad
\min\{1,w/\eta\}\le w/\eta.
\end{equation}
By \eqref{eq:sec4-small-conditional-moments},
\begin{equation}\label{eq:sec4-small-shell-moments}
(1+u)\sum_i\mathbb E^*|\boldsymbol\xi_i|^3
\le C_d\lambda(1+u).
\end{equation}
Apply \eqref{eq:sec4-small-interval-probability}, then
\eqref{eq:sec4-shell-weight-inequalities} and
\eqref{eq:sec4-small-shell-moments} to obtain \eqref{eq:sec4-small-shell-contribution}.
\end{proof}


\begin{proof}[Proof of Lemma~\ref{lem:sec4-solution}]
\emph{(i).}
To prove \eqref{eq:sec4-stein-field-envelope}, first bound
\(f_H\) and \(\partial_vf_H\) when the eigenvalues of \(v\) lie in
\([1/2,3/2]\). For \(y=(s^\top v^{-1}s)^{1/2}\),
\eqref{eq:sec4-stein-field} and
\eqref{eq:sec4-stein-matrix-derivative} give
\begin{equation}\label{eq:sec4-stein-field-radial-control}
|f_H(s,v)|\le C_dy|a_H(y)|,\qquad
\|\partial_vf_H(s,v)\|
 \le C_d\{y|a_H(y)|+y^2|a_H'(y)|\}.
\end{equation}
The value at \(y=0\) follows by continuity.

For \(H=H_{u,\eta}\), \(u\in\{x-\eta,x\}\),
\begin{equation}\label{eq:sec4-stein-centering-tail}
0\le1-\mathbb EH(\chi_d)\le\mathbb P(\chi_d\ge x-\eta)
\le C_d\mathbb P(\chi_d^2\ge x^2).
\end{equation}
The tail comparison in \eqref{eq:sec4-stein-centering-tail}
uses \(x\eta\le1/10\).
For $0\le y\le x-x^{-1}$, we have $H=1$ on $[0,y]$.
The substitution $u=yt$ therefore gives
\[
a_H(y)=(1-\mathbb EH(\chi_d))\int_0^1
 t^{d-1}\exp\{y^2(1-t^2)/2\}\,dt.
\]
Differentiating under the integral and using
$\int_0^1t^{d-1}e^{-y^2t^2/2}\,dt\le C_d(1+y)^{-d}$ yields
\begin{equation}\label{eq:sec4-stein-radial-below}
\begin{aligned}
|a_H(y)|&\le C_d\mathbb P(\chi_d^2\ge x^2)
 e^{y^2/2}(1+y)^{-d},\\
|a_H'(y)|&\le C_d\mathbb P(\chi_d^2\ge x^2)
 y e^{y^2/2}(1+y)^{-d},
\qquad 0\le y\le x-x^{-1}.
\end{aligned}
\end{equation}
For $y\ge x+\eta$, the negative-integral representation and the
substitution $u^2=y^2+2t$ give
\[
a_H(y)=\frac{\mathbb EH(\chi_d)}{y^2}\int_0^\infty
 \left(1+\frac{2t}{y^2}\right)^{(d-2)/2}e^{-t}\,dt.
\]
The integral is bounded, and its derivative is \(O(y^{-3})\),
uniformly for \(y\ge x+\eta\ge4\). Thus
\begin{equation}\label{eq:sec4-stein-radial-above}
|a_H(y)|+y|a_H'(y)|\le C_dy^{-2},\qquad y\ge x+\eta.
\end{equation}
By \eqref{eq:sec4-field-a-bound} and
\eqref{eq:sec4-field-radial-derivative},
\begin{equation}\label{eq:sec4-stein-radial-middle}
|a_H(y)|\le C_dx^{-2},\qquad |ya_H'(y)|\le C_d,
\qquad x-x^{-1}<y<x+\eta.
\end{equation}
Substituting \eqref{eq:sec4-stein-radial-below},
\eqref{eq:sec4-stein-radial-middle}, and
\eqref{eq:sec4-stein-radial-above} into
\eqref{eq:sec4-stein-field-radial-control} gives
\begin{equation}\label{eq:sec4-solution-field-envelope}
|f_H(s,v)|\le C_d A_{1,x}(y),\qquad
\|\partial_vf_H(s,v)\|\le C_d A_{3,x}(y).
\end{equation}
The chi-squared density gives
\begin{equation}\label{eq:sec4-stein-chi-tail}
c_dx^{d-2}e^{-x^2/2}
\le\mathbb P(\chi_d^2\ge x^2)
\le C_dx^{d-2}e^{-x^2/2},\qquad x\ge4.
\end{equation}
Indeed, substitute $u^2=x^2+2t$ in its tail integral; the remaining
integral is bounded above and below by positive constants depending
only on $d$. For $x/2\le y\le x-x^{-1}$,
\eqref{eq:sec4-stein-chi-tail} directly gives
$A_{3,x}(y)\le C_d(1+x)$. For $y<x/2$, the factor
$e^{-(x^2-y^2)/2}$ absorbs the polynomial factors in $x$ and $y$.
The second branch of $A_{3,x}$ equals $1+x$. Hence
$\sup_{y\ge0}A_{3,x}(y)\le C_d(1+x)$.
Apply \eqref{eq:sec4-solution-field-envelope} at \(K(v)\).
Since \(\|DK(v)\|\le C_d\) by \eqref{eq:sec4-extension},
the chain rule gives
\eqref{eq:sec4-stein-field-envelope}.

\emph{(ii).}
We first prove \eqref{eq:sec4-stein-path-matrix} by comparing
\(A_{3,x}\) along the path.
Use \(y=(s^\top v^{-1}s)^{1/2}\) and \(w_t\) from
\eqref{eq:sec4-generic-path-distance}. Since
\(|z|\le\varepsilon/x\le\varepsilon/4\),
\eqref{eq:sec4-deletion-path-spectrum} and
\eqref{eq:sec4-deletion-inverse-root} apply.

\noindent\emph{The matrix derivative along the path.}
For $y\le2x$, \eqref{eq:sec4-generic-path-distance} gives
\[
\sup_{0\le t\le1}|w_t-w_0|\le C_d|z|.
\]
On $0\le q\le2x+1$,
\[
\left|\frac{d}{dq}\log\{e^{q^2/2}(1+q)^{3-d}\}\right|
\le2x+1+|3-d|.
\]
Since $x|z|\le\varepsilon$, two radii in $[0,2x+1]$ that differ
by at most $C_d|z|$ have first-branch values of $A_{3,x}$ within a
fixed factor.
At the branch boundary, \eqref{eq:sec4-stein-chi-tail} gives
\begin{equation}\label{eq:sec4-stein-envelope-boundary}
c_d(1+x)\le
\mathbb P(\chi_d^2\ge x^2)
 e^{(x-x^{-1})^2/2}(1+x-x^{-1})^{3-d}
\le C_d(1+x).
\end{equation}
Hence \eqref{eq:sec4-stein-envelope-boundary} also compares the
two branches when an interval crosses $x-x^{-1}$.
If $y>2x$, the fixed choice of $\varepsilon$ and
\eqref{eq:sec4-generic-path-distance} give
$\sup_t|w_t-w_0|\le y/4$.
Every path radius then exceeds $3y/4>x+\eta$, where $A_{3,x}$ is
constant. Therefore,
\begin{equation}\label{eq:sec4-stein-path-envelope-comparison}
\sup_{0\le t\le1}A_{3,x}(|w_t|)\le C_dA_{3,x}(y).
\end{equation}
Since $K(v+tzz^\top)=v+tzz^\top$, the radius in
\eqref{eq:sec4-stein-field-envelope} equals $|w_t|$.
Combining \eqref{eq:sec4-stein-field-envelope} and
\eqref{eq:sec4-stein-path-envelope-comparison} proves
\eqref{eq:sec4-stein-path-matrix}.

\noindent\emph{The change in the Jacobian.}
To prove \eqref{eq:sec4-stein-path-jacobian}, put
$J_H(w)=\partial_sf_H(w,I_d)$.
Equations \eqref{eq:sec4-stein-radial-below},
\eqref{eq:sec4-field-a-bound},
\eqref{eq:sec4-field-radial-derivative}, and
\eqref{eq:sec4-stein-radial-above} give
\begin{equation}\label{eq:sec4-stein-jacobian-envelope}
\|J_H(w)\|_{\mathrm F}
+|w|\,|a_H(|w|)|+(1+|w|^2)|a_H'(|w|)|
\le C_dA_{3,x}(|w|).
\end{equation}
For $w\ne0$, let $e=w/|w|$. We use the first inequality in
\eqref{eq:sec4-field-jacobian-differential}, with $|w|$ as the
radial variable. At zero, \eqref{eq:sec4-stein-radial-below} gives
$a_H'(q)=O(q)$, so the differential estimate extends through zero.
For $|w|>x+\eta$, \eqref{eq:sec4-stein-radial-above} also gives
\begin{equation}\label{eq:sec4-stein-jacobian-tail}
\|J_H(w)\|_{\mathrm F}\le C_d(|w|)^{-2},\qquad
|w|\,|a_H(|w|)|+(1+|w|^2)|a_H'(|w|)|
\le C_d(|w|)^{-1}.
\end{equation}

Fix $t\in[0,1]$ and join $w_t$ to $w_1$ by a straight segment.
If $y\le2x$, the segment has length at most $C_d|z|$ and lies
within $C_d|z|$ of $w_0$. The comparison used for
\eqref{eq:sec4-stein-path-envelope-comparison} holds throughout
the segment. By \eqref{eq:sec4-field-jacobian-differential} and
\eqref{eq:sec4-stein-jacobian-envelope},
\begin{equation}\label{eq:sec4-stein-jacobian-slow-near}
\int_{\mathrm{segment}}
 \|dJ_H(w)-H'(|w|)\,d(|w|)\,ee^\top\|_{\mathrm F}
\le C_d|z|A_{3,x}(y).
\end{equation}
The radius along a straight segment is monotone on at most two
subintervals. Since $H$ is non-increasing and takes values in $[0,1]$,
\begin{equation}\label{eq:sec4-stein-jacobian-shell-variation}
\int_{\mathrm{segment}}|H'(|w|)|\,|d(|w|)|\le2,\qquad
\int_{\mathrm{segment}}|H'(|w|)|\,|d(|w|)|
\le\frac{C_d|z|}{\eta}.
\end{equation}
The integrals in \eqref{eq:sec4-stein-jacobian-shell-variation}
vanish unless the segment meets $(x-\eta,x+\eta)$ in radius.
In that case, $|y-x|<C_d|z|+\eta$.
Thus the contribution containing $H'$ is at most
\[
C_d\min\{1,C_d|z|/\eta\}
 \mathbf1_{\{|y-x|<C_d|z|+\eta\}}.
\]

If $y>2x$, every point on the segment has radius at least
$3y/4>x+\eta$, so $H'=0$ there.
Using \eqref{eq:sec4-stein-jacobian-tail} and
\eqref{eq:sec4-generic-path-distance}, we obtain
\begin{equation}\label{eq:sec4-stein-jacobian-slow-far}
\int_{\mathrm{segment}}\|dJ_H(w)\|_{\mathrm F}
\le\frac{C_d}{y}\bigl(|z|+|z|^2y\bigr)
\le C_d|z|\le C_d|z|A_{3,x}(y).
\end{equation}

Finally, since \(K(v+tzz^\top)=v+tzz^\top\),
\[
\partial_sF(s+tz,v+tzz^\top)
=(v+tzz^\top)^{-1/2}J_H(w_t)(v+tzz^\top)^{-1/2}.
\]
By \eqref{eq:sec4-deletion-inverse-root}, both inverse square-root
factors are bounded and change by at most $C_d|z|^2$.
Equations \eqref{eq:sec4-stein-jacobian-envelope} and
\eqref{eq:sec4-stein-path-envelope-comparison} bound $J_H(w_t)$ by
$C_dA_{3,x}(y)$. The change in the exterior factors is at most
$C_d|z|^2A_{3,x}(y)\le C_d|z|A_{3,x}(y)$.
Combining \eqref{eq:sec4-stein-jacobian-slow-near},
\eqref{eq:sec4-stein-jacobian-shell-variation}, and
\eqref{eq:sec4-stein-jacobian-slow-far} proves
\eqref{eq:sec4-stein-path-jacobian}, uniformly in $t$.
\end{proof}


\begin{proof}[Proof of Lemma~\ref{lem:sec4-preliminary}]
For any event \(A\) and \(\mathbb P(\mathcal N)>0\),
\begin{equation}\label{eq:sec4-conditioning-transfer-formula}
\left|\mathbb P(A\mid\mathcal N)-\mathbb P(A)\right|
\le\frac{\mathbb P(A\cap\mathcal N^c)
+\mathbb P(\mathcal N^c)\mathbb P(A)}{\mathbb P(\mathcal N)}.
\end{equation}
For \(A=\{R_n^2\ge x^2\}\), we need a lower bound for
\(\mathbb P(\mathcal N)\), the coarse bound
\eqref{eq:sec4-prelim-coarse}, and a bound for
\(\mathbb P(R_n^2\ge x^2,\mathcal N^c)\).
Let \(x\ge4\), and use the fixed \(\varepsilon\), \(\lambda_*\), and
cutoff \(b\) from the main proof. The union bound and
\eqref{eq:sec4-prelim-budget} give
\begin{equation}\label{eq:sec4-prelim-exceedance}
\sum_{i=1}^n\mathbb P\{|\boldsymbol\xi_i|>b\}
\le b^{-p}\sum_{i=1}^n\mathbb E|\boldsymbol\xi_i|^p
\le\varepsilon^{-p}\lambda\le\frac14.
\end{equation}
For \(J\subseteq\{1,\ldots,n\}\), write \(k=|J|\) and put
\begin{equation}\label{eq:sec4-prelim-deleted-sample}
\mathcal N_J=\bigcap_{i\notin J}\{|\boldsymbol\xi_i|\le b\},
\qquad
R_{-J}^2=\sup_{\theta\in\mathbb R^d}
\sum_{i\notin J}\{2\theta^\top\boldsymbol\xi_i
-(\theta^\top\boldsymbol\xi_i)^2\}.
\end{equation}
For \(J=\varnothing\), \(\mathcal N_J=\mathcal N\) and
\(R_{-J}^2=R_n^2\).
Write \(R_{-J}=(R_{-J}^2)^{1/2}\).
Equation \eqref{eq:sec4-prob-variational} identifies this statistic
with the deleted-sample statistic even when the empirical matrix is
singular. Since \(2a-a^2\le1\),
\begin{equation}\label{eq:sec4-prelim-deletion}
R_n^2\le R_{-J}^2+k.
\end{equation}

We shall prove, uniformly for \(0\le k\le x^2/2\),
\begin{equation}\label{eq:sec4-prelim-deleted-tail}
\mathbb P(R_{-J}^2\ge x^2-k\mid\mathcal N_J)
\le C_d\mathbb P(\chi_d^2\ge x^2-k).
\end{equation}
The probability of \(\mathcal N_J\) is at least \(3/4\) by
\eqref{eq:sec4-prelim-exceedance}. The covariance matrix of the
remaining observations is
\[
C_J=I_d-\sum_{i\in J}
\mathbb E(\boldsymbol\xi_i\boldsymbol\xi_i^\top).
\]
Indeed, \eqref{eq:sec4-prelim-budget} implies
\[
\|I_d-C_J\|_{\mathrm{op}}
\le k\left(\sum_{i=1}^n\mathbb E|\boldsymbol\xi_i|^p\right)^{2/p}
\le\frac12\lambda^{2/p}.
\]
Thus \(I_d/2\preceq C_J\preceq I_d\). Define
\[
\widetilde{\boldsymbol\xi}_i=C_J^{-1/2}\boldsymbol\xi_i
\quad(i\notin J),\qquad
W_J=\sum_{i\notin J}
\widetilde{\boldsymbol\xi}_i\widetilde{\boldsymbol\xi}_i^\top.
\]
Before conditioning, these vectors are independent and centered, their
covariance matrices sum to \(I_d\), and
\begin{equation}\label{eq:sec4-prelim-normalized-moments}
\sum_{i\notin J}\mathbb E|\widetilde{\boldsymbol\xi}_i|^p
\le2^{p/2}\sum_{i=1}^n\mathbb E|\boldsymbol\xi_i|^p.
\end{equation}
By \eqref{eq:sec4-prob-variational}, \(R_{-J}^2\) is unchanged by
this normalization.

We first control \(W_J\) under the conditional law. Independence
and the conditional expectation formula give
\[
\mathbb E(W_J\mid\mathcal N_J)
=\sum_{i\notin J}
\frac{\mathbb E[\widetilde{\boldsymbol\xi}_i
\widetilde{\boldsymbol\xi}_i^\top;
|\boldsymbol\xi_i|\le b]}
{\mathbb P\{|\boldsymbol\xi_i|\le b\}}.
\]
Subtracting the unconditional second moments, and using
\(\|C_J^{-1/2}\|_{\mathrm{op}}\le\sqrt2\) and
\eqref{eq:sec4-large-moment-threshold}, yields
\begin{equation}\label{eq:sec4-prelim-matrix-mean}
\begin{aligned}
\|\mathbb E(W_J\mid\mathcal N_J)-I_d\|_{\mathrm{op}}
&\le C_db^{2-p}\sum_{i=1}^n\mathbb E|\boldsymbol\xi_i|^p
 +C_db^2\sum_{i=1}^n\mathbb P\{|\boldsymbol\xi_i|>b\}\\
&\le C_d\varepsilon^{2-p}\lambda/x^2\le\frac18.
\end{aligned}
\end{equation}
Conditional on \(\mathcal N_J\), the remaining
observations are still independent and satisfy
\(|\widetilde{\boldsymbol\xi}_i|\le\sqrt2b\). Their outer products
obey
\[
(\widetilde{\boldsymbol\xi}_i
\widetilde{\boldsymbol\xi}_i^\top)^2
=|\widetilde{\boldsymbol\xi}_i|^2
\widetilde{\boldsymbol\xi}_i\widetilde{\boldsymbol\xi}_i^\top
\preceq2b^2\widetilde{\boldsymbol\xi}_i
\widetilde{\boldsymbol\xi}_i^\top.
\]
Consequently, each centered matrix summand has operator norm at most
\(4b^2\), and the operator norm of the sum of its conditional
second moments is at most
\(2b^2\|\mathbb E(W_J\mid\mathcal N_J)\|_{\mathrm{op}}
\le C_db^2\).
Matrix Bernstein's inequality therefore gives
\begin{equation}\label{eq:sec4-prelim-matrix}
\mathbb P\{\|W_J-I_d\|_{\mathrm{op}}>1/4\mid\mathcal N_J\}
\le C_d\exp(-c_d/b^2)\le C_d e^{-4x^2}.
\end{equation}
The fixed choice of \(\varepsilon\) gives
\eqref{eq:sec4-prelim-matrix}.

It remains to control the radial event on
\(\{\|W_J-I_d\|_{\mathrm{op}}\le1/4\}\).
Put \(t=\sqrt{x^2-k}\), so that \(\sqrt8\le t\le x\).
For fixed \(u\in\mathbb S^{d-1}\), put
\(X_i=u^\top\widetilde{\boldsymbol\xi}_i\) and
\(a=\sqrt{t^2-1}\). Before conditioning,
\(\mathbb EX_i=0\) and \(\sum_{i\notin J}\mathbb EX_i^2=1\).
Since \(2<p\le3\), \eqref{eq:sec4-prelim-normalized-moments} gives
\[
\begin{aligned}
&a^2\max_{i\notin J}\mathbb EX_i^2
\le\left\{(1+a)^p\sum_{i\notin J}\mathbb E|X_i|^p\right\}^{2/p}
\le2\lambda^{2/p}\le2\lambda^{2/3}\le1,\\
&(1+a)^2\sum_{i\notin J}\mathbb E[X_i^2;(1+a)|X_i|>1]\\
&\quad+(1+a)^3\sum_{i\notin J}\mathbb E[|X_i|^3;(1+a)|X_i|\le1]\\
&\quad\le(1+a)^p\sum_{i\notin J}\mathbb E|X_i|^p
\le2^{p/2}\lambda\le2^{3/2}\lambda.
\end{aligned}
\]
\begingroup
The preceding inequalities verify the variance and
truncated-moment conditions of the scalar self-normalized
result of \citet{JingShaoWang2003}, after decreasing
\(\lambda_*\). Since that result has absolute constants,
this choice can be made independently of \(p\).
Applying the result under the original law \(\mathbb P\) gives
\endgroup
\begin{equation}\label{eq:sec4-prelim-directional}
\mathbb P\left\{
\frac{u^\top\sum_{i\notin J}\widetilde{\boldsymbol\xi}_i}
{\left\{\sum_{i\notin J}
(u^\top\widetilde{\boldsymbol\xi}_i)^2\right\}^{1/2}}
\ge\sqrt{t^2-1}\right\}
\le C_d\overline\Phi(\sqrt{t^2-1}).
\end{equation}
Here \(\overline\Phi\) is the standard normal upper tail, and a ratio
with zero denominator is set to zero. Since
\(\mathbb P(\mathcal N_J)\ge3/4\), the conditional version of
\eqref{eq:sec4-prelim-directional} holds with another constant.

We now pass from fixed directions to the radial event.
Let \(\sigma\) denote normalized surface measure on
\(\mathbb S^{d-1}\). On the good matrix event, the transformation
\[
u=\frac{W_J^{-1/2}w}{|W_J^{-1/2}w|}
\]
has the surface-measure relation
\[
d\sigma(w)=\frac{\det(W_J)^{1/2}}
{(u^\top W_Ju)^{d/2}}\,d\sigma(u).
\]
The eigenvalues of \(W_J\) lie in \([3/4,5/4]\), so
\[
(3/5)^{d/2}\le
\frac{\det(W_J)^{1/2}}{(u^\top W_Ju)^{d/2}}
\le(5/3)^{d/2}.
\]
This transformation changes the self-normalized projection
into \(w^\top W_J^{-1/2}\sum_{i\notin J}
\widetilde{\boldsymbol\xi}_i\).
When \(R_{-J}\ge t\), rotate coordinates so that
\(W_J^{-1/2}\sum_{i\notin J}\widetilde{\boldsymbol\xi}_i
=R_{-J}e_1\). The directions
\(w_1\ge\sqrt{1-t^{-2}}\) then have projection at least
\(\sqrt{t^2-1}\), and
\[
\sigma\{w:w_1\ge\sqrt{1-t^{-2}}\}
=c_d\int_0^{\arcsin(1/t)}\sin^{d-2}\theta\,d\theta
\ge c_dt^{-(d-1)}.
\]
Thus, pointwise,
\begin{equation}\label{eq:sec4-prelim-angular}
\begin{aligned}
&\mathbf1_{\{R_{-J}^2\ge t^2,\,
\|W_J-I_d\|_{\mathrm{op}}\le1/4\}}\\
&\quad\le C_dt^{d-1}\int_{\mathbb S^{d-1}}
\mathbf1\left\{
\frac{u^\top\sum_{i\notin J}\widetilde{\boldsymbol\xi}_i}
{\left\{\sum_{i\notin J}
(u^\top\widetilde{\boldsymbol\xi}_i)^2\right\}^{1/2}}
\ge\sqrt{t^2-1}\right\}\,d\sigma(u).
\end{aligned}
\end{equation}
Taking conditional expectations in \eqref{eq:sec4-prelim-angular},
applying Tonelli's theorem, and using
\eqref{eq:sec4-prelim-directional} gives an upper bound
\(C_dt^{d-1}\overline\Phi(\sqrt{t^2-1})\).
For \(t\ge\sqrt8\), Mills' inequality and the chi-squared tail bound
give
\[
t^{d-1}\overline\Phi(\sqrt{t^2-1})
\le C_dt^{d-2}e^{-t^2/2}
\le C_d\mathbb P(\chi_d^2\ge t^2).
\]
Also, \(t\le x\) and \(d\ge2\) give
\(e^{-4x^2}\le e^{-t^2/2}\le\mathbb P(\chi_d^2\ge t^2)\).
Adding \eqref{eq:sec4-prelim-matrix} proves
\eqref{eq:sec4-prelim-deleted-tail}.

Let \(K_b=\sum_{i=1}^n\mathbf1_{\{|\boldsymbol\xi_i|>b\}}\).
By \eqref{eq:sec4-prelim-deletion}, independence, and
\eqref{eq:sec4-prelim-deleted-tail}, for \(0\le k=|J|\le x^2/2\),
\begin{equation}\label{eq:sec4-prelim-large-set}
\begin{aligned}
&\mathbb P\bigl(R_n^2\ge x^2,\mathcal N_J,
 |\boldsymbol\xi_i|>b\text{ for all }i\in J\bigr)\\
&\quad\le\prod_{i\in J}\mathbb P\{|\boldsymbol\xi_i|>b\}
\mathbb P(R_{-J}^2\ge x^2-k\mid\mathcal N_J)\\
&\quad\le C_de^{k/2}\mathbb P(\chi_d^2\ge x^2)
\prod_{i\in J}\mathbb P\{|\boldsymbol\xi_i|>b\}.
\end{aligned}
\end{equation}
The events on the left of \eqref{eq:sec4-prelim-large-set} are
disjoint for distinct \(J\).
Here we used
\(\mathbb P(\chi_d^2\ge x^2-k)
\le e^{k/2}\mathbb P(\chi_d^2\ge x^2)\), which follows from the
monotonicity of \(v\mapsto e^{v/2}\mathbb P(\chi_d^2\ge v)\)
for \(d\ge2\).
Taking \(J=\varnothing\) in \eqref{eq:sec4-prelim-large-set} gives
\begin{equation}\label{eq:sec4-prelim-zero-large}
\mathbb P(R_n^2\ge x^2,K_b=0)
\le C_d\mathbb P(\chi_d^2\ge x^2).
\end{equation}
Summing \eqref{eq:sec4-prelim-large-set} over
\(1\le|J|\le x^2/2\) and extending the sum on the right to all
nonempty \(J\) gives
\begin{equation}\label{eq:sec4-prelim-few-large}
\begin{aligned}
&\mathbb P(R_n^2\ge x^2,1\le K_b\le x^2/2)\\
&\quad\le C_d\mathbb P(\chi_d^2\ge x^2)
\sum_{k\ge1}\frac1{k!}
\left(e^{1/2}\sum_{i=1}^n
\mathbb P\{|\boldsymbol\xi_i|>b\}\right)^k\\
 &\quad\le C_d\lambda\mathbb P(\chi_d^2\ge x^2).
\end{aligned}
\end{equation}
The exponential series and \eqref{eq:sec4-prelim-exceedance}
give \eqref{eq:sec4-prelim-few-large}.

For the remaining event, put \(r=\lfloor x^2/2\rfloor+1\).
The factorial-moment bound for the sum of independent indicators gives
\begin{equation}\label{eq:sec4-prelim-count-tail}
\mathbb P(K_b>x^2/2)
=\mathbb P(K_b\ge r)
\le\frac1{r!}\left(\sum_{i=1}^n
\mathbb P\{|\boldsymbol\xi_i|>b\}\right)^r.
\end{equation}
Since \(\varepsilon\le1\),
\eqref{eq:sec4-large-moment-threshold} gives
\(\varepsilon^{-p}\lambda_*\le1/4<e^{-1}\).
Since \(r>x^2/2\), \eqref{eq:sec4-prelim-count-tail} and
\eqref{eq:sec4-prelim-exceedance} give
\begin{equation}\label{eq:sec4-prelim-many-large}
\mathbb P(K_b>x^2/2)
\le\varepsilon^{-p}\lambda e^{-(r-1)}
\le e\varepsilon^{-p}\lambda e^{-x^2/2}
\le C_d\lambda\mathbb P(\chi_d^2\ge x^2).
\end{equation}
Combining \eqref{eq:sec4-prelim-zero-large},
\eqref{eq:sec4-prelim-few-large}, and
\eqref{eq:sec4-prelim-many-large} proves
\eqref{eq:sec4-prelim-coarse}. Combining
\eqref{eq:sec4-prelim-few-large} and
\eqref{eq:sec4-prelim-many-large} gives
\begin{equation}\label{eq:sec4-prelim-large}
\mathbb P(R_n^2\ge x^2,\mathcal N^c)
\le C_d\lambda\mathbb P(\chi_d^2\ge x^2),
\end{equation}
since \(\mathcal N^c=\{K_b\ge1\}\).

Substituting \eqref{eq:sec4-prelim-large},
\eqref{eq:sec4-prelim-coarse}, and
\eqref{eq:sec4-prelim-exceedance} into
\eqref{eq:sec4-conditioning-transfer-formula}, with
\(A=\{R_n^2\ge x^2\}\), proves \eqref{eq:sec4-prelim-transfer}.
\end{proof}


\begin{proof}[Proof of Lemma~\ref{lem:sec4-truncated}]
Use Lemma~\ref{lem:sec4-conditional-moments} with
\(b=\varepsilon/x\), giving \eqref{eq:sec4-prob-moments}.
Write \(M=\mathbb E^*V_n\). We have
\[
\|M-I_d\|_{\mathrm F}\le C_d\lambda/x^2\le1/16,\qquad
\|M\|_{\mathrm F}\le\sqrt d+1/16,
\]
by \eqref{eq:sec4-prob-moments} and
\eqref{eq:sec4-large-moment-threshold}. Hence
\[
\mathcal G^c\subset
\{\|V_n-M\|_{\mathrm F}>3/16\}.
\]
Tail integration gives
\begin{equation}\label{eq:sec4-truncated-weighted-tail}
\begin{aligned}
\mathbb E^*[(1+\|V_n\|_{\mathrm F})\mathbf1_{\mathcal G^c}]
&\le(1+\|M\|_{\mathrm F}+3/16)
 \mathbb P^*(\|V_n-M\|_{\mathrm F}>3/16)\\
&\quad+\int_{3/16}^\infty
 \mathbb P^*(\|V_n-M\|_{\mathrm F}>t)\,dt.
\end{aligned}
\end{equation}

For fixed coordinates $r,s$, the centered entry $(V_n-M)_{rs}$
is a sum of independent centered variables under $\mathbb P^*$.
Each summand has absolute value at most $2b^2$, and their total
variance is at most
$\sum_i\mathbb E^*|\boldsymbol\xi_i|^4\le C_d\lambda/x^4$.
Put $T=64d x^2$. Then $2Tb^2=128d\varepsilon^2\le1$
by the fixed choice of $\varepsilon$.
Applying $e^u\le1+u+u^2$ for $|u|\le1$ to each centered summand,
and then using independence, gives
\begin{equation}\label{eq:sec4-truncated-entry-mgf}
\log\mathbb E^*\exp\{\pm T(V_n-M)_{rs}\}
\le T^2\sum_{i=1}^n\mathbb E^*|\boldsymbol\xi_i|^4
\le C_d\lambda.
\end{equation}
By \eqref{eq:sec4-truncated-entry-mgf},
\begin{equation}\label{eq:sec4-truncated-entry-cosh}
\mathbb E^*\bigl[\cosh\{T(V_n-M)_{rs}\}-1\bigr]
\le e^{C_d\lambda}-1\le C_d\lambda.
\end{equation}
If $\|V_n-M\|_{\mathrm F}>t$, some entry has absolute value greater
than $t/d$. Applying Markov's inequality to
\eqref{eq:sec4-truncated-entry-cosh} and taking a union bound over
the $d^2$ entries gives
\begin{equation}\label{eq:sec4-truncated-matrix-tail}
\mathbb P^*(\|V_n-M\|_{\mathrm F}>t)
\le C_d\lambda e^{-64x^2t},\qquad t\ge1/32.
\end{equation}
Indeed, $Tt/d=64x^2t\ge2x^2$, so
$\cosh(Tt/d)-1$ is bounded below by a constant times $e^{64x^2t}$.

Substituting \eqref{eq:sec4-truncated-matrix-tail} into
\eqref{eq:sec4-truncated-weighted-tail} gives
\[
\mathbb E^*[(1+\|V_n\|_{\mathrm F})\mathbf1_{\mathcal G^c}]
\le C_d\lambda e^{-12x^2}.
\]
Since $\mathbb P(\chi_d^2\ge x^2)\ge e^{-x^2/2}$ for $d\ge2$,
\eqref{eq:sec4-full-matrix} follows.

\end{proof}


\begin{proof}[Proof of Lemma~\ref{lem:sec4-conditional}]
For \eqref{eq:sec4-prob-bad-matrices}, the moment bounds
\eqref{eq:sec4-prob-moments} give
\begin{equation}\label{eq:sec4-prob-deleted-mean}
\begin{aligned}
\|\mathbb E^*V_{-i}-I_d\|_{\mathrm F}
&\le\|\mathbb E^*V_n-I_d\|_{\mathrm F}
 +\|\mathbb E^*(\boldsymbol\xi_i\boldsymbol\xi_i^\top)\|_{\mathrm F}\\
&\le C_d\lambda/x^2+b^2\le3/32,
\end{aligned}
\end{equation}
by \eqref{eq:sec4-large-moment-threshold}. Hence
\[
\mathcal G_i^c\subset
\{\|V_{-i}-\mathbb E^*V_{-i}\|_{\mathrm F}>1/32\},\qquad
\|\mathbb E^*V_{-i}\|_{\mathrm F}\le\sqrt d+3/32.
\]
Removing the \(i\)th summand preserves the entry bound \(2b^2\)
and the fourth-moment bound \(C_d\lambda/x^4\) used in
\eqref{eq:sec4-truncated-entry-mgf}.
The proof of \eqref{eq:sec4-truncated-matrix-tail} therefore gives
\begin{equation}\label{eq:sec4-prob-matrix-tail}
\sup_{1\le i\le n}\mathbb P^*
\{\|V_{-i}-\mathbb E^*V_{-i}\|_{\mathrm F}>t\}
\le C_d\lambda e^{-64x^2t},\qquad t\ge1/32.
\end{equation}
Apply the tail integration in \eqref{eq:sec4-truncated-weighted-tail}
at \(1/32\). Equations \eqref{eq:sec4-prob-deleted-mean} and
\eqref{eq:sec4-prob-matrix-tail} yield
\[
\begin{aligned}
\mathbb E^*[(1+\|V_{-i}\|_{\mathrm F})\mathbf1_{\mathcal G_i^c}]
&\le C_d\lambda\left\{e^{-2x^2}
 +\int_{1/32}^\infty e^{-64x^2t}\,dt\right\}\\
&\le C_d\lambda e^{-2x^2},
\end{aligned}
\]
which proves \eqref{eq:sec4-prob-bad-matrices} uniformly in \(i\).

For \eqref{eq:sec4-prob-envelope-integrals}, we first bound the
tails after deleting one or two observations.
Since $2u-u^2\le1$, \eqref{eq:sec4-prob-variational} yields
\begin{equation}\label{eq:sec4-prob-deletion}
R_n^2\le R_{-i}^2+1,\qquad
R_{-i}^2\le R_{-\{i,j\}}^2+1\quad(i\ne j),
\end{equation}
where the latter statistic is obtained by deleting both observations.

For
$J\subset\{1,\ldots,n\}$ with $1\le|J|\le2$, define
\[
C_J=I_d-\sum_{j\in J}
          \mathbb E(\boldsymbol\xi_j\boldsymbol\xi_j^\top).
\]
The moment assumption gives
\[
\|C_J-I_d\|_{\mathrm{op}}
\le\sum_{j\in J}\mathbb E|\boldsymbol\xi_j|^2
\le2\left(\sum_{j=1}^n\mathbb E|\boldsymbol\xi_j|^p\right)^{2/p}
\le2\lambda^{2/p}.
\]
\begingroup
Since \(\lambda^{2/p}\le\lambda^{2/3}\),
\eqref{eq:sec4-large-moment-threshold} makes every such \(C_J\)
positive definite and gives
\[
 (21/20)^p\|C_J^{-1/2}\|_{\mathrm{op}}^p
 \le(21/20)^3(1-2\lambda_*^{2/3})^{-3/2}\le2.
\]
\endgroup
The vectors $C_J^{-1/2}\boldsymbol\xi_\ell$, $\ell\notin J$, are
centered and independent, with total covariance $I_d$. For
$0\le t\le x+1/x$, we have $(1+t)/(1+x)\le21/20$, and hence
\begin{equation}\label{eq:sec4-prob-deleted-budget}
(1+t)^p\sum_{\ell\notin J}
 \mathbb E|C_J^{-1/2}\boldsymbol\xi_\ell|^p
\le2\lambda.
\end{equation}
By \eqref{eq:sec4-prob-variational}, their self-normalized statistic is
$R_{-J}^2$. For $t\ge4$, apply \eqref{eq:sec4-prelim-coarse}
with parameter $2\lambda\le\lambda_*$. For $0\le t\le4$, use
$\mathbb P(R_{-J}^2\ge t^2)\le1$ and
$\mathbb P(\chi_d^2\ge t^2)\ge\mathbb P(\chi_d^2\ge16)>0$.
Hence, by \eqref{eq:sec4-prelim-coarse},
\begin{equation}\label{eq:sec4-prob-deleted-coarse}
\mathbb P(R_{-J}^2\ge t^2)
\le C_d\mathbb P(\chi_d^2\ge t^2),\qquad
0\le t\le x+1/x,\quad 1\le|J|\le2.
\end{equation}

We also use the elementary chi-squared tail comparison
\begin{equation}\label{eq:sec4-prob-chi-shift}
\mathbb P(\chi_d^2\ge u-v)
\le e^{v/2}\mathbb P(\chi_d^2\ge u),\qquad 0\le v\le u.
\end{equation}
Indeed, for $d\ge2$ the function $e^{u/2}\mathbb P(\chi_d^2\ge u)$
is nondecreasing, as is seen from the representation
\[
e^{u/2}\mathbb P(\chi_d^2\ge u)
=c_d\int_0^\infty(u+s)^{d/2-1}e^{-s/2}\,ds.
\]

The coarse bound \eqref{eq:sec4-prob-deleted-coarse} and
\eqref{eq:sec4-conditioning-bound} imply
\begin{equation}\label{eq:sec4-prob-deleted-upper}
\mathbb P^*(R_{-i}\ge t)
\le\frac{\mathbb P(R_{-i}\ge t)}{\mathbb P(\mathcal N)}
\le C_d\mathbb P(\chi_d^2\ge t^2),\qquad 0\le t\le x+1/x.
\end{equation}

We use \eqref{eq:sec4-prob-deleted-upper} to estimate the envelopes
in \eqref{eq:sec4-stein-envelopes}. On $y>x-x^{-1}$,
\eqref{eq:sec4-stein-envelopes},
\eqref{eq:sec4-prob-deleted-upper}, and
\eqref{eq:sec4-prob-chi-shift} give
\begin{equation}\label{eq:sec4-prob-envelope-above}
\begin{aligned}
\mathbb E^*[A_{k,x}(R_{-i})
              \mathbf1_{\{R_{-i}>x-x^{-1}\}}]
&\le(1+x)^{k-2}\mathbb P^*(R_{-i}\ge x-x^{-1})\\
&\le C_d(1+x)^{k-2}\mathbb P(\chi_d^2\ge x^2).
\end{aligned}
\end{equation}
For $0\le y<x-x^{-1}$, \eqref{eq:sec4-stein-envelopes} gives
$A_{k,x}(0)=\mathbb P(\chi_d^2\ge x^2)$ and
\[
|A_{k,x}'(y)|
\le C_d\mathbb P(\chi_d^2\ge x^2)
      e^{y^2/2}(1+y)^{k-d+1},\qquad 0\le y<x-x^{-1}.
\]
The fundamental theorem of calculus applied to the first branch in
\eqref{eq:sec4-stein-envelopes} gives, for every $u\ge0$,
\[
A_{k,x}(u)\mathbf1_{\{u\le x-x^{-1}\}}
\le A_{k,x}(0)+\int_0^{x-x^{-1}}
       |A_{k,x}'(y)|\mathbf1_{\{u\ge y\}}\,dy.
\]
Taking expectations at $u=R_{-i}$ and applying Tonelli's theorem,
\eqref{eq:sec4-prob-deleted-upper}, and the elementary bound
$\mathbb P(\chi_d^2\ge y^2)\le C_d(1+y)^{d-2}e^{-y^2/2}$ gives
\begin{equation}\label{eq:sec4-prob-envelope-below}
\begin{aligned}
\mathbb E^*[A_{k,x}(R_{-i})
               \mathbf1_{\{R_{-i}\le x-x^{-1}\}}]
&\le A_{k,x}(0)+\int_0^{x-x^{-1}}
      |A_{k,x}'(y)|\mathbb P^*(R_{-i}\ge y)\,dy\\
&\le C_d\mathbb P(\chi_d^2\ge x^2)
       \left\{1+\int_0^{x-x^{-1}}(1+y)^{k-1}\,dy\right\}\\
&\le C_d\mathbb P(\chi_d^2\ge x^2)(1+x)^k.
\end{aligned}
\end{equation}
Combining \eqref{eq:sec4-prob-envelope-above} and
\eqref{eq:sec4-prob-envelope-below} proves
\eqref{eq:sec4-prob-envelope-integrals}, since $A_{k,x}\ge0$ and
$\mathbf1_{\mathcal G_i}\le1$.
\end{proof}

\begin{proof}[Proof of Corollary~\ref{lem:sec4-deleted-comparison}]
By the triangle inequality,
\[
\begin{aligned}
&\left|\mathbb P^*(R_{-i}\ge t)-\mathbb P(\chi_d^2\ge t^2)\right|\\
&\quad\le\left|\mathbb P(R_{-i}\ge t)-\mathbb P(\chi_d^2\ge t^2)\right|\\
&\qquad+\left|\mathbb P^*(R_{-i}\ge t)-\mathbb P(R_{-i}\ge t)\right|.
\end{aligned}
\]
Equations \eqref{eq:sec4-prob-deleted-budget},
\eqref{eq:sec4-prob-variational}, and \eqref{eq:sec4-error} give
\begin{equation}\label{eq:sec4-prob-deleted-original}
\left|\mathbb P(R_{-i}\ge t)-\mathbb P(\chi_d^2\ge t^2)\right|
\le\mathcal E(2\lambda)\mathbb P(\chi_d^2\ge t^2),
\qquad 0\le t\le x+1/x.
\end{equation}
Use \(\mathcal N_{\{i\}}\) from
\eqref{eq:sec4-prelim-deleted-sample}.
Because $R_{-i}$ uses only the deleted sample, its distribution under
$\mathbb P^*$ is its original distribution conditional on $\mathcal N_{\{i\}}$.
By \eqref{eq:sec4-conditioning-bound},
$\mathbb P(\mathcal N_{\{i\}})\ge\mathbb P(\mathcal N)\ge3/4$.
For $3\le t\le x+1/x$, delete a second observation on
$\mathcal N_{\{i\}}^c$. By
\eqref{eq:sec4-prob-deletion}, independence,
\eqref{eq:sec4-prob-deleted-coarse}, and
\eqref{eq:sec4-prob-chi-shift},
\begin{equation}\label{eq:sec4-prob-deleted-large-event}
\begin{aligned}
\mathbb P(R_{-i}\ge t,\mathcal N_{\{i\}}^c)
&\le\sum_{j\ne i}\mathbb P(|\boldsymbol\xi_j|>b)
       \mathbb P(R_{-\{i,j\}}^2\ge t^2-1)\\
&\le C_d\lambda\,\mathbb P(\chi_d^2\ge t^2).
\end{aligned}
\end{equation}
The bound \eqref{eq:sec4-prob-deleted-coarse} applies at
$\sqrt{t^2-1}$ because $0\le\sqrt{t^2-1}\le x+1/x$.
Apply \eqref{eq:sec4-conditioning-transfer-formula} with
\(A=\{R_{-i}\ge t\}\) and \(\mathcal N\) replaced by
\(\mathcal N_{\{i\}}\). Equations
\eqref{eq:sec4-prob-deleted-large-event},
\eqref{eq:sec4-prelim-exceedance}, and
\eqref{eq:sec4-prob-deleted-coarse} give
\begin{equation}\label{eq:sec4-prob-deleted-transfer}
\left|\mathbb P^*(R_{-i}\ge t)-\mathbb P(R_{-i}\ge t)\right|
\le C_d\lambda\,\mathbb P(\chi_d^2\ge t^2),
\qquad 3\le t\le x+1/x.
\end{equation}
Combining \eqref{eq:sec4-prob-deleted-transfer} with
\eqref{eq:sec4-prob-deleted-original} proves
\eqref{eq:sec4-prob-deleted-comparison}.
\end{proof}

\begin{proof}[Proof of Corollary~\ref{cor:phewas-bonferroni}]
Independence of the two subjects in each pair gives
\[
 \mathbb E Z_{1j}
 =\mathbb E(YX_j)-(\mathbb E Y)(\mathbb E X_j)
 =\operatorname{Cov}(Y,X_j).
\]
If \(\mathcal H_0\) is empty, \eqref{eq:apps-fwer} is immediate.
Otherwise, fix a true null index \(j\). Since the pairs are fixed and disjoint,
\(Z_{1j},\ldots,Z_{Nj}\) are independent and identically distributed
with mean zero. Set
\[
 \xi_{kj}=N^{-1/2}\Gamma_j^{-1/2}Z_{kj},\qquad 1\le k\le N.
\]
Then
\[
 \sum_{k=1}^N\mathbb E(\xi_{kj}\xi_{kj}^{\top})=I_d,
 \qquad
 \sum_{k=1}^N\mathbb E|\xi_{kj}|^p
 =N^{-(p-2)/2}\mathbb E|\Gamma_j^{-1/2}Z_{1j}|^p
 \le K N^{-(p-2)/2}.
\]
The invertible change of dual variable
\(\vartheta=N^{-1/2}\Gamma_j^{-1/2}\theta\) shows that
\[
 \sup_{\theta\in\mathbb R^d}
 \left\{2\theta^{\top}\sum_{k=1}^N\xi_{kj}
 -\theta^{\top}\Bigl(\sum_{k=1}^N\xi_{kj}\xi_{kj}^{\top}\Bigr)\theta\right\}
 =U_j.
\]
Theorem~\ref{thm:main-CMD-T} therefore gives, uniformly over true null indices,
\[
 \left|\frac{\mathbb P(U_j\ge x^2)}{Q_d(x^2)}-1\right|
 \le C_d K N^{-(p-2)/2}(1+x)^p
 \quad\text{if}\quad
 K N^{-(p-2)/2}(1+x)^p\le c_d.
\]

Let \(q_{\alpha,M}\) be defined by
\(Q_d(q_{\alpha,M})=\alpha/M\), where \(0<\alpha<1\).
The exponential Markov inequality and
\(\mathbb E\exp(\chi_d^2/4)=2^{d/2}\) imply
\[
 Q_d(t)\le 2^{d/2}e^{-t/4},\qquad
 q_{\alpha,M}\le 4\log(M/\alpha)+2d\log 2.
\]
Consequently,
\[
 (1+\sqrt{q_{\alpha,M}})^p
 \le C_d\{1+\log(M/\alpha)\}^{p/2},
\]
with a constant uniform for \(2<p\le3\). After reducing the constant
\(c_d\) in the corollary, its stated condition permits the preceding
moderate deviation result to be applied at \(x=\sqrt{q_{\alpha,M}}\).
Writing \(\mathcal H_0\) for the true null indices, the union bound gives
\begin{align*}
 \mathbb P\bigl(U_j\ge q_{\alpha,M}\text{ for some }j\in\mathcal H_0\bigr)
 &\le \sum_{j\in\mathcal H_0}\mathbb P(U_j\ge q_{\alpha,M})\\
 &\le \frac{|\mathcal H_0|}{M}\alpha
 \left[1+C_dK N^{-(p-2)/2}
       \{1+\log(M/\alpha)\}^{p/2}\right]\\
 &\le \alpha
 \left[1+C_dK N^{-(p-2)/2}
       \{1+\log(M/\alpha)\}^{p/2}\right].
\end{align*}
No independence across the indices \(j\) is needed.

To verify the sufficient conditions on the original observations,
expand the paired product once more. Without imposing the null hypothesis,
\begin{align*}
 \Gamma_j
 &=\frac12\mathbb E\bigl\{(Y-\mathbb E Y)(Y-\mathbb E Y)^{\top}
                    (X_j-\mathbb E X_j)^2\bigr\}\\
 &\quad+\frac12\operatorname{Var}(X_j)\operatorname{Cov}(Y)\\
 &\succeq\frac{\operatorname{Var}(X_j)}{2}\operatorname{Cov}(Y).
\end{align*}
Suppose that \(|Y-\mathbb E Y|\le B\) almost surely,
\(\operatorname{Cov}(Y)\succeq\kappa I_d\) for some \(\kappa>0\),
and \(\operatorname{Var}(X_j)>0\) for every true null index.
Then \(\Gamma_j\) is positive definite, and
\[
 \mathbb E|Z_{1j}|^p
 \le B^p\mathbb E|X_{2j}-X_{1j}|^p
 \le 2^pB^p\mathbb E|X_j-\mathbb E X_j|^p.
\]
It follows that
\[
 \mathbb E|\Gamma_j^{-1/2}Z_{1j}|^p
 \le 2^{3p/2}B^p\kappa^{-p/2}
       \frac{\mathbb E|X_j-\mathbb E X_j|^p}
            {\{\operatorname{Var}(X_j)\}^{p/2}}.
\]
A uniform bound on the displayed standardized moment ratio therefore
verifies the required moment assumption. No positive definiteness of
\(\operatorname{Var}\{(Y-\mathbb E Y)(X_j-\mathbb E X_j)\}\) is required.

Finally, fix a coordinate under a fixed alternative distribution with
\(\mathbb E Z_{1j}\ne0\), \(\mathbb E|Z_{1j}|^2<\infty\), and
\(\mathbb E(Z_{1j}Z_{1j}^{\top})\) positive definite.
The strong law gives, almost surely,
\[
 \frac1N\sum_{k=1}^N Z_{kj}\longrightarrow\mathbb E Z_{1j},
 \qquad
 \frac1N\sum_{k=1}^N Z_{kj}Z_{kj}^{\top}
       \longrightarrow\mathbb E(Z_{1j}Z_{1j}^{\top}).
\]
The sample second-moment matrix is therefore eventually invertible, and
maximizing the quadratic expression defining \(U_j\) yields
\[
 \frac{U_j}{N}\longrightarrow
 (\mathbb E Z_{1j})^{\top}
 \{\mathbb E(Z_{1j}Z_{1j}^{\top})\}^{-1}
 \mathbb E Z_{1j}>0
 \qquad\text{almost surely}.
\]
For fixed \(\alpha\), the quantile bound gives \(q_{\alpha,M}=o(N)\)
whenever \(\log M=o(N)\), proving that the rejection probability tends to one.
\end{proof}

\section*{Concurrent work}
\citet{JingLiuWangZhou2026} posted their preprint while we were
preparing this manuscript. They use a geometric argument to establish
a Cram\'er-type moderate deviation result for the one-sample Hotelling
statistic based on independent and identically distributed
observations with a finite third moment.
We use Stein's method to obtain explicit finite-sample estimates,
combining a radial Stein identity with a recursive error estimate.
Our results apply to independent observations with possibly different
distributions under finite \(p\)th moments for \(2<p\le3\).
We also establish uniform Berry--Esseen bounds and give explicit
constructions showing sharpness of the moderate deviation
error order over independent arrays.

\begin{funding}
Shao Q.M. was partially supported by 2025JC09A021, National Nature Science Foundation of China NSFC 12031005 and Shenzhen Outstanding Talents Training Fund, China.
Liu S.H. was partially supported by the Fundamental Research Funds for the Central Universities DUT25RC(3)133.
\end{funding}

\section*{AI Statement}
During the preparation of this manuscript, the authors used ChatGPT and
Codex (OpenAI) for language editing, LaTeX formatting, and assistance
with proof exposition and organization. These tools were also used to
check references, notation, cross-references, and mathematical derivations.
The authors reviewed and revised the manuscript, verified the mathematical
arguments, and take full responsibility for its content.

\end{document}